\documentclass[a4paper,12pt]{article}
\usepackage{amsmath}
\usepackage{amssymb}
\usepackage{amsthm}
\usepackage{amsfonts}
\usepackage{mathrsfs}
\usepackage{graphicx}
\usepackage{hyperref}
\usepackage{float}
\usepackage{slashed}
\usepackage{cite}

\newtheorem{theorem}{Theorem}

\newtheorem{definition}{Definition}

\newtheorem{preposition}{Preposition}
\newtheorem{assumption}{Assumption}

\numberwithin{equation}{section}

\newtheoremstyle{named}{}{}{\itshape}{}{\bfseries}{.}{.5em}{\thmnote{#3 }#1}
\theoremstyle{named}

\hypersetup{ pdftitle={Report 2},
	pdfauthor={Mulyanto}, 
	pdfkeywords={Maxwell Higgs, Schwarzchilc}, 
	bookmarksnumbered, pdfstartview={FitH}, urlcolor=blue, }

\begin{document}
	\pagestyle{plain}
	
	
	
	
	\title{\LARGE\textbf{Decay Estimate of Maxwell-Higgs System on Schwarzschild Black Holes}}
	
	\author{{Mulyanto, Fiki Taufik Akbar, and Bobby Eka Gunara\footnote{Corresponding author} }\\ 
		\textit{\small Theoretical Physics Laboratory,}\\
		\textit{\small Theoretical High Energy Physics Research Division,}\\
		\textit{\small Faculty of Mathematics and Natural Sciences,}\\
		\textit{\small Institut Teknologi Bandung}\\
		\textit{\small Jl. Ganesha no. 10 Bandung, Indonesia, 40132}\\ \\
		\small email: 30220313@mahasiswa.itb.ac.id, ftakbar@itb.ac.id, bobby@itb.ac.id}

	\date{}
	
	\maketitle
	
	
	
	
\begin{abstract}
		
In this paper, we prove the decay  estimate of Maxwell-Higgs system on four dimensional Schwarzschild spacetimes. We show that if the field equations support a Morawetz type estimate supported around the trapped surface, the uniform decay properties in the entire exterior of the Schwarzschild black holes can be obtained by using Sobolev inequalities and energy estimates. Our results also consider various forms of physical potential such as the mass terms, $\phi^4$-theory, sine Gordon potential, and Toda potential.		
\end{abstract}

	
	
	
\section{Introduction}
\label{sec:Introduction}
\subsection{Motivation and the relation to the existing literature}
The behaviour of solutions to the Maxwell Klein-Gordon (MKG) or Maxwell-Higgs (MH) equation is of considerable interest in physics, as it provides an insight into the fundamental properties of electromagnetic fields and their interaction with charged particles. The decay and the global existence of solutions to the MKG equation have been studied extensively over the past few decades, and many significant results have been obtained. In four dimensional Minkowski spacetime, the global existence of the theory has been obtained by considering various gauge conditions \cite{Klainerman, Selberg, yuan}. In particular, there are several works obtaining the time decay solutions of MKG theory on four dimensional flat spacetime, see for example, \cite{Petrescu, klainerman2, fang}.

On curved spacetimes, for example, the author in \cite{Blue} obtained the time decay of Maxwell equations on Schwarzschild spacetimes in which one introduces the middle components $F_{vw}$ and $F_{\theta\varphi}$, where $F$ is the Maxwell field strength and $(v,w)$ are Eddington-Finkelstein coordinates
\begin{eqnarray}
  w=t-r^{*} ~ ,\label{EFw}\\
  v=t+r^{*}  ~ , \label{EFv}
\end{eqnarray}
where $r^{*}$ is the Regge-Wheeler tortoise coordinate  defined by
\begin{eqnarray}\label{RW}
    r^{*}=r+2m \text{log} \left(r-2m\right)  ~ .
\end{eqnarray}
 Several years later, some authors in \cite{Ander1} have modified the method in \cite{Blue} by introducing a new Morawetz local decay estimate for the Maxwell equation on the exterior of a Schwarzschild black hole such that the solutions are bounded above an  energy functional. In nonlinear cases, however, if we want to extend this to the case of Yang-Mills fields or MH equation, we are unable to separate the middle components from the other components. Then, using the conformal energy method, the author in \cite{Ghanem1} later proved the uniform decay without separating the middle components of the Maxwell fields on Schwarzschild spacetimes. In order to solve the problem, one must assume that the non-stationary solutions of the Maxwell fields verify the Morawetz type estimate at the zero-derivative level. 

On the other hand, in the complex scalar field case, there has been the work in \cite{dafermos} demonstrating the decay for free scalar wave equation solutions, namely, $\Box\phi=0$, in the exterior of the Schwarzschild black holes. While in \cite{Blue2} by obtaining the pointwise decay estimates for the inhomogeneous wave equation on a Schwarzschild background, the authors have showed a similar decay result for initial data that vanishes on the bifurcated sphere, with decay rate  weaker than \cite{dafermos} along the event horizon. There are also some works on other spacetimes such as Kerr geometry in \cite{Finster2} and Reissner-Nordstr\"{o}m–de Sitter black hole in \cite{costa}.  

So far, there is no study of the decay estimate of MH theory on curved spacetimes, in particular, Schwarzschild spacetimes. If one tries to modify the method of \cite{dafermos} for the free scalar wave equation to MH system on Schwarzschild spacetimes, by using the energy-momentum tensor, two problems arise. First, one needs to control the flux energy near the surface  at $r=3m$ that involves the middle components of the Maxwell ﬁeld. Second, one must also define a weighted Sobolev-type norms of energy and prove that the energy is finite. However, in MH theory it cannot be showed that the defined energy functional is finite since there are some terms containing the third derivative of the complex scalar field. It is also interesting to extend the analysis of \cite{Ghanem1} to determine the decay of Maxwell-Higgs theory. However, due to the coupling between the $U(1)$ gauge field and the complex scalar field expressed in the covariant derivative of the scalar field, we have to modify the method in \cite{Ghanem1} in order to get the estimate of the coupling term. Thus, we cannot directly employ the setup in \cite{Ghanem1} and \cite{dafermos} to our case.

\subsection{Summary of results and article organization}

 The aim of this paper is to consider MH theory on Schwarzschild spacetimes. Here, we prove the decay estimate for the MH system using the energy-momentum tensor directly together with the appropriate Sobolev inequalities with addition of general potential $P\left(\phi,\bar{\phi}\right)$. Our starting point is the Lagrangian of MH system which can be written down as	
\begin{equation}\label{lagrangian}
		\mathcal{L}=-\frac{1}{4}{{F}_{\mu \nu }}{{F}^{\mu \nu }}-{{D}_{\mu }}\phi \overline{{{D}^{\mu }}\phi }-P\left( \phi ,\overline{\phi } \right)~,
\end{equation}
where $F_{\mu \nu } \equiv {{\partial }_{\mu }}A_{\nu }-{{\partial }_{\nu }}A_{\mu }$ is the gauge field strength and $D_{\mu}{\phi}\equiv\partial _{\mu }{\phi }-iA_{\mu}{\phi }$, is the covariant derivative of the complex scalar fields. The Lagrangian (\ref{lagrangian}) describe the interaction of a two-form function $F_{\mu\nu}$ and a complex scalar field $\phi$ on 4-dimensional Schwarzschild spacetimes with the standard coordinate $x^{\mu} = (t,x^i)$ where $\mu=0,1,2,3$ and $i=1,2,3$, endowed with metric

\begin{equation}\label{metric}
d{{s}^{2}}=-\left( 1-\frac{2m}{r} \right)d{{t}^{2}}+\frac{1}{\left( 1-\frac{2m}{r} \right)}d{{r}^{2}}+{{r}^{2}}d{{\Omega }^{2}}~.
\end{equation}

The fields equations of motions and the energy-momentum tensor for the MH system with respect to the Lagrangian (\ref{lagrangian}) respectively are 
\begin{eqnarray}\label{eom1}
	{{\partial }^{\alpha }}{{F}_{\alpha \gamma }}=-i\left( {{D}_{\gamma }}\phi \overline{\phi }-\phi \overline{{{D}_{\gamma }}\phi } \right)~,
\end{eqnarray}
\begin{eqnarray}\label{eom2}
	{{D}^{\mu }}{{D}_{\mu }}\phi ={{\partial }_{{\bar{\phi }}}}P~,
\end{eqnarray}
\begin{equation}
	{{{T}}^{\mu \nu }}=F_{\gamma }^{\nu }{{F}^{\mu \gamma }}-\frac{1}{4}{{g}^{\mu \nu }}{{F}^{\alpha \beta }}{{F}_{\alpha \beta }}+{{D}^{\mu }}\phi \overline{{{D}^{\nu }}\phi }-\frac{1}{2}~{{g}^{\mu \nu }}~{{D}_{\gamma }}\phi ~\overline{{{D}^{\gamma }}\phi }-{{g}^{\mu \nu }}P~.
\end{equation}
We split our proof of the decay of MH system on the Schwarzschild spacetimes into two parts. The first part proves the decay for far away from the horizon ($r>3m$), and the second part considers the region near the horizon ($2m<r<3m$).

In the first part, by modifying the idea in \cite{Ghanem1}, we get the system energy $E^{\hat{t}}$ generated by the time-like vector field in \eqref{energytimelike}. Then, a Morawetz vector field is defined from a different direction than the Killing field as in \eqref{K}, which can be turned in the direction of the Killing vector field. In \eqref{JK}, we construct an energy functional derived by applying the divergence theorem to the Morawetz vector field contracted with MH energy momentum tensor. We also have the flux $E^K$ of the vector field through the closed surface as in equation \eqref{EK}. In preposition \ref{prep2}, we demonstrate that \eqref{JK} can be estimated in the term of \eqref{EK}. Furthermore, in preposition \ref{prep3}, we prove that the flux can be bounded by the initial energy of the system $\hat{E}^{MH}$. To prove the decay estimates, we use the Sobolev inequality on the manifold for the gauge field, the complex scalar field, and their interaction that represents by the covariant derivative of the scalar field. The energy estimate method, combined with the H\"{o}lder inequality, is employed to obtain the decay of the Maxwell-Higgs system for far away from the horizon.

In the second part, to get the decay near the horizon, we define a vector field $H$ as in equation \eqref{H}. A spacetime integral  supported on a bounded region in space near the event horizon is construct to control the $L^2$ norm of certain components of the fields as in \eqref{IH}. We also define the flux of \eqref{H} as in definition \ref{defnear1}, see \eqref{eH1} and \eqref{eH2}. We demonstrate that the flux bound by the energy generated from a time-like vector field in preposition \ref{prep4}. Near horizon, one need to choose condition as in \eqref{h1}-\eqref{h4} such that  \eqref{IH} bound by the initial energy as in preposition \ref{prep8}. Then, with the help of Cauchy stability, the flux can be estimated from the initial data on the Cauchy hypersurface as in preposition \ref{prep9}. Finally, to obtain the decay for the normalised components, the Sobolev inequality is used together with the energy estimation method, which is bound to hypersurfaces of $w = constant$ and a fixed length in $v$.

Thus, the summary of our works can be stated in the following two theorems

\begin{theorem}[Far Away From Horizon] \label{teorema1}
	Let $(M, g)$ be a maximally extended Schwarzschild spacetimes. Let $\Sigma_{t=t_0}$ denote a Cauchy hypersurface defined by $t = t_0$. Let $F_{\mu\nu}\left(w, v,\theta,\varphi,\right)$ be the Maxwell field components and $\phi$ be the complex scalar field in a coordinate system described as the solution of the Cauchy problem of the MH equations in (\ref{eom1}) and (\ref{eom2}). Suppose that for all $t_i = (1.1)^it_0$ and for away from the horizon (outside photon sphere region, $r\geq 3m$), the decay of the MH equation is obtained as follows, 
\begin{eqnarray}\label{decayF1}
	\left| {{F}_{\hat{\mu }\hat{\nu }}} \right|\left( w,v,\theta ,\varphi  \right)\lesssim \frac{{\hat{E}^{MH}}}{1+\left| v \right|}~,
\end{eqnarray}
\begin{eqnarray}\label{decayF2}
	\left| {{F}_{\hat{\mu }\hat{\nu }}} \right|\left( w,v,\theta ,\varphi  \right)\lesssim \frac{{\hat{E}^{MH}}}{1+\left| w \right|}~,
\end{eqnarray}
\begin{eqnarray}\label{decayscalar}
\left| \phi  \right|\lesssim \frac{{\hat{E}^{MH}}}{1+\left| v \right|}~.
\end{eqnarray}
We also have the decay of covariant derivative of the scalar field as follows
\begin{eqnarray}\label{decayDphi}
	\left| D \phi  \right|\lesssim \frac{\left(\hat{E}^{MH}\right)^2+\left(\hat{E}^{MH}\right)^{5/4}}{1+\left|v\right|}
\end{eqnarray}
 where
\begin{eqnarray}
	\hat{E}^{MH}&=&\sum\limits_{i=0}^{1}{\sum\limits_{j=0}^{5}{E_{{{r}^{j}}{{\left( \slashed{\mathcal{L}} \right)}^{j}}{{\left( {{\mathcal{L}}_{t}} \right)}^{i}}}^{\hat{t}}\left( {{t}_{0}} \right)}}+\sum\limits_{i=0}^{1}{\sum\limits_{j=0}^{4}{E_{{{r}^{j}}{{\left( \slashed{\mathcal{L}} \right)}^{j}}{{\left( {{\mathcal{L}}_{t}} \right)}^{i}}}^{K}\left( {{t}_{0}} \right)}} \notag\\ 
	&+&E_{{{r}^{6}}{{\left( \slashed{\mathcal{L}} \right)}^{6}}}^{\hat{t}}\left( {{t}_{0}} \right)+E_{{{r}^{5}}{{\left( \slashed{\mathcal{L}} \right)}^{5}}}^{ K }\left( {{t}_{0}} \right).
\end{eqnarray}
 with $\slashed{\mathcal{L}}$ is the Lie derivative restricted on the 2-spheres.
 \end{theorem}
 
 \begin{theorem}[Near Horizon]\label{teorema2}
Let $H$ be a vector field as in \eqref{H} and $\mu=\frac{2m}{r}$. Suppose that we have \eqref{h1}-\eqref{h4}, then the decay estimate is obtained in the photon sphere region ($2m<r<3m$) as
 \begin{eqnarray}\label{F1}
    {{\left| {{F}_{\hat{v}\hat{w}}} \right|}^{2}}\lesssim {{\left( \frac{w}{{{v}_{+}}} \right)}^{2}}{{E}_{1}}~,
\end{eqnarray}
\begin{eqnarray}\label{F2}
    {{\left| {{F}_{{{e}_{1}}{{e}_{2}}}} \right|}^{2}}\lesssim {{\left( \frac{w}{{{v}_{+}}} \right)}^{2}}{{E}_{1}}~,
\end{eqnarray}
\begin{eqnarray}\label{F3}
    \left| {{F}_{\hat{v}\hat{\varphi }}} \right|^2\lesssim {{\left( \frac{w}{{{v}_{+}}} \right)}^{2}}{{E}_{2}}~,
\end{eqnarray}
\begin{eqnarray}\label{F4}
    {{\left| \sqrt{1-\mu }{{F}_{\hat{w}\hat{\theta }}} \right|}^{2}}\lesssim {{\left( \frac{w}{{{v}_{+}}} \right)}^{2}}{{E}_{2}}~,
\end{eqnarray}
\begin{eqnarray}\label{F5}
    {{\left| \sqrt{1-\mu }{{F}_{\hat{w}\hat{\varphi }}} \right|}^{2}}\lesssim {{\left( \frac{w}{{{v}_{+}}} \right)}^{2}}{{E}_{2}}~,
\end{eqnarray}
\begin{eqnarray}\label{pi1}
    \left| \phi  \right|\lesssim {{\left( 1+{{\left( \frac{w}{{{v}_{+}}} \right)}^{2}} \right)}^{1/2}}{{E}_{3}}^{1/2}~,
\end{eqnarray}
with ${{v}_{+}}=\max \left\{ 1,v \right\}$. More completely, we also get the decay for the covariant derivative of scalar field,
\begin{eqnarray}\label{Dpi}
    \left| D\phi  \right|\lesssim E_{4}^{1/2}+\left( 1+{{\left( \frac{w}{{{v}_{+}}} \right)}^{2}} \right){{E}_{3}}^{1/2}{{E}_{1}}^{1/2}~,
\end{eqnarray}
where $E_1, E_2, E_3,$ and $E_4$ are the energy that depend only on $E^{\hat{t}}\left(t_0\right)$ and $E^K\left(t_0\right)$ as in \eqref{E1}, \eqref{E2}, \eqref{E3}, \eqref{E4} respectively.
\end{theorem}

We organize this paper as follows. The energy functional and the estimates for the energy with respect to the time-like Killing field, the Morawetz vector field, and the reduced energy are discussed in Section \ref{sec:MH energy}. In Section \ref{sec:decayfaraway}, we prove the decay estimate for MH outside the photon sphere region. At the end, we discuss the decay estimates inside the photon sphere region in Section \ref{sec:decaynear}.

\section{The energy }
\label{sec:MH energy}

Let $V^\nu$ be a vector field. We first introduce
\begin{equation}\label{J}
	{{J}_{\mu }} \equiv {{V}^{\nu }}{T_{\mu \nu }}~.
\end{equation}
Considering the Maxwell field $F$ and Higgs field $\phi$ as it satisfies (\ref{eom1}) and (\ref{eom2}), it can be shown
\begin{equation}
  \nabla^\nu T_{\mu\nu}\left(F,\phi\right)=0 ~. 
\end{equation}
Then, by applying the divergence theorem to $J$ in the region $\mathcal{B}={{J}^{+}}\left( {{\Sigma }_{t_1}} \right)\cap {{J}^{-}}\left( {\Sigma }_{{t_{2}}} \right)\cap \mathcal{D}$, where $\mathcal{D}$ is the outer-communication of Schwarzschild black holes, bounded by a null hypersurface $N$, as well as to the past by ${{\Sigma }_{{{t}_{1}}}}$ and to the future by ${{\Sigma }_{{{t}_{2}}}}$, one can demonstrate
\begin{eqnarray}\label{con}
\int_{\mathcal{B}}{{{\pi }^{\mu \nu }}{T_{\mu \nu }}(V)d{{V}_{\mathcal{B}}}}=\int_{{{\Sigma }_{{{{{t}'}}_{1}}}}}{{{J}_{\mu }}{{n}^{\mu }}d{{V}_{{{\Sigma }_{{{{{t}}}_{1}}}}}}}-\int_{{{\Sigma }_{{t}2}}}{{{J}_{\mu }}{{n}^{\mu }}d{{V}_{{{\Sigma }_{{{{{t}}}_{2}}}}}}}-\int_{N}{{{J}_{\mu }}n_{N}^{\mu }d{{V}_{N}}}~,
\end{eqnarray}	
where $n_\mu$ are the unit normal to the hypersurfaces $\Sigma$ and $\pi^{\mu\nu} \equiv \nabla^{\mu}V^{\nu}$. 

Let us consider the timelike vector field $V=\frac{\partial}{\partial t}$ such that we have
\begin{eqnarray}\label{Ttime}
	{{T}_{\alpha \beta }}{{\pi }^{\alpha \beta }}\left( \frac{\partial }{\partial t} \right)&=&-\left\{ \frac{1}{{{r}^{2}}}{{\left| {{F }_{w\theta }} \right|}^{2}}+\frac{1}{{{r}^{2}}{{\sin }^{2}}\theta }{{\left| {{F }_{w\varphi }} \right|}^{2}}+{{\left| {{D}_{w}}\phi  \right|}^{2}} \right\}\frac{2}{\left( 1-\mu  \right)}{{\partial }_{v}}{{t}^{w}} \notag\\
	&-&\left\{ \frac{1}{{{r}^{2}}}{{\left| {{F }_{v\theta }} \right|}^{2}}+\frac{1}{{{r}^{2}}{{\sin }^{2}}\theta }{{\left| {{F }_{v\varphi }} \right|}^{2}}+{{\left| {{D}_{v}}\phi  \right|}^{2}} \right\}\frac{2}{\left( 1-\mu  \right)}{{\partial }_{w}}{{t}^{v}} \notag\\
	&-&2\left\{ \frac{1}{{{(1-\mu )}^{2}}}{{\left| {{F }_{vw}} \right|}^{2}}+\frac{1}{4{{r}^{4}}{{\sin }^{2}}\theta }{{\left| {{F }_{\theta \varphi }} \right|}^{2}} \right\}\left( {{\partial }_{v}}{{t}^{v}}+{{\partial }_{w}}{{t}^{w}}+\frac{\left( 3\mu -2 \right)}{2r}\left[ {{t}^{v}}-{{t}^{w}} \right] \right) \notag\\ 
	&-&\left\{ \frac{2}{\left( 1-\mu  \right)}{{D}_{v}}\phi \overline{{{D}_{w}}\phi }+\frac{1}{2{{r}^{2}}}{{\left| {{D}_{\theta }}\phi  \right|}^{2}}+\frac{1}{2{{r}^{2}}{{\sin }^{2}}\theta }{{\left| {{D}_{\varphi }}\Phi  \right|}^{2}}+P \right\} \notag\\ 
	&\times& \left( {{\partial }_{v}}{{t}^{v}}+{{\partial }_{w}}{{t}^{w}}+\frac{\mu }{2r}\left[ {{t}^{v}}-{{t}^{w}} \right] \right)+\frac{\left( 1-\mu  \right)\left( 2+{{\sin }^{2}}\theta  \right)}{4r}P\left[ {{t}^{v}}-{{t}^{w}} \right] \notag\\ 
	&=&0~.  
\end{eqnarray}
where $\mu=\frac{2m}{r}$ and $t^ v$ is the time-like vector field. We have used the Eddington-Finkelstein coordinates $v$ and $w$ by the formulas in \eqref{EFw} and \eqref{EFv}, followed by the Regge-Wheeler tortoise coordinate $r^*$ as in \eqref{RW}.
\begin{definition}\label{definisi1}
The energy $E^{\hat{t}}$ is defined as
\begin{eqnarray}
    E^{\hat{t}}\left( t={{t}_{i}} \right)\equiv\int_{{{\Sigma }_{{{{{t}'}}_{i}}}}}{{{J}_{\mu }}\left(\frac{\partial}{\partial t}\right)^{\mu}d{{V}_{{{\Sigma }_{{{{{t}'}}_{i}}}}}}}~.
\end{eqnarray}
\end{definition}
Thus, from (\ref{con}) we get the conservation of energy generated by the time-like vector field $\frac{\partial }{\partial t}$
\begin{eqnarray}
	E^{\hat{t}}\left( t={{t}_{2}} \right)-E^{\hat{t}}\left( t={{t}_{1}} \right)=0~, 
\end{eqnarray}
where 
\begin{eqnarray}\label{energytimelike}
	E^{\hat{t}}\left( t={{t}_{i}} \right)&=&\int_{{{\Sigma }_{{{t}_{i}}}}}{\left( 1-\mu  \right){{r}^{2}}d{{r}^{*}}d{{\sigma }^{2}}\left\{ {{\left| {{D}_{{\hat{t}}}}\phi  \right|}^{2}}+{{\left| {{D}_{{\hat{\theta }}}}\phi  \right|}^{2}}+{{\left| {{D}_{{\hat{\varphi }}}}\phi  \right|}^{2}}+{{\left|P \right|}} \right.} \notag\\ 
	&+&\left. \frac{1}{2}\left( {{\left| {{F }_{\hat{t}{{{\hat{r}}}^{*}}}} \right|}^{2}}+{{\left| {{F }_{\hat{t}\hat{\theta }}} \right|}^{2}}+{{\left| {{F }_{\hat{t}\hat{\varphi }}} \right|}^{2}}+{{\left| {{F }_{{{{\hat{r}}}^{*}}\hat{\theta }}} \right|}^{2}}+{{\left| {{F }_{{{{\hat{r}}}^{*}}\hat{\varphi }}} \right|}^{2}}+{{\left| {{F }_{\hat{\varphi }\hat{\theta }}} \right|}^{2}} \right) \right\}~. \notag\\
\end{eqnarray}
Here we used the notation $\frac{{\hat{\partial }}}{\partial t}=\frac{1}{\sqrt{1-\mu }}\frac{\partial }{\partial t}$. Furthermore, on the Schwarzschild spacetimes  the angular momentum operators are Killing vectors so that $\mathcal{L}_{\Omega_j} F$ and $\mathcal{L}_{\Omega_j} \phi $ satisfy $\nabla^\nu T_{\mu\nu}\left(\mathcal{L}_{\Omega_j} F,\mathcal{L}_{\Omega_j} \phi\right)=0$, where $\Omega_j$ is a basis of angular momentum operators. Thus, we also have $E^{\hat{t}}$ evaluated on the pair $\left( {{\mathcal{L}}_{{{\Omega }_{j}}}}F,{{\mathcal{L}}_{{{\Omega }_{j}}}}\phi  \right)$,

\begin{eqnarray}\label{energyangularij}
	E^{\hat{t}}_{\mathcal{L}_{\Omega_j}}\left( t_i\right)&=&\int_{{{\Sigma }_{{{t}_{i}}}}}{\left( 1-\mu  \right){{r}^{2}}d{{r}^{*}}d{{\sigma }^{2}}\left\{ {{\left| \mathcal{L}_{\Omega_j}\left({{D}_{{\hat{t}}}}\phi \right) \right|}^{2}}+{{\left| \mathcal{L}_{\Omega_j}\left({{D}_{{\hat{\theta }}}} \phi \right) \right|}^{2}}+{{\left| \mathcal{L}_{\Omega_j} \left({{D}_{{\hat{\varphi }}}}\phi  \right)\right|}^{2}}+{{\left|\mathcal{L}_{\Omega_j} P \right|}} \right.} \notag\\ 
	&+&\left. \frac{1}{2}\left( {{\left| {{\mathcal{L}_{\Omega_j} F }_{\hat{t}{{{\hat{r}}}^{*}}}} \right|}^{2}}+{{\left| {{\mathcal{L}_{\Omega_j} F }_{\hat{t}\hat{\theta }}} \right|}^{2}}+{{\left| {{\mathcal{L}_{\Omega_j} F }_{\hat{t}\hat{\varphi }}} \right|}^{2}}+{{\left| {{\mathcal{L}_{\Omega_j} F }_{{{{\hat{r}}}^{*}}\hat{\theta }}} \right|}^{2}}+{{\left| {{\mathcal{L}_{\Omega_j} F }_{{{{\hat{r}}}^{*}}\hat{\varphi }}} \right|}^{2}}+{{\left| {{\mathcal{L}_{\Omega_j} F }_{\hat{\varphi }\hat{\theta }}} \right|}^{2}} \right) \right\}~. \notag\\
\end{eqnarray}

\subsection{The Morawetz vector field energy}
\label{sec:Morawetz}
\begin{definition}\label{definisi2}
   We introduce the Morawetz vector field $K$  as
   \begin{eqnarray}\label{K}
	K \equiv -{{w}^{2}}\frac{\partial }{\partial w}-{{v}^{2}}\frac{\partial }{\partial v} ~ .  
\end{eqnarray}
\end{definition}
\begin{definition}\label{definisi3}
   The quantity $\mathcal{J}^{K}$ is the the space-time integral derived by applying the divergence theorem to the vector field $K$ contracted with the MH energy momentum tensor 
\begin{eqnarray}
	{{\mathcal{J}}^{K}}\left( {{t}_{i}}\le t\le {{t}_{i+1}} \right)\equiv\int_{{{t}_{i}}}^{{{t}_{i+1}}}{\int_{{{r}^{*}}=-\infty }^{{{r}^{*}}=\infty }{\int\limits_{{{S}^{2}}}{{{\pi }^{\alpha \beta }}\left( K \right){{T}_{\alpha \beta }}\left( F,\phi  \right)\text{dVol}}}}~.
\end{eqnarray} 
\end{definition}
\begin{definition}\label{definisi4}
The functional energy $E^K$ is the flux of the vector field $K$ through the closed surface defined by
\begin{eqnarray}
	{{E}^{K}}\left( {{t}_{i}}\le t\le {{t}_{i+1}} \right)\equiv\int_{{{r}^{*}}=-\infty }^{{{r}^{*}}=\infty }{\int\limits_{{{S}^{2}}}{{{J}_{\alpha }}\left( K \right){{n}^{\alpha }}\text{dVo}{{\text{l}}_{t={{t}_{i}}}}}}~.
\end{eqnarray}
\end{definition}
Thus, one can show
\begin{eqnarray}\label{JK}
	{{\mathcal{J}}^{K}}\left( {{t}_{i}}\le t\le {{t}_{i+1}} \right)&=&\int_{{{t}_{i}}}^{{{t}_{i+1}}}{\int_{{{r}^{*}}=-\infty }^{{{r}^{*}}=\infty }{\int\limits_{{{S}^{2}}}{\left[ \left\{ \frac{1}{{{(1-\mu )}^{2}}}{{\left| {{F }_{vw}} \right|}^{2}}+\frac{1}{4{{r}^{4}}{{\sin }^{2}}\theta }{{\left| {{F }_{\theta \varphi }} \right|}^{2}} \right\}\left( 4+\frac{\left( 3\mu -2 \right)}{r}{{r}^{*}} \right) \right.}}} \notag\\ 
	&+&\left\{ \frac{2}{\left( 1-\mu  \right)}{{D}_{v}}\phi \overline{{{D}_{w}}\phi }+\frac{1}{2{{r}^{2}}}{{\left| {{D}_{\theta }}\phi  \right|}^{2}}+\frac{1}{2{{r}^{2}}{{\sin }^{2}}\theta }{{\left| {{D}_{\varphi }}\phi  \right|}^{2}}+P \right\}\left( 2+\frac{\mu }{r}{{r}^{*}} \right) \notag\\ 
	&+& \left.\frac{\left( 1-\mu  \right)\left( 2+{{\sin }^{2}}\theta  \right)}{2r}P{{r}^{*}} \right]2t{{r}^{2}}\left( 1-\mu  \right)d{{\sigma }^{2}}d{{r}^{*}}dt~,  
\end{eqnarray}
\begin{eqnarray}\label{EK}
	{{E}^{K}}\left( {{t}_{i}}\le t\le {{t}_{i+1}} \right)&=&\int_{{{r}^{*}}=-\infty }^{{{r}^{*}}=\infty }{\int\limits_{{{S}^{2}}}{\left( {{w}^{2}}\left[ \frac{1}{{{r}^{2}}}{{\left| {{F}_{w\theta }} \right|}^{2}}+\frac{1}{{{r}^{2}}{{\sin }^{2}}\theta }{{\left| {{F }_{w\varphi }} \right|}^{2}}+{{\left| {{D}_{w}}\phi  \right|}^{2}} \right] \right.}} \notag\\ 
	&+&{{v}^{2}}\left[ \frac{1}{{{r}^{2}}}{{\left| {{F }_{v\theta }} \right|}^{2}}+\frac{1}{{{r}^{2}}{{\sin }^{2}}\theta }{{\left| {{F}_{v\varphi }} \right|}^{2}}+{{\left| {{D}_{v}}\phi  \right|}^{2}} \right] \notag\\ 
	&+&\left( {{v}^{2}}+{{w}^{2}} \right)\left[ \frac{1}{{{(1-\mu )}^{2}}}{{\left| {{F}_{vw}} \right|}^{2}}+\frac{(1-\mu )}{4{{r}^{4}}{{\sin }^{2}}\theta }{{\left| {{F }_{\varphi \theta }} \right|}^{2}}+{{D}_{v}}\phi \overline{{{D}_{w}}\phi } \right. \notag\\ 
	&+&\left. \left. \frac{\left( 1-\mu  \right)}{4{{r}^{2}}}{{\left| {{D}_{\theta }}\phi  \right|}^{2}}+\frac{\left( 1-\mu  \right)}{4{{r}^{2}}{{\sin }^{2}}\theta }{{\left| {{D}_{\varphi }}\phi  \right|}^{2}}+\frac{\left( 1-\mu  \right)}{2}P\right] \right){{r}^{2}}d{{r}^{*}}d{{\sigma }^{2}}~.\notag\\  
\end{eqnarray}
It can be observed that from \eqref{JK}, ${{\mathcal{J}}^{K}}\left( {{t}_{i}}\le t\le {{t}_{i+1}} \right)$ can only be positive on the interval $[r_0,R_0]$ where $r_0> 2m$ and $R_0>3m$.

To control the estimate of $\mathcal{J}^{K}$ and $E^K$, a radial vector field $G$ is required to estimate the conformal energy and its derivative. 
\begin{definition}
The radial vector field $G$ is defined as
\begin{eqnarray}\label{G}
	G \equiv-f\left( {{r}^{*}} \right)\frac{\partial }{\partial w}+f\left( {{r}^{*}} \right)\frac{\partial }{\partial v}~,
\end{eqnarray}
where $f$ is bounded function. 
\end{definition}
To accomplish this, it is necessary to have the space-time integral $\mathcal{J}^{G}$, which is referred to as the functional energy and contains positive terms.
\begin{definition}\label{definisi6}
Introducing $\mathcal{J}^{G}$ whose form is given by
\begin{eqnarray}
	{{\mathcal{J}}^{G}}\left( {{t}_{i}}\le t\le {{t}_{i+1}} \right)&\equiv&\int_{{{t}_{i}}}^{{{t}_{i+1}}}{\int_{{{r}^{*}}=r_{0}^{*}}^{{{r}^{*}}=R_{0}^{*}}{\int\limits_{{{S}^{2}}}{\left\{ {{\left| {{F }_{\hat{v}\hat{w}}} \right|}^{2}}+{{\left| {{F }_{\hat{\varphi }\hat{\theta }}} \right|}^{2}}+{{\left| {{D}_{{\hat{v}}}}\phi  \right|}^{2}}+{{\left| {{D}_{{\hat{w}}}}\phi  \right|}^{2}} \right.}}} \notag\\ 
	&+&\left. {{\left| {{D}_{{\hat{\theta }}}}\phi  \right|}^{2}}+{{\left| {{D}_{{\hat{\varphi }}}}\phi  \right|}^{2}}+P \right\}d{{r}^{*}}d{{\sigma }^{2}}dt~.
\end{eqnarray}
\end{definition} 
\begin{definition}
The quantity $E^G$ is the flux of the vector field $G$ through the closed surface defined by

\begin{eqnarray}\label{defEG}
	{{E}^{G}}\left( {{t}_{i}}\le t\le {{t}_{i+1}} \right)\equiv\int_{{{r}^{*}}=-\infty }^{{{r}^{*}}=\infty }{\int\limits_{{{S}^{2}}}{{{J}_{\alpha }}\left( G \right){{n}^{\alpha }}\text{dVo}{{\text{l}}_{t={{t}_{i}}}}}}~.
\end{eqnarray}
\end{definition}

\subsection{The energy estimate}
\label{sec:Energy}
\begin{preposition}\label{prep1}
    Let $f(r^{*})$ in \eqref{G} be a bounded function, then we can attain
    \begin{eqnarray}
        \left|E^G \left(t_i\right)\right|\lesssim\left|\Tilde{E}^{\hat{t}} \left(t_i\right)\right|\lesssim\left|{E}^{\hat{t}} \left(t_i\right)\right|
    \end{eqnarray}
\end{preposition}
\begin{proof}
    From \eqref{defEG}, one can show that
  \begin{eqnarray}\label{EG}
	{{E}^{G}}\left( t \right)=\int_{{{r}^{*}}=-\infty }^{{{r}^{*}}=\infty }{\int\limits_{{{S}^{2}}}{-f\left\{ \frac{1}{{{r}^{2}}}{{F }_{t\theta }}{{F }_{{{r}^{*}}\theta }}+\frac{1}{{{r}^{2}}{{\sin }^{2}}\theta }{{F }_{t\varphi }}{{F }_{{{r}^{*}}\varphi }}+{{D}_{t}}\phi \overline{{{D}_{{{r}^{*}}}}\phi } \right\}{{r}^{2}}d{{r}^{*}}d{{\sigma }^{2}}}}~.
\end{eqnarray}  
It can be demonstrated that
    \begin{eqnarray}
\left| {{E}^{G}}\left( {{t}_{i}} \right) \right|&=&\int_{{{r}^{*}}=-\infty }^{{{r}^{*}}=\infty }{\int\limits_{{{S}^{2}}}{\left| -f \right|\left\{ \frac{1}{{{r}^{2}}}\left| {{F}_{t\theta }}{{F}_{{{r}^{*}}\theta }} \right|+\frac{1}{{{r}^{2}}{{\sin }^{2}}\theta }\left| {{F}_{t\varphi }}{{F}_{{{r}^{*}}\varphi }} \right|+\left| {{D}_{t}}\phi \overline{{{D}_{{{r}^{*}}}}\phi } \right| \right\}{{r}^{2}}d{{r}^{*}}d{{\sigma }^{2}}}} \notag\\ 
 &\lesssim& \int_{{{r}^{*}}=-\infty }^{{{r}^{*}}=\infty }{\int\limits_{{{S}^{2}}}{\left| -f \right|\left\{ \frac{1}{{{r}^{2}}\left( 1-\mu  \right)}{{\left| {{F}_{t\theta }} \right|}^{2}}+\frac{1}{{{r}^{2}}\left( 1-\mu  \right)}{{\left| {{F}_{{{r}^{*}}\theta }} \right|}^{2}}+\frac{1}{\left( 1-\mu  \right){{r}^{2}}{{\sin }^{2}}\theta }{{\left| {{F}_{t\varphi }} \right|}^{2}} \right.}} \notag\\ 
 && + \left. \frac{1}{\left( 1-\mu  \right){{r}^{2}}{{\sin }^{2}}\theta }{{\left| {{F}_{{{r}^{*}}\varphi }} \right|}^{2}}+\frac{1}{\left( 1-\mu  \right)}{{\left| {{D}_{t}}\phi  \right|}^{2}}+\frac{1}{\left( 1-\mu  \right)}{{\left| {{D}_{{{r}^{*}}}}\phi  \right|}^{2}} \right\}\frac{1}{\left( 1-\mu  \right)}{{r}^{2}}d{{r}^{*}}d{{\sigma }^{2}} \notag\\ 
 &\lesssim&\left| {{{\tilde{E}}}^{{\hat{t}}}}\left( t={{t}_{i}} \right) \right| \notag\\ 
 &\lesssim& \left| {{E}^{{\hat{t}}}}\left( t={{t}_{i}} \right) \right| ~. 
\end{eqnarray}
\end{proof}
\begin{preposition}\label{prep2}
	Suppose $\mathcal{J}^K$ as in definition \ref{definisi3}. Then, for $t_{i+1}\le t_{i}+0.1t_i$ and $\left| {{r}^{*}}\left( {{r}_{0}} \right) \right|+\left| {{r}^{*}}\left( R \right) \right|\le 0.4{{t}_{i}}$, we have
	\begin{eqnarray}\label{esJK}
		{{\mathcal{J}}^{K}}\left( {{t}_{i}}\le t\le {{t}_{i+1}} \right)&\lesssim& \text{ }{{t}_{i+1}}{{\mathcal{J}}^{G}}\left( {{t}_{i}}\le t\le {{t}_{i+1}} \right)\left( {{r}_{0}}\le r\le {{R}_{0}} \right) \notag\\ 
		&\lesssim& \text{ }{{t}_{i+1}}\left\{ \frac{1}{t_{i}^{2}}{{E}^{K}}\left( {{t}_{i}} \right)+\frac{1}{t_{i}^{2}}\sum\limits_{j=1}^{3}{E_{{{\mathcal{L}}_{{{\Omega }_{j}}}}}^{K}\left( {{t}_{i}} \right)+}\frac{1}{t_{i+1}^{2}}{{E}^{K}}\left( {{t}_{i+1}} \right)+\frac{1}{t_{i+1}^{2}}\sum\limits_{j=1}^{3}{E_{{{\mathcal{L}}_{{{\Omega }_{j}}}}}^{K}\left( {{t}_{i+1}} \right)} \right\} ~. \notag\\
	\end{eqnarray}
\end{preposition}
\begin{proof}
Remark that from (\ref{JK}), $\left( 4+\frac{\left( 3\mu -2 \right)}{r}{{r}^{*}} \right)$ is positive only in a bounded interval $[r_0,R_0]$ where $r_0> 2m$ and $R_0>3m$. Then, 
	\begin{eqnarray}
		\mathcal{J}^K\left( {{t}_{i}}\le t\le {{t}_{i+1}} \right)&\lesssim& {{t}_{i+1}}\int_{{{t}_{i}}}^{{{t}_{i+1}}}{\int_{{{r}^{*}}={{r}_{0}}}^{{{r}^{*}}={{R}_{0}}}{\int\limits_{{{S}^{2}}}{\left\{ {{\left| {{F }_{\hat{v}\hat{w}}} \right|}^{2}}+\frac{1}{4}{{\left| {{F }_{\hat{\theta }\hat{\varphi }}} \right|}^{2}}+{{\left| {{D}_{{\hat{v}}}}\phi  \right|}^{2}}+{{\left| {{D}_{{\hat{w}}}}\phi  \right|}^{2}} \right.}}} \notag\\ 
		&& +\left. \frac{1}{2}{{\left| {{D}_{{\hat{\theta }}}}\phi  \right|}^{2}}+\frac{1}{2}{{\left| {{D}_{{\hat{\varphi }}}}\phi  \right|}^{2}}+P \right\}d{{\sigma }^{2}}d{{r}^{*}}dt \notag\\ 
		&\lesssim& \text{ }{{t}_{i+1}}\mathcal{J}^G\left( {{t}_{i}}\le t\le {{t}_{i+1}} \right)\left( r_{0}^{*}\le {{r}^{*}}\le R_{0}^{*} \right)~.  
	\end{eqnarray}
Using the inequality (8) in \cite{Ghanem1}, we have
	\begin{eqnarray}\label{asumeq}
		{{\mathcal{J}}^{G}}\left( {{t}_{i}}\le t\le {{t}_{i+1}} \right)\lesssim \left| {{\Tilde{E}}^{\hat{t}}}\left( {{t}_{i+1}} \right) \right|+\left| {{\Tilde{E}}^{\hat{t}}}\left( {{t}_{i}} \right) \right|+\sum\limits_{j=1}^{3}{\left( \left| \Tilde{E}_{{{\mathcal{L}}_{{{\Omega }_{j}}}}}^{\hat{t}}\left( {{t}_{i+1}} \right) \right|+\left| \Tilde{E}_{{{\mathcal{L}}_{{{\Omega }_{j}}}}}^{\hat{t}}\left( {{t}_{i}} \right) \right| \right)}~,\notag\\
	\end{eqnarray}
 where ${\Tilde{E}}^{\hat{t}}$ is called "reduced energy" given by
 \begin{eqnarray}\label{energyasumtimelike}
	\Tilde{E}^{\hat{t}}\left( t={{t}_{i}} \right) &\equiv& \int_{{{\Sigma }_{{{t}_{i}}}}}{\left( 1-\mu  \right){{r}^{2}}d{{r}^{*}}d{{\sigma }^{2}}\left\{ {{\left| {{D}_{{\hat{t}}}}\phi  \right|}^{2}}+{{\left| {{D}_{{\hat{\theta }}}}\phi  \right|}^{2}}+{{\left| {{D}_{{\hat{\varphi }}}}\phi  \right|}^{2}}+{{\left|P \right|}} \right.} \notag\\ 
	&+&\left. \frac{1}{2}\left( {{\left| {{F }_{\hat{t}\hat{\theta }}} \right|}^{2}}+{{\left| {{F }_{\hat{t}\hat{\varphi }}} \right|}^{2}}+{{\left| {{F }_{{{{\hat{r}}}^{*}}\hat{\theta }}} \right|}^{2}}+{{\left| {{F }_{{{{\hat{r}}}^{*}}\hat{\varphi }}} \right|}^{2}} \right) \right\}~. \notag\\
\end{eqnarray}
Using \eqref{asumeq}, we obtain
	\begin{eqnarray}\label{es1}
		{{\mathcal{J}}^{K}}\left( {{t}_{i}}\le t\le {{t}_{i+1}} \right)\lesssim{{t}_{i+1}}\left\{ \left| {{\Tilde{E}}^{\hat{t}}}\left( {{t}_{i+1}} \right) \right|+\left| {{\Tilde{E}}^{\hat{t}}}\left( {{t}_{i}} \right) \right|+\sum\limits_{j=1}^{3}{\left( \left| \Tilde{E}_{{{\mathcal{L}}_{{{\Omega }_{j}}}}}^{\hat{t}}\left( {{t}_{i+1}} \right) \right|+\left| \Tilde{E}_{{{\mathcal{L}}_{{{\Omega }_{j}}}}}^{\hat{t}}\left( {{t}_{i}} \right) \right| \right)}\right\}~.\notag\\
	\end{eqnarray}
Let $\hat{F}$ denotes the anti-symmetric tensor which satisfies
	\begin{equation}
		{{\hat{F }}_{\hat{\mu }\hat{\nu }}}\left( {{t}_{i}},{{r}^{*}},\theta ,\varphi  \right)=\hat{\chi }\left( \frac{2{{r}^{*}}}{{{t}_{i}}} \right){{F }_{\hat{\mu }\hat{\nu }}}\left( {{t}_{i}},{{r}^{*}},\theta ,\varphi  \right)~,
	\end{equation}
where $\hat{\chi }$ is a smooth cut-off function equal to one on $\left[-1,1\right]$ and zero outside $\left[-\frac{3}{2},\frac{3}{2}\right]$. We could see that for ${{r}^{*}}\le -\frac{3{{t}_{i}}}{4}$ and ${{r}^{*}}\ge \frac{3{{t}_{i}}}{4}$ 
\begin{eqnarray}
		{{\hat{F }}_{kl}}\left( {{t}_{i}},{{r}^{*}},\theta ,\varphi  \right)=0~,
	\end{eqnarray}
where $\left( k,l \right)\in \left\{ \left( {{r}^{*}},\theta  \right),\left( {{r}^{*}},\varphi  \right),\left( t,\theta  \right),\left( t,\varphi  \right) \right\}$. For $-\frac{{{t}_{i}}}{2}\le {{r}^{*}}\le \frac{{{t}_{i}}}{2}$ and $t_{i}\le t\le t_{i+1}$, we have
	\begin{eqnarray}
		{{\hat{F }}_{\mu \nu }}\left( {{t}_{i}},{{r}^{*}},\theta ,\varphi  \right)={{F }_{\mu \nu }}\left( {{t}_{i}},{{r}^{*}},\theta ,\varphi  \right)~.
	\end{eqnarray}

Now consider region $t_i<t<t_{i+1}$ and $r_0<r<R_0$ where $t_{i+1}\le t_{i}+0.1t_i$ and $\left| {{r}^{*}}\left( {{r}_{0}} \right) \right|+\left| {{r}^{*}}\left( R \right) \right|\le 0.4{{t}_{i}}$, on $t=t_i$ we have ${{\hat{F }}_{\mu \nu }}={{F }_{\mu \nu }}$ in the specified region. At $t=t_{i+1}$ we have ${{\hat{F}}_{\mu \nu }}\left( t={{t}_{i+1}} \right)=0$ if $r_{0}^{*}\ge -\frac{3{{t}_{i}}}{4}-0.1{{t}_{i}}=-0.85{{t}_{i}}$ and $r_{0}^{*}\le \frac{3{{t}_{i}}}{4}+0.1{{t}_{i}}=0.85{{t}_{i}}$. Therefore, by using the boundedness of $\hat{\chi}$, we have
	\begin{eqnarray}\label{redenergy}
	\Tilde{E}^{\hat{t}}\left( {{t}_{i+1}} \right)&\lesssim& \int_{{{r}^{*}}=-0.85{{t}_{i}}}^{{{r}^{*}}=0.85{{t}_{i}}}{\int_{{{S}^{2}}}{\left( {{\left| {{{\hat{F }}}_{\hat{v}\hat{\theta }}} \right|}^{2}}+{{\left| {{{\hat{F }}}_{\hat{v}\hat{\varphi }}} \right|}^{2}}+{{\left| {{{\hat{F }}}_{\hat{w}\hat{\theta }}} \right|}^{2}}+{{\left| {{{\hat{F }}}_{\hat{w}\hat{\varphi }}} \right|}^{2}}+{{\left| {{D}_{{\hat{v}}}}\hat{\phi}  \right|}^{2}} \right.}} \notag\\ 
		&& + \left. {{\left| {{D}_{{\hat{w}}}}\hat{\phi} \right|}^{2}}+{{\left| {{D}_{{\hat{\theta }}}}\hat{\phi}  \right|}^{2}}+{{\left| {{D}_{{\hat{\varphi }}}}\hat{\phi}  \right|}^{2}}+\left|P\right| \right){{r}^{2}}\left( 1-\mu  \right)d{{r}^{*}}d{{\sigma }^{2}}~.  
	\end{eqnarray}
Next, from equation (\ref{EK}), we can prove
	\begin{eqnarray}\label{lema1}
		&&\int_{{{r}^{*}}=r_{1}^{*}}^{{{r}^{*}}=r_{2}^{*}}{\int\limits_{{{S}^{2}}}{\left( \frac{1}{{{r}^{2}}(1-\mu )}{{\left| {{F }_{w\theta }} \right|}^{2}}+\frac{1}{(1-\mu ){{r}^{2}}{{\sin }^{2}}\theta }{{\left| {{F }_{w\varphi }} \right|}^{2}} \right.}}+\frac{1}{{{r}^{2}}(1-\mu )}{{\left| {{F }_{v\theta }} \right|}^{2}} \notag\\ 
		&& + \frac{1}{(1-\mu ){{r}^{2}}{{\sin }^{2}}\theta }{{\left| {{F }_{v\varphi }} \right|}^{2}}+\frac{1}{{{(1-\mu )}^{2}}}{{\left| {{F }_{vw}} \right|}^{2}}+\frac{(1-\mu )}{4{{r}^{4}}{{\sin }^{2}}\theta }{{\left| {{F }_{\varphi \theta }} \right|}^{2}} \notag\\ 
		&& +\left. \frac{{{\left| {{D}_{w}}\phi  \right|}^{2}}}{(1-\mu )}+\frac{{{\left| {{D}_{v}}\phi  \right|}^{2}}}{(1-\mu )}+\frac{1}{4{{r}^{2}}}{{\left| {{D}_{\theta }}\phi  \right|}^{2}}+\frac{1}{4{{r}^{2}}{{\sin }^{2}}\theta }{{\left| {{D}_{\varphi }}\phi  \right|}^{2}}+\frac{\left|P\right|}{2} \right)\left( 1-\mu  \right){{r}^{2}}d{{r}^{*}}d{{\sigma }^{2}} \notag\\ 
		&\lesssim& \frac{{{E}^{K}}\left( t \right)}{{{\min }_{w\in \left\{ t \right\}\cap \left\{ r_{1}^{*}\le {{r}^{*}}\le r_{2}^{*} \right\}}}{{w}^{2}}}+\frac{{{E}^{K}}\left( t \right)}{{{\min }_{v\in \left\{ t \right\}\cap \left\{ r_{1}^{*}\le {{r}^{*}}\le r_{2}^{*} \right\}}}{{v}^{2}}}~.	
	\end{eqnarray}
Thus, by utilizing \eqref{lema1}, it is possible to establish a bound for \eqref{redenergy}, which leads to the following inequality
	\begin{eqnarray}\label{Et}
		\tilde{E}^{\hat{t}}\left( {{t_i}} \right)&\lesssim& \frac{{{E}^{K}}\left( t_i \right)}{t_{i}^{2}}~. 
	\end{eqnarray}
In similiar way, we get
	\begin{eqnarray}\label{EtL}
		\tilde{E}_{{{\mathcal{L}}_{\Omega j}}}^{\hat{t}}\left( {{t_i}} \right)\lesssim \frac{1}{t_{i}^{2}}\sum\limits_{j=1}^{3}{E_{{{\mathcal{L}}_{{{\Omega }_{j}}}}}^{K}\left( {{t_i}} \right)}~. 
	\end{eqnarray}
Combining the result in (\ref{Et}) and (\ref{EtL}) to (\ref{es1}), one can get the estimate of $\mathcal{J}^K$ in the term of $E^K$.
\end{proof}

\begin{preposition}\label{prep3}
	Let $E^K$ given in equation (\ref{EK}). Then, for $t_{i+1}\le t_{i}+0.1t_i$ and $\left| {{r}^{*}}\left( {{r}_{0}} \right) \right|+\left| {{r}^{*}}\left( R \right) \right|\le 0.4{{t}_{i}}$, we have
	\begin{eqnarray}\label{EKproof}
		{{E}^{K}}\left( t \right)\lesssim E^{MH}~,
	\end{eqnarray}
	where 
	\begin{eqnarray}
  {{E}^{MH}}&\equiv&\sum\limits_{j=1}^{3}{\sum\limits_{k=1}^{3}{\sum\limits_{l=1}^{3}{\left( {{E}^{{\hat{t}}}}+E_{{{\mathcal{L}}_{\Omega j}}}^{{\hat{t}}}+E_{{{\mathcal{L}}_{\Omega j}}{{\mathcal{L}}_{\Omega k}}}^{{\hat{t}}}+E_{{{\mathcal{L}}_{\Omega j}}{{\mathcal{L}}_{\Omega k}}{{\mathcal{L}}_{\Omega l}}}^{{\hat{t}}} \right)}}}\left( t={{t}_{0}} \right) \notag\\ 
 & +&\sum\limits_{k=1}^{3}{\sum\limits_{l=1}^{3}{\left( {E}^K+E_{{{\mathcal{L}}_{\Omega k}}}^K+E_{{{\mathcal{L}}_{\Omega k}}{{\mathcal{L}}_{\Omega l}}}^K \right)}}\left( t={{t}_{0}} \right) ~. 
	\end{eqnarray}
\end{preposition}
\begin{proof}
For the vector field $K$ on region $t\in \left[ {{t}_{i}},{{t}_{i+1}} \right]$, it can be shown that
	\begin{eqnarray}\label{ek1}
		{{E}^{K}}\left( {{t}_{i+1}} \right)&\le& {{\mathcal{J}}^{K}}\left( {{t}_{0}}\le t\le {{t}_{i+1}} \right)+{{E}^{K}}\left( {{t}_{0}} \right) \notag\\ 
		&\le& {{t}_{i+1}}{{\mathcal{J}}^{G}}\left( {{t}_{0}}\le t\le {{t}_{i+1}} \right)\left( {{r}_{0}}\le r\le {{R}_{0}} \right)+{{E}^{K}}\left( {{t}_{0}} \right)~, 
	\end{eqnarray}
	and
	\begin{eqnarray}\label{ek2}  
		E_{{{\mathcal{L}}_{{{\Omega }_{j}}}}}^{K}\left( {{t}_{i+1}} \right)&\le& \mathcal{J}_{{{\mathcal{L}}_{{{\Omega }_{j}}}}}^{K}\left( {{t}_{0}}\le t\le {{t}_{i+1}} \right)+E_{{{\mathcal{L}}_{{{\Omega }_{j}}}}}^{K}\left( {{t}_{0}} \right) \notag\\ 
		&\le& {{t}_{i+1}}\mathcal{J}_{{{\mathcal{L}}_{{{\Omega }_{j}}}}}^{G}\left( {{t}_{0}}\le t\le {{t}_{i+1}} \right)\left( {{r}_{0}}\le r\le {{R}_{0}} \right)+E_{{{\mathcal{L}}_{{{\Omega }_{j}}}}}^{K}\left( {{t}_{0}} \right)~.  
	\end{eqnarray}
	Recall that from \eqref{asumeq},(\ref{ek1}), and (\ref{ek2}), one can obtain the inequality
	\begin{eqnarray}
		{{\mathcal{J}}^{G}}\left(t_i \le t\le t_{i+1}\right)&\le& \frac{1}{{{t}_{i}}}\sum\limits_{j=1}^{3}{\mathcal{J}_{,\text{ }{{\mathcal{L}}_{{{\Omega }_{j}}}}}^{G}\left( {{t}_{0}}\le t\le {{t}_{i+1}} \right)\left( {{r}_{0}}\le r\le {{R}_{0}} \right)}+\frac{1}{t_{i}^{2}}\sum\limits_{j=1}^{3}{E_{,\text{ }{{\mathcal{L}}_{{{\Omega }_{j}}}}}^{K}\left( {{t}_{0}} \right)} \notag\\ 
		&+&\frac{1}{{{t}_{i+1}}}\sum\limits_{j=1}^{3}{\mathcal{J}_{,\text{ }{{\mathcal{L}}_{{{\Omega }_{j}}}}}^{G}\left( {{t}_{0}}\le t\le {{t}_{i+1}} \right)\left( {{r}_{0}}\le r\le {{R}_{0}} \right)}+\frac{1}{t_{i+1}^{2}}\sum\limits_{j=1}^{3}{E_{,\text{ }{{\mathcal{L}}_{{{\Omega }_{j}}}}}^{K}\left( {{t}_{0}} \right)}~, \notag\\ 
	\end{eqnarray}
	where
	\begin{eqnarray}
		\mathcal{J}_{,\text{ }{{\mathcal{L}}_{{{\Omega }_{j}}}} }^{G}&=&{{\mathcal{J}}^{G}}+\mathcal{J}_{{{\mathcal{L}}_{{{\Omega }_{j}}}}}^{G}~,\notag\\
		E_{,\text{ }{{\mathcal{L}}_{{{\Omega }_{j}}}} }^{K}&=&{{E}^{K}}+E_{{{\mathcal{L}}_{{{\Omega }_{j}}}}}^{K}~.
	\end{eqnarray}
Furthermore, we can estimate $\mathcal{J}^K$ in the following form
	\begin{eqnarray}
		{{\mathcal{J}}^{K}}\left( {{t}_{i}}\le t\le {{t}_{i+1}} \right)&\lesssim& {{t}_{i+1}}{\mathcal{J}^{G}}\left( {{t}_{i}}\le t\le {{t}_{i+1}} \right)\left( {{r}_{0}}\le r\le {{R}_{0}} \right)\notag\\
  &\lesssim&\sum\limits_{j=1}^{3}{\mathcal{J}_{,{{\mathcal{L}}_{{{\Omega }_{j}}}} }^{G}\left( {{t}_{0}}\le t\le {{t}_{i+1}} \right)\left( {{r}_{0}}\le r\le {{R}_{0}} \right)}+\frac{1}{{{t}_{i+1}}}\sum\limits_{j=1}^{3}{E_{,{{\mathcal{L}}_{{{\Omega }_{j}}}} }^{K}\left( {{t}_{0}} \right)}~.\notag\\
	\end{eqnarray}
	Since ${{t}_{i+1}}=1.1{{t}_{i}}$ then ${{t}_{i+1}}={{\left( 1.1 \right)}^{i+1}}{{t}_{0}}$, implying $\sum\nolimits_{i}{\frac{1}{{{t}_{i+1}}}=}\sum\nolimits_{i}{\frac{1}{{{\left( 1.1 \right)}^{i+1}}{{t}_{0}}}}\lesssim 1$, such that
	\begin{eqnarray}
		{{\mathcal{J}}^{K}}\left( {{t}_{0}}\le t\le {{t}_{i+1}} \right)&=&\sum\limits_{i=0}^{i+1}{{{\mathcal{J}}^{K}}\left( {{t}_{i}}\le t\le {{t}_{i+1}} \right)} \notag\\ 
		&\lesssim& \left( i+1 \right)\text{ }\sum\limits_{j=1}^{3}{\mathcal{J}_{,{{\mathcal{L}}_{{{\Omega }_{j}}}}}^{G}\left( {{t}_{0}}\le t\le {{t}_{i+1}} \right)\left( {{r}_{0}}\le r\le {{R}_{0}} \right)}+\sum\limits_{j=1}^{3}{E_{,{{\mathcal{L}}_{{{\Omega }_{j}}}}}^{K}\left( {{t}_{0}} \right)}~, \notag\\ 
	\end{eqnarray}
we attain
	\begin{eqnarray}
		{{E}^{K}}\left( {{t}_{i+1}} \right)\lesssim \left( i+1 \right)\text{ }\sum\limits_{j=1}^{3}{\mathcal{J}_{,{{\mathcal{L}}_{{{\Omega }_{j}}}}}^{G}\left( {{t}_{0}}\le t\le {{t}_{i+1}} \right)\left( {{r}_{0}}\le r\le {{R}_{0}} \right)}+\sum\limits_{j=1}^{3}{E_{,{{\mathcal{L}}_{{{\Omega }_{j}}}}\Psi }^{K}\left( {{t}_{0}} \right)}~,
	\end{eqnarray}
	\begin{eqnarray}
		E_{{{\mathcal{L}}_{{{\Omega }_{j}}}}}^{K}\left( {{t}_{i+1}} \right)&\lesssim& \left( i+1 \right)\text{ }\sum\limits_{j=1}^{3}{\sum\limits_{l=1}^{3}{\mathcal{J}_{{{\mathcal{L}}_{{{\Omega }_{l}}}},{{\mathcal{L}}_{{{\Omega }_{l}}}}{{\mathcal{L}}_{{{\Omega }_{j}}}}}^{G}}\left( {{t}_{0}}\le t\le {{t}_{i+1}} \right)\left( {{r}_{0}}\le r\le {{R}_{0}} \right)}+\sum\limits_{j=1}^{3}{\sum\limits_{l=1}^{3}{E_{{{\mathcal{L}}_{{{\Omega }_{l}}}},{{\mathcal{L}}_{{{\Omega }_{l}}}}{{\mathcal{L}}_{{{\Omega }_{j}}}}}^{K}\left( {{t}_{0}} \right)}}~,\notag\\  
	\end{eqnarray}
	\begin{eqnarray}
		{{\mathcal{J}}^{K}}\left( {{t}_{i}}\le t\le {{t}_{i+1}} \right)\left( -\infty \le {{r}^{*}}\le \infty  \right)&\lesssim& {{t}_{i+1}}{{\mathcal{J}}^{G}}\left( {{t}_{i}}\le t\le {{t}_{i+1}} \right)\left( {{r}_{0}}\le r\le {{R}_{0}} \right) \notag\\ 
		& \lesssim&\frac{i+1}{{{t}_{i+1}}}\sum\limits_{l=1}^{3}{\sum\limits_{j=1}^{3}{\mathcal{J}_{,{{\mathcal{L}}_{{{\Omega }_{l}}}},{{\mathcal{L}}_{{{\Omega }_{l}}}}{{\mathcal{L}}_{{{\Omega }_{j}}}}}^{G}\left( {{t}_{0}}\le t\le {{t}_{i+1}} \right)\left( {{r}_{0}}\le r\le {{R}_{0}} \right)}} \notag\\ 
		&+&\frac{1}{{{t}_{i+1}}}\sum\limits_{l=1}^{3}{\sum\limits_{j=1}^{3}{E_{,{{\mathcal{L}}_{{{\Omega }_{l}}}},{{\mathcal{L}}_{{{\Omega }_{l}}}}{{\mathcal{L}}_{{{\Omega }_{j}}}}}^{K}\left( {{t}_{0}} \right)}}  ~,
	\end{eqnarray}
	which leads to 
	\begin{eqnarray}\label{ekti1}
		{{E}^{K}}\left( {{t}_{i+1}} \right)&\lesssim& {{\mathcal{J}}^{K}}\left( {{t}_{0}}\le t\le {{t}_{i+1}} \right)\left( -\infty \le {{r}^{*}}\le \infty  \right)+{{E}^{K}}\left( {{t}_{0}} \right) \notag\\ 
		& \lesssim& \sum\limits_{l=1}^{3}{\sum\limits_{j=1}^{3}{\left[ \mathcal{J}_{,{{\mathcal{L}}_{{{\Omega }_{l}}}},{{\mathcal{L}}_{{{\Omega }_{l}}}}{{\mathcal{L}}_{{{\Omega }_{j}}}}}^{G}\left( {{t}_{0}}\le t\le {{t}_{i+1}} \right)\left( {{r}_{0}}\le r\le {{R}_{0}} \right)+E_{,{{\mathcal{L}}_{{{\Omega }_{l}}}},{{\mathcal{L}}_{{{\Omega }_{l}}}}{{\mathcal{L}}_{{{\Omega }_{j}}}}}^{K}\left( {{t}_{0}} \right) \right]}}~. \notag\\ 
	\end{eqnarray}
From (\ref{asumeq}), we have
	\begin{eqnarray}\label{jgkecil}
		{{\mathcal{J}}^{G}}\left( {{t}_{0}}\le t\le {{t}_{i+1}} \right)\left( {{r}_{0}}\le r\le {{R}_{0}} \right)\lesssim \sum\limits_{j=1}^{3}{E_{,{{\mathcal{L}}_{{{\Omega }_{j}}}}}^{\left( \hat{t} \right)}}~.
	\end{eqnarray}
	Thus, substituting \eqref{jgkecil} into \eqref{ekti1}, one can obtain,
	\begin{eqnarray}\label{esEK}
		{{E}^{K}}\left( {{t}_{i+1}} \right)&\lesssim&\sum\limits_{j=1}^{3}{\sum\limits_{k=1}^{3}{\sum\limits_{l=1}^{3}{\left( {{E}^{{\hat{t}}}}+E_{{{\mathcal{L}}_{\Omega j}}}^{{\hat{t}}}+E_{{{\mathcal{L}}_{\Omega j}}{{\mathcal{L}}_{\Omega k}}}^{{\hat{t}}}+E_{{{\mathcal{L}}_{\Omega j}}{{\mathcal{L}}_{\Omega k}}{{\mathcal{L}}_{\Omega l}}}^{{\hat{t}}} \right)}}}\left( t={{t}_{0}} \right)\notag\\ 
 & +&\sum\limits_{k=1}^{3}{\sum\limits_{l=1}^{3}{\left( {E}^K+E_{{{\mathcal{L}}_{\Omega k}}}^K+E_{{{\mathcal{L}}_{\Omega k}}{{\mathcal{L}}_{\Omega l}}}^K \right)}}\left( t={{t}_{0}} \right)~.
	\end{eqnarray}
Our proof worked not only on $t_{i+1}=(1.1)^{i+1}t_0$, but also for all $t>t_0$ since for $t_0>1$ and $1<a<2$ there is always exist $i$ and $a$ such that $t=a^i t_0$. So, for all $t$, we have \eqref{EKproof} 
\end{proof}

\section{Proof of Theorem 1}
\label{sec:decayfaraway}
In this section, we prove the decay estimate of MH system outside photon sphere region $r>3m$.

\subsection{Decay for the gauge field}
Let us consider the Sobolev inequality for the electromagnetic field in the region $w\ge 1,r\ge R$, where $R$ is fixed, 
\begin{eqnarray}\label{emineq}
{{r}^{2}}{{\left| F  \right|}^{2}}\le \int_{{{S}^{2}}}{{{r}^{2}}{{\left| F  \right|}^{2}}d{{\sigma }^{2}}}+\int_{{{S}^{2}}}{{{r}^{2}}{{\left| \slashed{\mathcal{L}}\left| F  \right| \right|}^{2}}d{{\sigma }^{2}}}+\int_{{{S}^{2}}}{{{r}^{2}}{{\left| \slashed{\mathcal{L}}\slashed{\mathcal{L}}\left| F  \right| \right|}^{2}}d{{\sigma }^{2}}}~.
\end{eqnarray}
Let ${{r}_{F}}$ be a value of $r$ such that $R\le {{r}_{F }}\le R+1$, then	
\begin{eqnarray}
\int_{{{S}^{2}}}{{{r}^{2}}{{\left| F  \right|}^{2}}\left( t,r,\theta,\varphi \right)d{{\sigma }^{2}}}&\lesssim& \int_{{{S}^{2}}}{{{r_F}^{2}}{{\left| F \right|}^{2}}\left( t,{{r}_{F }},\theta,\varphi \right)d{{\sigma }^{2}}}+\int_{{{S}^{2}}}{\int_{{{{\bar{r}}}^{*}}=r_{F }^{*}}^{{{{\bar{r}}}^{*}}={{r}^{*}}}{{{\nabla }_{{{r}^{*}}}}\left[ {{r}^{2}}{{\left| F  \right|}^{2}} \right]}\left( t,r,\theta,\varphi \right)d{{{\bar{r}}}^{*}}d{{\sigma }^{2}}} \notag\\ 
&\lesssim& \int_{{{S}^{2}}}{r_{F }^{2}{{\left| F \right|}^{2}}\left( t,{{r}_{F }},\theta,\varphi \right)d{{\sigma }^{2}}}+\int_{{{S}^{2}}}{\int_{{{{\bar{r}}}^{*}}=r_{F }^{*}}^{{{{\bar{r}}}^{*}}={{r}^{*}}}{2r{{\left| F \right|}^{2}}}\left( t,r,\theta,\varphi \right)\left( 1-\mu  \right)d{{{\bar{r}}}^{*}}d{{\sigma }^{2}}} \notag\\ 
&+&\int_{{{S}^{2}}}{\int_{{{{\bar{r}}}^{*}}=r_{F }^{*}}^{{{{\bar{r}}}^{*}}={{r}^{*}}}{{{r}^{2}}{{\nabla }_{{{r}^{*}}}}{{\left| F  \right|}^{2}}}\left( t,r,\theta,\varphi \right)d{{{\bar{r}}}^{*}}d{{\sigma }^{2}}}~.  
\end{eqnarray}	
By using (\ref{lema1}), we can estimate 
\begin{eqnarray}
	\int_{{{S}^{2}}}{r_{F }^{2}{{\left| F  \right|}^{2}}\left( t,{{r}_{F }},\theta,\varphi \right)d{{\sigma }^{2}}}\le \frac{{{E}^{K}}\left( t \right)}{{{t}^{2}}}~. 
\end{eqnarray}	
In the same way, it follows that
\begin{eqnarray}
\int\limits_{{{S}^{2}}}r_{{{\slashed{\mathcal{L}}}_{F}}}^{2}{{\left| \slashed{\mathcal{L}}F \right|}^{2}}\left( t,{{r}_{{{\slashed{\mathcal{L}}}_{F}}}},\theta ,\varphi  \right)d{{\sigma }^{2}}&=&\sum\limits_{j=1}^{3}{\frac{1}{r_{{{\slashed{\mathcal{L}}}_{F}}}^{2}}\int_{{{S}^{2}}}{r_{{{\slashed{\mathcal{L}}}_{F}}}^{2}{{\left| {{\mathcal{L}}_{{{\Omega }_{j}}}}F \right|}^{2}}\left( t,{{r}_{{{\slashed{\mathcal{L}}}_{F}}}},\theta ,\varphi  \right)d{{\sigma }^{2}}}} \notag\\ 
&\lesssim&\sum\limits_{j=1}^{3}{\frac{1}{{{R}^{2}}}}\int_{{{S}^{2}}}{r_{{{\slashed{\mathcal{L}}}_{F}}}^{2}{{\left| {{\mathcal{L}}_{{{\Omega }_{j}}}}F \right|}^{2}}\left( t,{{r}_{{{\slashed{\mathcal{L}}}_{F}}}},\theta ,\varphi  \right)d{{\sigma }^{2}}} \notag\\ 
&\lesssim& \sum\limits_{j=1}^{3}{\frac{E_{{{\mathcal{L}}_{{{\Omega }_{j}}}}}^{K}\left( t \right)}{{{t}^{2}}}}~,  
\end{eqnarray}
and
\begin{eqnarray}
	\int\limits_{{{S}^{2}}}{r_{\slashed{\mathcal{L}}{{\slashed{\mathcal{L}}}_{F}}}^{2}{{\left| \slashed{\mathcal{L}}\slashed{\mathcal{L}}F \right|}^{2}}\left( t,{{r}_{{{\slashed{\mathcal{L}}}_{F}}}},\theta ,\varphi  \right)d{{\sigma }^{2}}}\le \sum\limits_{i=1}^{3}{\sum\limits_{j=1}^{3}{\frac{E_{{{\mathcal{L}}_{{{\Omega }_{i}}}}{{\mathcal{L}}_{{{\Omega }_{j}}}}}^{K}\left( t \right)}{{{t}^{2}}}}}~.
\end{eqnarray}
On the other hand,
\begin{eqnarray}
	\int_{{{r}^{*}}=r_{F }^{*}}^{{{r}^{*}}={{r}^{*}}}{\int\limits_{{{S}^{2}}}{\bar{r}{{\left| F  \right|}^{2}}}}\left( 1-\mu  \right)\left( t,\bar{r},\theta,\varphi \right)d{{{\bar{r}}}^{*}}d{{\sigma }^{2}}&\lesssim&\int_{{{r}^{*}}=r_{F }^{*}}^{{{r}^{*}}={{r}^{*}}}{\frac{{\bar{r}}}{R}\int\limits_{{{S}^{2}}}{\bar{r}{{\left| F  \right|}^{2}}}}\left( 1-\mu  \right)\left( t,\bar{r},\theta,\varphi \right)d{{{\bar{r}}}^{*}}d{{\sigma }^{2}} \notag\\ 
	&\lesssim& \frac{{{E}^{K}}\left( t \right)}{{{t}^{2}}}+\frac{{{E}^{K}}\left( t \right)}{{{w}^{2}}}~.  
\end{eqnarray}
By the same method as above, we get
\begin{eqnarray}	\int_{{{r}^{*}}=r_{{{\slashed{\mathcal{L}}}_{F}}}^{*}}^{{{r}^{*}}={{r}^{*}}}{\int\limits_{{{S}^{2}}}{\bar{r}{{\left| \slashed{\mathcal{L}}F \right|}^{2}}}}\left( 1-\mu  \right)\left( t,\bar{r},\theta ,\varphi  \right)d{{\bar{r}}^{*}}d{{\sigma }^{2}}\lesssim \sum\limits_{j=1}^{3}{\frac{E_{{{\mathcal{L}}_{{{\Omega }_{j}}}}}^{K}\left( t \right)}{{{t}^{2}}}}+\sum\limits_{j=1}^{3}{\frac{E_{{{\mathcal{L}}_{{{\Omega }_{j}}}}}^{K}\left( t \right)}{{{w}^{2}}}}~,
\end{eqnarray}
\begin{eqnarray}
	\int_{{{r}^{*}}=r_{\slashed{\mathcal{L}}{{\slashed{\mathcal{L}}}_{F}}}^{*}}^{{{r}^{*}}={{r}^{*}}}{\int\limits_{{{S}^{2}}}{\bar{r}{{\left| \slashed{\mathcal{L}}\slashed{\mathcal{L}}F \right|}^{2}}}}\left( 1-\mu  \right)\left( t,\bar{r},\theta ,\varphi  \right)d{{\bar{r}}^{*}}d{{\sigma }^{2}}\lesssim \sum\limits_{i=1}^{3}{\sum\limits_{j=1}^{3}{\frac{E_{{{\mathcal{L}}_{{{\Omega }_{i}}}}{{\mathcal{L}}_{{{\Omega }_{j}}}}}^{K}\left( t \right)}{{{t}^{2}}}}}+\sum\limits_{i=1}^{3}{\sum\limits_{j=1}^{3}{\frac{E_{{{\mathcal{L}}_{{{\Omega }_{i}}}}{{\mathcal{L}}_{{{\Omega }_{j}}}}}^{K}\left( t \right)}{{{w}^{2}}}}}~.\notag\\
\end{eqnarray}
Next, we want to estimate the term
\begin{eqnarray}
	\int_{{{S}^{2}}}{\int_{{{{\bar{r}}}^{*}}=r_{F }^{*}}^{{{{\bar{r}}}^{*}}={{r}^{*}}}{{{r}^{2}}{{\nabla }_{{{r}^{*}}}}{{\left| F  \right|}^{2}}}\left( t,r,\theta,\varphi \right)d{{{\bar{r}}}^{*}}d{{\sigma }^{2}}} &\lesssim& {{\left[ \int_{{{S}^{2}}}{\int_{{{{\bar{r}}}^{*}}=r_{F }^{*}}^{{{{\bar{r}}}^{*}}={{r}^{*}}}{{{r}^{2}}{{\left| {{\nabla }_{{{{\hat{r}}}^{*}}}}F  \right|}^{2}}}\left( t,r,\theta,\varphi \right)\left( 1-\mu  \right)d{{{\bar{r}}}^{*}}d{{\sigma }^{2}}} \right]}^{1/2}} \notag\\ 
	&\cdot& {{\left[ \int_{{{S}^{2}}}{\int_{{{{\bar{r}}}^{*}}=r_{F }^{*}}^{{{{\bar{r}}}^{*}}={{r}^{*}}}{{{r}^{2}}{{\left| F  \right|}^{2}}}\left( t,r,w \right)d{{{\bar{r}}}^{*}}d{{\sigma }^{2}}} \right]}^{1/2}}~. \notag\\ 
\end{eqnarray}	
By using the Bianchi identity, one can show that
\begin{eqnarray}
{{\left| {{\nabla }_{{{{\hat{r}}}^{*}}}}F \right|}^{2}}\left( 1-\mu  \right)\lesssim {{\left| F \right|}^{2}}+{{\left| {{\mathcal{L}}_{t}}F \right|}^{2}}+\sum\limits_{j=1}^{3}{{{\left| {{\mathcal{L}}_{{{\Omega }_{j}}}}F \right|}^{2}}}~.
\end{eqnarray}	
Hence, we have
\begin{eqnarray}
	\int_{{{S}^{2}}}{\int_{{{{\bar{r}}}^{*}}=r_{F }^{*}}^{{{{\bar{r}}}^{*}}={{r}^{*}}}{{{r}^{2}}{{\nabla }_{{{r}^{*}}}}{{\left| F  \right|}^{2}}}\left( t,r,\theta,\varphi \right)d{{{\bar{r}}}^{*}}d{{\sigma }^{2}}} &\lesssim&\frac{{{E}^{K}}\left( t \right)+E_{{{\mathcal{L}}_{t}}}^{K}\left( t \right)+\sum\limits_{j=1}^{3}{E_{{{\mathcal{L}}_{{{\Omega }_{j}}}}}^{K}\left( t \right)}}{{{t}^{2}}}\notag\\
	&+&\frac{{{E}^{K}}\left( t \right)+E_{{{\mathcal{L}}_{t}}}^{K}\left( t \right)+\sum\limits_{j=1}^{3}{E_{{{\mathcal{L}}_{{{\Omega }_{j}}}}}^{K}\left( t \right)}}{{{w}^{2}}}~.
\end{eqnarray}	
In a similar way, it can be demonstrated that
\begin{eqnarray}
	\int_{{{S}^{2}}}{\int_{{{{\bar{r}}}^{*}}=r_{\slashed{\mathcal{L}}F }^{*}}^{{{{\bar{r}}}^{*}}={{r}^{*}}}{{{r}^{2}}{{\nabla }_{{{r}^{*}}}}{{\left| \slashed{\mathcal{L}}F  \right|}^{2}}}\left( t,r,\theta,\varphi \right)d{{{\bar{r}}}^{*}}d{{\sigma }^{2}}} &\lesssim&\sum\limits_{i=1}^{3}{\sum\limits_{j=1}^{3}{\frac{E_{{{\mathcal{L}}_{F}},{{\mathcal{L}}_{\Omega t}}{{\mathcal{L}}_{{{\Omega }_{i}}}},{{\mathcal{L}}_{{{\Omega }_{i}}}}{{\mathcal{L}}_{{{\Omega }_{j}}}}}^{K}\left( t \right)}{{{t}^{2}}}}}\notag\\
	&+&\sum\limits_{i=1}^{3}{\sum\limits_{j=1}^{3}{\frac{E_{{{\mathcal{L}}_{F}},{{\mathcal{L}}_{\Omega t}}{{\mathcal{L}}_{{{\Omega }_{i}}}},{{\mathcal{L}}_{{{\Omega }_{i}}}}{{\mathcal{L}}_{{{\Omega }_{j}}}}}^{K}\left( t \right)}{{{w}^{2}}}}}~,
\end{eqnarray}
and
\begin{eqnarray}
	\int_{{{S}^{2}}}{\int_{{{{\bar{r}}}^{*}}=r_{\slashed{\mathcal{L}}\slashed{\mathcal{L}}F }^{*}}^{{{{\bar{r}}}^{*}}={{r}^{*}}}{{{r}^{2}}{{\nabla }_{{{r}^{*}}}}{{\left| \slashed{\mathcal{L}}\slashed{\mathcal{L}}F  \right|}^{2}}}\left( t,r,\theta,\varphi \right)d{{{\bar{r}}}^{*}}d{{\sigma }^{2}}} &\lesssim&\frac{\sum\limits_{i=1}^{3}{\sum\limits_{j=1}^{3}{\sum\limits_{l=1}^{3}{E_{{{\mathcal{L}}_{{{\Omega }_{i}}}}{{\mathcal{L}}_{{{\Omega }_{j}}}},{{\mathcal{L}}_{{{\Omega }_{t}}}}{{\mathcal{L}}_{{{\Omega }_{i}}}}{{\mathcal{L}}_{{{\Omega }_{j}}}},{{\mathcal{L}}_{\Omega l}}{{\mathcal{L}}_{{{\Omega }_{i}}}}{{\mathcal{L}}_{{{\Omega }_{j}}}}}^{K}\left( t \right)}}}}{{{t}^{2}}}\notag\\
	&+&\frac{\sum\limits_{i=1}^{3}{\sum\limits_{j=1}^{3}{\sum\limits_{l=1}^{3}{E_{{{\mathcal{L}}_{{{\Omega }_{i}}}}{{\mathcal{L}}_{{{\Omega }_{j}}}},{{\mathcal{L}}_{{{\Omega }_{t}}}}{{\mathcal{L}}_{{{\Omega }_{i}}}}{{\mathcal{L}}_{{{\Omega }_{j}}}},{{\mathcal{L}}_{\Omega l}}{{\mathcal{L}}_{{{\Omega }_{i}}}}{{\mathcal{L}}_{{{\Omega }_{j}}}}}^{K}\left( t \right)}}}}{{{w}^{2}}}~.\notag\\
\end{eqnarray}

Finally, we can estimate \eqref{emineq} as follows
\begin{eqnarray}
	{{r}^{2}}{{\left| F \right|}^{2}}\left( t,r,\theta ,\varphi  \right)&\lesssim&\frac{\sum\limits_{j=1}^{3}{E_{,{{\mathcal{L}}_{t}},{{\mathcal{L}}_{{{\Omega }_{j}}}}}^{K}\left( t \right)}}{{{t}^{2}}}+\frac{\sum\limits_{j=1}^{3}{E_{,{{\mathcal{L}}_{t}},{{\mathcal{L}}_{{{\Omega }_{j}}}}}^{K}\left( t \right)}}{{{w}^{2}}}\notag\\
	&+&\frac{\sum\limits_{i=1}^{3}{\sum\limits_{j=1}^{3}{E_{,{{\mathcal{L}}_{t}}{{\mathcal{L}}_{{{\Omega }_{i}}}},{{\mathcal{L}}_{{{\Omega }_{i}}}}{{\mathcal{L}}_{{{\Omega }_{j}}}}}^{K}\left( t \right)}}}{{{t}^{2}}}+\frac{\sum\limits_{i=1}^{3}{\sum\limits_{j=1}^{3}{E_{,{{\mathcal{L}}_{t}}{{\mathcal{L}}_{{{\Omega }_{i}}}},{{\mathcal{L}}_{{{\Omega }_{i}}}}{{\mathcal{L}}_{{{\Omega }_{j}}}}}^{K}\left( t \right)}}}{{{w}^{2}}}\notag\\
	&+&\frac{\sum\limits_{i=1}^{3}{\sum\limits_{j=1}^{3}{\sum\limits_{l=1}^{3}{E_{{{\mathcal{L}}_{{{\Omega }_{i}}}}{{\mathcal{L}}_{{{\Omega }_{j}}}},{{\mathcal{L}}_{{{\Omega }_{t}}}}{{\mathcal{L}}_{{{\Omega }_{i}}}}{{\mathcal{L}}_{{{\Omega }_{j}}}},{{\mathcal{L}}_{\Omega l}}{{\mathcal{L}}_{{{\Omega }_{i}}}}{{\mathcal{L}}_{{{\Omega }_{j}}}}}^{K}\left( t \right)}}}}{{{t}^{2}}}\notag\\
	&+&\frac{\sum\limits_{i=1}^{3}{\sum\limits_{j=1}^{3}{\sum\limits_{l=1}^{3}{E_{{{\mathcal{L}}_{{{\Omega }_{i}}}}{{\mathcal{L}}_{{{\Omega }_{j}}}},{{\mathcal{L}}_{{{\Omega }_{t}}}}{{\mathcal{L}}_{{{\Omega }_{i}}}}{{\mathcal{L}}_{{{\Omega }_{j}}}},{{\mathcal{L}}_{\Omega l}}{{\mathcal{L}}_{{{\Omega }_{i}}}}{{\mathcal{L}}_{{{\Omega }_{j}}}}}^{K}\left( t \right)}}}}{{{w}^{2}}}~.\notag\\
\end{eqnarray}	
Therefore, we have	
\begin{eqnarray}\label{esF}
	\left| {{F}_{\hat{\mu }\hat{\nu }}} \right|\left( w,v,\theta ,\varphi  \right)&\lesssim& \frac{{{\left[ E_{,{{\mathcal{L}}_{t}},{{\mathcal{L}}_{{{\Omega }_{i}}}},{{\mathcal{L}}_{t}}{{\mathcal{L}}_{{{\Omega }_{i}}}},{{\mathcal{L}}_{{{\Omega }_{i}}}}{{\mathcal{L}}_{{{\Omega }_{j}}}},{{\mathcal{L}}_{t}}{{\mathcal{L}}_{{{\Omega }_{i}}}}{{\mathcal{L}}_{{{\Omega }_{j}}}},{{\mathcal{L}}_{l}}{{\mathcal{L}}_{{{\Omega }_{i}}}}{{\mathcal{L}}_{{{\Omega }_{j}}}}}^{K} \right]}^{1/2}}}{rt}\notag\\
	&+&\frac{{{\left[ E_{,{{\mathcal{L}}_{t}},{{\mathcal{L}}_{{{\Omega }_{i}}}},{{\mathcal{L}}_{t}}{{\mathcal{L}}_{{{\Omega }_{i}}}},{{\mathcal{L}}_{{{\Omega }_{i}}}}{{\mathcal{L}}_{{{\Omega }_{j}}}},{{\mathcal{L}}_{t}}{{\mathcal{L}}_{{{\Omega }_{i}}}}{{\mathcal{L}}_{{{\Omega }_{j}}}},{{\mathcal{L}}_{l}}{{\mathcal{L}}_{{{\Omega }_{i}}}}{{\mathcal{L}}_{{{\Omega }_{j}}}}}^{K} \right]}^{1/2}}}{rw}~.
\end{eqnarray}
Consider the region where $t\ge 1$ such that we have $r+t\le rt$ for $R$ fixed. Since $w\ge 1$, consequently $v={{r}^{*}}+t\lesssim 2r+2t\lesssim rt$ and $v\lesssim r+t\lesssim r+w\le rw$. Then 
\begin{eqnarray}\label{region1}
	\left| {{F}_{\hat{\mu }\hat{\nu }}} \right|\left( w,v,\theta ,\varphi  \right)\le \frac{{{\left[ E_{,{{\mathcal{L}}_{t}},{{\mathcal{L}}_{{{\Omega }_{i}}}},{{\mathcal{L}}_{t}}{{\mathcal{L}}_{{{\Omega }_{i}}}},{{\mathcal{L}}_{{{\Omega }_{i}}}}{{\mathcal{L}}_{{{\Omega }_{j}}}},{{\mathcal{L}}_{t}}{{\mathcal{L}}_{{{\Omega }_{i}}}}{{\mathcal{L}}_{{{\Omega }_{j}}}},{{\mathcal{L}}_{l}}{{\mathcal{L}}_{{{\Omega }_{i}}}}{{\mathcal{L}}_{{{\Omega }_{j}}}}}^{K} \right]}^{1/2}}}{v}~.
\end{eqnarray}
For the region $t\le 1$, $w\ge 1,r\ge R$ is a bounded compact region, so in this region one can obtain
\begin{eqnarray}\label{comregion}
	\left| F  \right|\lesssim E_{{{\mathcal{L}}_{{{\Omega }_{j}}}},{{\mathcal{L}}_{{{\Omega }_{i}}}}{{\mathcal{L}}_{{{\Omega }_{j}}}}}^{\hat{t}}\left( t={{t}_{0}} \right)~.
\end{eqnarray}
Recall that from (\ref{EKproof}), we have
\begin{eqnarray}\label{ineqEK}
	{{E}^{K}}\lesssim \sum\limits_{j=0}^{3}{E_{{{r}^{j}}{{\left( \slashed{\mathcal{L}} \right)}^{j}}}^{\hat{t}}\left( {{t}_{0}} \right)}+\sum\limits_{j=0}^{2}{E_{{{r}^{j}}{{\left( \slashed{\mathcal{L}} \right)}^{j}}}^{K}\left( {{t}_{0}} \right)}
\end{eqnarray}
Finally, combining (\ref{region1}),(\ref{comregion}), and (\ref{ineqEK}), we get
\begin{eqnarray}
		E_{,{{\mathcal{L}}_{t}},{{\mathcal{L}}_{{{\Omega }_{i}}}},{{\mathcal{L}}_{t}}{{\mathcal{L}}_{{{\Omega }_{i}}}},{{\mathcal{L}}_{{{\Omega }_{i}}}}{{\mathcal{L}}_{{{\Omega }_{j}}}},{{\mathcal{L}}_{t}}{{\mathcal{L}}_{{{\Omega }_{i}}}}{{\mathcal{L}}_{{{\Omega }_{j}}}},{{\mathcal{L}}_{l}}{{\mathcal{L}}_{{{\Omega }_{i}}}}{{\mathcal{L}}_{{{\Omega }_{j}}}}}^{K}\lesssim\left(\hat{E}^{MH}\right)^{2},
\end{eqnarray}
where 
\begin{eqnarray}
	\hat{E}^{MH}&=&\sum\limits_{i=0}^{1}{\sum\limits_{j=0}^{5}{E_{{{r}^{j}}{{\left( \slashed{\mathcal{L}} \right)}^{j}}{{\left( {{\mathcal{L}}_{t}} \right)}^{i}}}^{\hat{t}}\left( {{t}_{0}} \right)}}+\sum\limits_{i=0}^{1}{\sum\limits_{j=0}^{4}{E_{{{r}^{j}}{{\left( \slashed{\mathcal{L}} \right)}^{j}}{{\left( {{\mathcal{L}}_{t}} \right)}^{i}}}^{K}\left( {{t}_{0}} \right)}} \notag\\ 
	&+&E_{{{r}^{6}}{{\left( \slashed{\mathcal{L}} \right)}^{6}}}^{\hat{t}}\left( {{t}_{0}} \right)+E_{{{r}^{5}}{{\left( \slashed{\mathcal{L}} \right)}^{5}}}^{ K }\left( {{t}_{0}} \right).
\end{eqnarray}
Thus, combining the result in \eqref{region1} and \eqref{comregion}, we attain \eqref{decayF1}.

Now, let us focus on the region $w \le -1, r \geq R, \left|t\right| \geq 1$, we have $r^* =\frac{v-w}{2}$. Thus, for $w \le -1$, we have, $r^*\geq v$, and hence $r \geq v$. Therefore, $\frac{1}{r}\lesssim\frac{1}{v}$. Since in this region we have $\left|w\right| \geq 1$, and $\left|t\right| \geq 1$, similarly, by employing some computation, we obtain the estimates as in (\ref{esF}) that leads to \eqref{decayF1}.

We also have for fixed $R$ and $\left|t\right| \geq 1$, $\left|w\right|=\left|t-r^*\right|\lesssim\left|r^*\right|+\left|t\right|\lesssim\left|rt\right|$. Thus, in this region, by same method, one can get \eqref{decayF2}. 

Let us consider the region $w \le -1, r \geq R,$ and $-1 \le t \le 1$. Let $\tilde{t}=t-2$, when $-1 \le t \le 1$, we have $-3 \le \tilde{t} \le -1$. Suppose $\tilde{w}=\tilde{t}-r^*=t-r^*-2$, for $w\le-1$ we get $\tilde{w}\le-3$. Thus, the region is in the new system of coordinates included in the region $\tilde{w} \le -1, r \geq R, \tilde{t}\le-1$. $\frac{\partial}{\partial t}$ is a Killing vector field, therefore, the time translation will keep the metric invariant in the new system of coordinates $\left\{\tilde{t}, r, \theta, \varphi\right\}$. Consequently, we will have the same results proven previously.

Next, we consider the region $-1\le w\le1$, $r\geq R$. Let $\overline{r}=r+2$, then when $-1 \le w \le 1$, we have $-3 \le \overline{w} \le -1$ and when $r^*\geq R^*$, $\overline{r}\geq R^*+2\geq R$. Thus the region $-1 \le w \le 1, r \geq R$ is included in the new system of coordinates in the region $\overline{w} \le -1, \overline{r}\geq R^*$. As a result, we will see the same outcomes as before.

To estimate the potential $A_\mu$ we need to use again the Sobolev inequality
\begin{eqnarray}\label{estiA}
{{r}^{2}}{{\left| A \right|}^{2}}\le C\int_{{{S}^{2}}}{{{r}^{2}}{{\left| A \right|}^{2}}d{{\sigma }^{2}}}+C\int_{{{S}^{2}}}{{{r}^{2}}{{\left| \slashed{\nabla }A \right|}^{2}}d{{\sigma }^{2}}}+C\int_{{{S}^{2}}}{{{r}^{2}}{{\left| \slashed{\nabla }\slashed{\nabla }A \right|}^{2}}d{{\sigma }^{2}}}~.
\end{eqnarray}
Then, we can estimate the first term of \eqref{estiA} as
\begin{eqnarray}\label{firsttermestiA}
  \int_{{{S}^{2}}}{{{r}^{2}}{{\left| A \right|}^{2}}\left( t,r,\theta ,\varphi  \right)d{{\sigma }^{2}}}&\le& \int_{{{S}^{2}}}{{{r_A}^{2}}{{\left| A \right|}^{2}}\left( t,{{r}_{A}},\theta ,\varphi  \right)d{{\sigma }^{2}}}+\int_{{{S}^{2}}}{\int_{{{{\bar{r}}}^{*}}=r_{A}^{*}}^{{{{\bar{r}}}^{*}}={{r}^{*}}}{{{\nabla }_{{{r}^{*}}}}\left[ {{r}^{2}}{{A}^{2}} \right]}\left( t,r,\theta ,\varphi  \right)d{{{\bar{r}}}^{*}}d{{\sigma }^{2}}} \notag\\ 
 &\le& \int_{{{S}^{2}}}{r_{A}^{2}{{\left| A \right|}^{2}}\left( t,{{r}_{A}},\theta ,\varphi  \right)d{{\sigma }^{2}}}+\int_{{{S}^{2}}}{\int_{{{{\bar{r}}}^{*}}=r_{A}^{*}}^{{{{\bar{r}}}^{*}}={{r}^{*}}}{2r{{A}^{2}}}\left( 1-\mu  \right)\left( t,r,\theta ,\varphi  \right)d{{{\bar{r}}}^{*}}d{{\sigma }^{2}}} \notag\\ 
 &+&\int_{{{S}^{2}}}{\int_{{{{\bar{r}}}^{*}}=r_{A}^{*}}^{{{{\bar{r}}}^{*}}={{r}^{*}}}{{{r}^{2}}{{\nabla }_{{{r}^{*}}}}\left( {{A}^{2}} \right)}\left( t,r,\theta ,\varphi  \right)d{{{\bar{r}}}^{*}}d{{\sigma }^{2}}} ~. 
\end{eqnarray}
By using Holder inequality, we can show that
\begin{eqnarray}
  \frac{\partial }{\partial t}\left\| A \right\|_{{{L}^{2}}}^{2}&=&\int_{{{r}^{*}}=r_{1}^{*}}^{{{r}^{*}}=r_{2}^{*}}{\int\limits_{{{S}^{2}}}{\frac{\partial }{\partial t}\left( {{A}^{2}} \right)}}\text{ }\left( 1-\mu  \right)^{1/2}{{r}^{2}}d{{{\bar{r}}}^{*}}d{{\sigma }^{2}} \notag\\ 
 {{\left\| A \right\|}_{{{L}^{2}}}}\frac{\partial {{\left\| A \right\|}_{{{L}^{2}}}}}{\partial t}&\lesssim& 2{{\left\| A \right\|}_{{{L}^{2}}}}{{\left( \int_{{{r}^{*}}=r_{1}^{*}}^{{{r}^{*}}=r_{2}^{*}}{\int\limits_{{{S}^{2}}}{{{\left| {{\partial }_{t}}A \right|}^{2}}}}\text{ }{{r}^{2}}\left( 1-\mu  \right)^{1/2}d{{{\bar{r}}}^{*}}d{{\sigma }^{2}} \right)}^{1/2}} \notag\\ 
 \frac{\partial {{\left\| A \right\|}_{{{L}^{2}}}}}{\partial t}&\lesssim& {{\left( \int_{{{r}^{*}}=r_{1}^{*}}^{{{r}^{*}}=r_{2}^{*}}{\int\limits_{{{S}^{2}}}{{{\left| {{\partial }_{t}}A \right|}^{2}}}}\text{ }{{r}^{2}}\left( 1-\mu  \right)^{1/2}d{{{\bar{r}}}^{*}}d{{\sigma }^{2}} \right)}^{1/2}} \notag\\ 
{{\left\| A \right\|}_{{{L}^{2}}}}&\lesssim& \int{\frac{{{E}^{K}}\left( t \right)}{{{t}^{2}}}dt}  
\end{eqnarray}
where,
\begin{eqnarray}
    \left\| A \right\|_{{{L}^{2}}}^{2}=\int_{{{r}^{*}}=r_{1}^{*}}^{{{r}^{*}}=r_{2}^{*}}{\int\limits_{{{S}^{2}}}{{{\left| A \right|}^{2}}}}{{r}^{2}}\left( 1-\mu  \right)^{1/2}d{{\bar{r}}^{*}}d{{\sigma }^{2}}.
\end{eqnarray}
Then, we get
\begin{eqnarray}\label{A}
{{\left\| A \right\|}_{{{L}^{2}}}}\lesssim \int{\frac{{{E}^{K}}\left( t \right)}{{{t}^{2}}}dt}~. 
\end{eqnarray}
Similarly to the case of $F$, one can obtain
\begin{eqnarray}\label{Aint}
    \int_{{{S}^{2}}}{{{r_A}^{2}}{{\left| A \right|}^{2}}\left( t,{{r}_{A}},\theta ,\varphi  \right)d{{\sigma }^{2}}}\lesssim \frac{{{\left( {\hat{E}^{MH}} \right)}^{2}}}{{{t}^{2}}}
\end{eqnarray}
In the same way, it follows that

\begin{eqnarray}
\int\limits_{{{S}^{2}}}r_{{{\slashed{\mathcal{L}}}_{A}}}^{2}{{\left| \slashed{\mathcal{L}}A \right|}^{2}}\left( t,{{r}_{{{\slashed{\mathcal{L}}}_{A}}}},\theta ,\varphi  \right)d{{\sigma }^{2}}&=&\sum\limits_{j=1}^{3}{\frac{1}{r_{{{\slashed{\mathcal{L}}}_{A}}}^{2}}\int_{{{S}^{2}}}{r_{{{\slashed{\mathcal{L}}}_{A}}}^{2}{{\left| {{\mathcal{L}}_{{{\Omega }_{j}}}}A \right|}^{2}}\left( t,{{r}_{{{\slashed{\mathcal{L}}}_{A}}}},\theta ,\varphi  \right)d{{\sigma }^{2}}}} \notag\\ 
&\lesssim&\sum\limits_{j=1}^{3}{\frac{1}{{{R}^{2}}}}\int_{{{S}^{2}}}{r_{{{\slashed{\mathcal{L}}}_{A}}}^{2}{{\left| {{\mathcal{L}}_{{{\Omega }_{j}}}}A \right|}^{2}}\left( t,{{r}_{{{\slashed{\mathcal{L}}}_{A}}}},\theta ,\varphi  \right)d{{\sigma }^{2}}} \notag\\ 
&\lesssim& \int\sum\limits_{j=1}^{3}{\frac{E_{{{\mathcal{L}}_{{{\Omega }_{j}}}}}^{K}\left( t \right)}{{{t}^{2}}}} dt\notag\\
&\lesssim&\frac{{{\left( {\hat{E}^{MH}} \right)}^{2}}}{{{t}^{2}}}~,
\end{eqnarray}

\begin{eqnarray}
	\int\limits_{{{S}^{2}}}{r_{\slashed{\mathcal{L}}{{\slashed{\mathcal{L}}}_{A}}}^{2}{{\left| \slashed{\mathcal{L}}\slashed{\mathcal{L}}F \right|}^{2}}\left( t,{{r}_{{{\slashed{\mathcal{L}}}_{F}}}},\theta ,\varphi  \right)d{{\sigma }^{2}}}\lesssim \frac{{{\left( {\hat{E}^{MH}} \right)}^{2}}}{{{t}^{2}}}~.
\end{eqnarray}
On the other hand,
\begin{eqnarray}
	\int_{{{r}^{*}}=r_{F }^{*}}^{{{r}^{*}}={{r}^{*}}}{\int\limits_{{{S}^{2}}}{\bar{r}{{\left| A  \right|}^{2}}}}\left( 1-\mu  \right)\left( t,\bar{r},\theta,\varphi \right)d{{{\bar{r}}}^{*}}d{{\sigma }^{2}}&\lesssim&\int_{{{r}^{*}}=r_{F }^{*}}^{{{r}^{*}}={{r}^{*}}}{\frac{{\bar{r}}}{R}\int\limits_{{{S}^{2}}}{\bar{r}{{\left| A  \right|}^{2}}}}\left( 1-\mu  \right)\left( t,\bar{r},\theta,\varphi \right)d{{{\bar{r}}}^{*}}d{{\sigma }^{2}} \notag\\ 
	&\lesssim& \frac{{{\left( {\hat{E}^{MH}} \right)}^{2}}}{{{t}^{2}}}+\frac{{{\left( {\hat{E}^{MH}} \right)}^{2}}}{{{w}^{2}}}~.
\end{eqnarray}
By the same method, one can get
\begin{eqnarray}
	\int_{{{r}^{*}}=r_{{{\slashed{\mathcal{L}}}_{F}}}^{*}}^{{{r}^{*}}={{r}^{*}}}{\int\limits_{{{S}^{2}}}{\bar{r}{{\left| \slashed{\mathcal{L}}F \right|}^{2}}}}\left( 1-\mu  \right)\left( t,\bar{r},\theta ,\varphi  \right)d{{\bar{r}}^{*}}d{{\sigma }^{2}}\lesssim\frac{{{\left( {\hat{E}^{MH}} \right)}^{2}}}{{{t}^{2}}}+\frac{{{\left( {\hat{E}^{MH}} \right)}^{2}}}{{{w}^{2}}}~,
\end{eqnarray}
\begin{eqnarray}
	\int_{{{r}^{*}}=r_{\slashed{\mathcal{L}}{{\slashed{\mathcal{L}}}_{F}}}^{*}}^{{{r}^{*}}={{r}^{*}}}{\int\limits_{{{S}^{2}}}{\bar{r}{{\left| \slashed{\mathcal{L}}\slashed{\mathcal{L}}F \right|}^{2}}}}\left( 1-\mu  \right)\left( t,\bar{r},\theta ,\varphi  \right)d{{\bar{r}}^{*}}d{{\sigma }^{2}}\lesssim \frac{{{\left( {\hat{E}^{MH}} \right)}^{2}}}{{{t}^{2}}}+\frac{{{\left( {\hat{E}^{MH}} \right)}^{2}}}{{{w}^{2}}}~.
\end{eqnarray}
Thus, after some calculation as in $F$, we attain
\begin{eqnarray}\label{Aa}
    \left| A \right|\left( t,r,\theta ,\varphi  \right)&\lesssim& \frac{{\hat{E}^{MH}}}{rt}+\frac{{\hat{E}^{MH}}}{rw}\notag\\
    &\lesssim& \frac{\left(\hat{E}^{MH}\right)^2}{1+\left| v \right|}
\end{eqnarray}

\subsection{Decay for the scalar field}
Next, for the uniform bound of scalar fields, let us first again consider the region $w\ge 1,r\ge R$, where $R$ is fixed. As in the preceding subsection, we have Sobolev inequality as follows
\begin{eqnarray}\label{scalarineq}
{{r}^{2}}{{\left| \phi  \right|}^{2}}\le \int_{{{S}^{2}}}{{{r}^{2}}{{\left| \phi  \right|}^{2}}d{{\sigma }^{2}}}+\int_{{{S}^{2}}}{{{r}^{2}}{{\left| \slashed{\mathcal{L}}\left| \phi  \right| \right|}^{2}}d{{\sigma }^{2}}}+\int_{{{S}^{2}}}{{{r}^{2}}{{\left| \slashed{\mathcal{L}}\slashed{\mathcal{L}}\left| \phi  \right| \right|}^{2}}d{{\sigma }^{2}}}~.
\end{eqnarray}
Let ${{r}_{\phi }}$ be a value of $r$ such that $R\le {{r}_{\phi }}\le R+1$, we have 
\begin{eqnarray}
\int_{{{S}^{2}}}{{{r}^{2}}{{\left| \phi  \right|}^{2}}\left( t,r,\theta,\varphi \right)d{{\sigma }^{2}}}&\lesssim& \int_{{{S}^{2}}}{{{r_{\phi}}^{2}}{{\left| \phi  \right|}^{2}}\left( t,{{r}_{\phi }},\theta,\varphi \right)d{{\sigma }^{2}}}+\int_{{{S}^{2}}}{\int_{{{{\bar{r}}}^{*}}=r_{\phi }^{*}}^{{{{\bar{r}}}^{*}}={{r}^{*}}}{{{\nabla }_{{{r}^{*}}}}\left[ {{r}^{2}}{{\left| \phi  \right|}^{2}} \right]}\left( t,r,\theta,\varphi \right)d{{{\bar{r}}}^{*}}d{{\sigma }^{2}}} \notag\\ 
&\lesssim& \int_{{{S}^{2}}}{r_{\phi }^{2}{{\left| \phi  \right|}^{2}}\left( t,{{r}_{\phi }},\theta,\varphi \right)d{{\sigma }^{2}}}+\int_{{{S}^{2}}}{\int_{{{{\bar{r}}}^{*}}=r_{\phi }^{*}}^{{{{\bar{r}}}^{*}}={{r}^{*}}}{2r{{\left| \phi  \right|}^{2}}}\left( t,r,\theta,\varphi \right)\left( 1-\mu  \right)d{{{\bar{r}}}^{*}}d{{\sigma }^{2}}} \notag\\ 
&+&\int_{{{S}^{2}}}{\int_{{{{\bar{r}}}^{*}}=r_{\phi }^{*}}^{{{{\bar{r}}}^{*}}={{r}^{*}}}{{{r}^{2}}{{\nabla }_{{{r}^{*}}}}{{\left| \phi  \right|}^{2}}}\left( t,r,\theta,\varphi \right)d{{{\bar{r}}}^{*}}d{{\sigma }^{2}}}  
\end{eqnarray}
Similar to \eqref{A}, we obtain
\begin{eqnarray}
\int_{{{S}^{2}}}{r_{\phi }^{2}{{\left| \phi  \right|}^{2}}\left( t,{{r}_{\phi }},\theta,\varphi \right)d{{\sigma }^{2}}}\le \frac{{{E}^{K}}\left( t \right)}{{{t}^{2}}}. 
\end{eqnarray}
On the other hand,
\begin{eqnarray}
\int_{{{r}^{*}}=r_{\phi }^{*}}^{{{r}^{*}}={{r}^{*}}}{\int\limits_{{{S}^{2}}}{\bar{r}{{\left| \phi  \right|}^{2}}}}\left( 1-\mu  \right)\left( t,\bar{r},\theta,\varphi \right)d{{{\bar{r}}}^{*}}d{{\sigma }^{2}}&\lesssim&\int_{{{r}^{*}}=r_{\phi }^{*}}^{{{r}^{*}}={{r}^{*}}}{\frac{{\bar{r}}}{R}\int\limits_{{{S}^{2}}}{\bar{r}{{\left| \phi  \right|}^{2}}}}\left( 1-\mu  \right)\left( t,\bar{r},\theta,\varphi \right)d{{{\bar{r}}}^{*}}d{{\sigma }^{2}} \notag\\ 
&\lesssim& \frac{{{E}^{K}}\left( t \right)}{{{t}^{2}}}+\frac{{{E}^{K}}\left( t \right)}{{{w}^{2}}},  
\end{eqnarray}
and
\begin{eqnarray}
\int_{{{S}^{2}}}{\int_{{{{\bar{r}}}^{*}}=r_{\phi }^{*}}^{{{{\bar{r}}}^{*}}={{r}^{*}}}{{{r}^{2}}{{\nabla }_{{{r}^{*}}}}{{\left| \phi  \right|}^{2}}}\left( t,r,\theta,\varphi \right)d{{{\bar{r}}}^{*}}d{{\sigma }^{2}}} &\lesssim& {{\left[ \int_{{{S}^{2}}}{\int_{{{{\bar{r}}}^{*}}=r_{\phi }^{*}}^{{{{\bar{r}}}^{*}}={{r}^{*}}}{{{r}^{2}}{{\left| {{\nabla }_{{{{\hat{r}}}^{*}}}}\phi  \right|}^{2}}}\left( t,r,\theta,\varphi \right)\left( 1-\mu  \right)d{{{\bar{r}}}^{*}}d{{\sigma }^{2}}} \right]}^{1/2}} \notag\\ 
&\cdot& {{\left[ \int_{{{S}^{2}}}{\int_{{{{\bar{r}}}^{*}}=r_{\phi }^{*}}^{{{{\bar{r}}}^{*}}={{r}^{*}}}{{{r}^{2}}{{\left| \phi  \right|}^{2}}}\left( t,r,\theta,\varphi \right)d{{{\bar{r}}}^{*}}d{{\sigma }^{2}}} \right]}^{1/2}} \notag\\ 
&\lesssim& \frac{{{E}^{K}}\left( t \right)}{{{t}^{2}}}.  
\end{eqnarray}
Thus, we get the estimate for the first term of equation (\ref{scalarineq})
\begin{eqnarray}
\int_{{{S}^{2}}}{{{r}^{2}}{{\left| \phi  \right|}^{2}}\left( t,r,\theta,\varphi \right)d{{\sigma }^{2}}}\lesssim \frac{{{E}^{K}}\left( t \right)}{{{t}^{2}}}+\frac{{{E}^{K}}\left( t \right)}{{{w}^{2}}}~.
\end{eqnarray}
Similarly, estimates for the second and third term of (\ref{scalarineq}) can be derived as follows 
\begin{eqnarray}
\int_{{{S}^{2}}}{{{r}^{2}}{{\left| \slashed{\mathcal{L}}\phi  \right|}^{2}}\left( t,r,\theta,\varphi \right)d{{\sigma }^{2}}}\le \sum\limits_{j=1}^{3}{\frac{E_{{{\mathcal{L}}_{{{\Omega }_{j}}}}}^{K}\left( t \right)}{{{t}^{2}}}+\frac{E_{{{\mathcal{L}}_{{{\Omega }_{j}}}}}^{K}\left( t \right)}{{{w}^{2}}}}~,
\end{eqnarray}
\begin{eqnarray}
\int_{{{S}^{2}}}{{{r}^{2}}{{\left| \slashed{\mathcal{L}}\slashed{\mathcal{L}}\phi  \right|}^{2}}\left( t,r,\theta,\varphi \right)d{{\sigma }^{2}}}\le \sum\limits_{i=1}^{3}{\sum\limits_{j=1}^{3}{\left[ \frac{E_{{{\mathcal{L}}_{{{\Omega }_{i}}}}{{\mathcal{L}}_{{{\Omega }_{j}}}}}^{K}\left( t \right)}{{{t}^{2}}}+\frac{E_{{{\mathcal{L}}_{{{\Omega }_{i}}}}{{\mathcal{L}}_{{{\Omega }_{j}}}}}^{K}\left( t \right)}{{{w}^{2}}} \right]}}~.
\end{eqnarray}
Therefore, it can be demonstrated that
\begin{eqnarray}
\left| \phi  \right|\lesssim \frac{{{\left[ E_{{{\mathcal{L}}_{{{\Omega }_{j}}}},{{\mathcal{L}}_{{{\Omega }_{i}}}}{{\mathcal{L}}_{{{\Omega }_{j}}}}}^{K}\left( t \right) \right]}^{1/2}}}{rt}+\frac{{{\left[ E_{{{\mathcal{L}}_{{{\Omega }_{j}}}},{{\mathcal{L}}_{{{\Omega }_{i}}}}{{\mathcal{L}}_{{{\Omega }_{j}}}}}^{K}\left( t \right) \right]}^{1/2}}}{rw}~,
\end{eqnarray}
where
\begin{eqnarray}
E_{{{\mathcal{L}}_{{{\Omega }_{j}}}},{{\mathcal{L}}_{{{\Omega }_{i}}}}{{\mathcal{L}}_{{{\Omega }_{j}}}}}^{K}\left( t \right)={{E}^{K}}\left( t \right)+\sum\limits_{j=1}^{3}{E_{{{\mathcal{L}}_{{{\Omega }_{j}}}}}^{K}\left( t \right)}+\sum\limits_{i=1}^{3}{\sum\limits_{j=1}^{3}{E_{{{\mathcal{L}}_{{{\Omega }_{i}}}}{{\mathcal{L}}_{{{\Omega }_{j}}}}}^{K}\left( t \right)}}~.
\end{eqnarray}
For $R$ fixed, consider the region where $t\ge 1$ such that we have $r+t\le rt$. Since $w\ge 1$, consequently $v={{r}^{*}}+t\lesssim 2r+t\lesssim 2r+2t\lesssim rt$ and $v\lesssim r+t\lesssim r+t-{{r}^{*}}\lesssim r+w\le rw$. Then 
\begin{eqnarray}\label{scalarregion}
\left| \phi  \right|\lesssim \frac{{{\left( E_{{{\mathcal{L}}_{{{\Omega }_{j}}}},{{\mathcal{L}}_{{{\Omega }_{i}}}}{{\mathcal{L}}_{{{\Omega }_{j}}}}}^{K}\left( t \right) \right)}^{1/2}}}{v}~.
\end{eqnarray}
For the region $t\le 1$, $w\ge 1,r\ge R$ is a bounded compact region, so in this region
\begin{eqnarray}\label{scalarcomregion}
\left| \phi  \right|\lesssim E_{{{\mathcal{L}}_{{{\Omega }_{j}}}},{{\mathcal{L}}_{{{\Omega }_{i}}}}{{\mathcal{L}}_{{{\Omega }_{j}}}}}^{\hat{t}}\left( t={{t}_{0}} \right)~.
\end{eqnarray}
Thus, combining the result in \eqref{scalarregion} and \eqref{scalarcomregion}, we attain \eqref{decayscalar}.

\subsection{Decay Estimate for \texorpdfstring{$D\phi$}{}}

Now, we want to get the estimate for $\left|D\phi\right|$, from the definition, we have
\begin{eqnarray}
	\left| D\phi  \right|\le \left| \nabla \phi  \right|+\left| A \right|\left| \phi  \right|
\end{eqnarray}
Similarly as in the previous result of $\phi$, we obtain
\begin{eqnarray}\label{curvaturpi}
	{{r}^{2}}{{\left| \nabla \phi  \right|}^{2}}\le \int_{{{S}^{2}}}{{{r}^{2}}{{\left| \nabla \phi  \right|}^{2}}d{{\sigma }^{2}}}+\int_{{{S}^{2}}}{{{r}^{2}}{{\left| \slashed{\mathcal{L}}\left| \nabla \phi  \right| \right|}^{2}}d{{\sigma }^{2}}}+\int_{{{S}^{2}}}{{{r}^{2}}{{\left| \slashed{\mathcal{L}}\slashed{\mathcal{L}}\left| \nabla \phi  \right| \right|}^{2}}d{{\sigma }^{2}}}
\end{eqnarray}
The first term of (\ref{curvaturpi}) can be estimated as follows
\begin{eqnarray}\label{firsttermcurpi}
	\int_{{{S}^{2}}}{{{r}^{2}}{{\left| \nabla \phi  \right|}^{2}}\left( t,r,\theta,\varphi \right)d{{\sigma }^{2}}}&\lesssim& \int_{{{S}^{2}}}{{{r_{\phi}}^{2}}{{\left| \nabla \phi  \right|}^{2}}d{{\sigma }^{2}}}+\int_{{{S}^{2}}}{\int_{{{{\bar{r}}}^{*}}=r_{\phi }^{*}}^{{{{\bar{r}}}^{*}}={{r}^{*}}}{{{\nabla }_{{{r}^{*}}}}\left[ {{r}^{2}}{{\left| \nabla \phi  \right|}^{2}} \right]}d{{{\bar{r}}}^{*}}d{{\sigma }^{2}}} \notag\\ 
	&\lesssim& \int_{{{S}^{2}}}{r_{\phi }^{2}{{\left| \nabla \phi  \right|}^{2}}d{{\sigma }^{2}}}+\int_{{{S}^{2}}}{\int_{{{{\bar{r}}}^{*}}=r_{\phi }^{*}}^{{{{\bar{r}}}^{*}}={{r}^{*}}}{2r{{\left| \nabla \phi  \right|}^{2}}}\left( 1-\mu  \right)d{{{\bar{r}}}^{*}}d{{\sigma }^{2}}} \notag\\ 
	&+&\int_{{{S}^{2}}}{\int_{{{{\bar{r}}}^{*}}=r_{\phi }^{*}}^{{{{\bar{r}}}^{*}}={{r}^{*}}}{{{r}^{2}}{{\nabla }_{{{r}^{*}}}}{{\left| \nabla \phi  \right|}^{2}}}d{{{\bar{r}}}^{*}}d{{\sigma }^{2}}}~.  	
\end{eqnarray}
Let  ${{r}_{\phi }}$ be a constant such that $R\le {{r}_{\phi }}\le R+1$ and we have estimate for the first and second term of (\ref{firsttermcurpi}),
\begin{eqnarray}
	\int_{{{S}^{2}}}{r_{\phi }^{2}{{\left| \nabla \phi  \right|}^{2}}\left( t,{{r}_{\phi }},\theta,\varphi \right)d{{\sigma }^{2}}}&\lesssim& \frac{{{E}^{K}}\left( t \right)}{{{t}^{2}}},\\
	\int_{{{r}^{*}}=r_{\phi }^{*}}^{{{r}^{*}}={{r}^{*}}}{\int\limits_{{{S}^{2}}}{\bar{r}{{\left| \nabla \phi  \right|}^{2}}}}\left( 1-\mu  \right)\left( t,\bar{r},\theta,\varphi \right)d{{\bar{r}}^{*}}d{{\sigma }^{2}}&\lesssim& \frac{{{E}^{K}}\left( t \right)}{{{t}^{2}}}+\frac{{{E}^{K}}\left( t \right)}{{{w}^{2}}}.
\end{eqnarray}
We can write the last term (\ref{firsttermcurpi}) as
\begin{eqnarray}
	\int_{{{S}^{2}}}{\int_{{{{\bar{r}}}^{*}}=r_{\phi }^{*}}^{{{{\bar{r}}}^{*}}={{r}^{*}}}{{{r}^{2}}{{\nabla }_{{{r}^{*}}}}{{\left| \nabla \phi  \right|}^{2}}}\left( t,r,\theta,\varphi \right)d{{{\bar{r}}}^{*}}d{{\sigma }^{2}}}&=&\int_{{{S}^{2}}}{\int_{{{{\bar{r}}}^{*}}=r_{\phi }^{*}}^{{{{\bar{r}}}^{*}}={{r}^{*}}}{2{{r}^{2}}\left| \nabla \phi  \right|{{\nabla }_{{{{\hat{r}}}^{*}}}}\left| \nabla \phi  \right|}\left( t,r,\theta,\varphi \right)d{{{\bar{r}}}^{*}}d{{\sigma }^{2}}} \notag\\ 
	&\lesssim& {{\left( \int_{{{S}^{2}}}{\int_{{{{\bar{r}}}^{*}}=r_{\phi }^{*}}^{{{{\bar{r}}}^{*}}={{r}^{*}}}{{{r}^{2}}{{\left| \nabla \nabla \phi  \right|}^{2}}}\left( t,r,\theta,\varphi \right)d{{{\bar{r}}}^{*}}d{{\sigma }^{2}}} \right)}^{1/2}} \notag\\ 
	&\cdot& {{\left( \int_{{{S}^{2}}}{\int_{{{{\bar{r}}}^{*}}=r_{\phi }^{*}}^{{{{\bar{r}}}^{*}}={{r}^{*}}}{{{r}^{2}}{{\left| \nabla \phi  \right|}^{2}}}\left( t,r,\theta,\varphi \right)d{{{\bar{r}}}^{*}}d{{\sigma }^{2}}} \right)}^{1/2}} \notag\\ 
	&\lesssim& \frac{{{E}^{K}}\left( t \right)}{{{t}^{2}}}{{\left( \int_{{{S}^{2}}}{\int_{{{{\bar{r}}}^{*}}=r_{\phi }^{*}}^{{{{\bar{r}}}^{*}}={{r}^{*}}}{{{r}^{2}}{{\left| \nabla \nabla \phi  \right|}^{2}}}\left( t,r,\theta,\varphi \right)d{{{\bar{r}}}^{*}}d{{\sigma }^{2}}} \right)}^{1/2}}~. \notag\\ 
\end{eqnarray}
Using the equation of motion in (\ref{eom2}), we get the inequality

\begin{eqnarray}\label{secondderivatifphi}
	\int_{{{S}^{2}}}{\int_{{{{\bar{r}}}^{*}}=r_{\phi }^{*}}^{{{{\bar{r}}}^{*}}={{r}^{*}}}{{{r}^{2}}{{\left| \nabla \nabla \phi  \right|}^{2}}}\left( t,r,\theta,\varphi \right)d{{{\bar{r}}}^{*}}d{{\sigma }^{2}}}&\lesssim& \int_{{{S}^{2}}}{\int_{{{{\bar{r}}}^{*}}=r_{\phi }^{*}}^{{{{\bar{r}}}^{*}}={{r}^{*}}}{{{r}^{2}}{{\left| \nabla_{\phi }P \right|}^{2}}}\left( t,r,\theta,\varphi \right)d{{{\bar{r}}}^{*}}d{{\sigma }^{2}}} \notag\\ 
	&+&\int_{{{S}^{2}}}{\int_{{{{\bar{r}}}^{*}}=r_{\phi }^{*}}^{{{{\bar{r}}}^{*}}={{r}^{*}}}{{{r}^{2}}{{\left| \nabla A\phi  \right|}^{2}}}\left( t,r,\theta,\varphi \right)d{{{\bar{r}}}^{*}}d{{\sigma }^{2}}} \notag\\ 
	&+&\int_{{{S}^{2}}}{\int_{{{{\bar{r}}}^{*}}=r_{\phi }^{*}}^{{{{\bar{r}}}^{*}}={{r}^{*}}}{{{r}^{2}}{{\left| A\nabla \phi  \right|}^{2}}}\left( t,r,\theta,\varphi \right)d{{{\bar{r}}}^{*}}d{{\sigma }^{2}}} \notag\\ 
	&+&\int_{{{S}^{2}}}{\int_{{{{\bar{r}}}^{*}}=r_{\phi }^{*}}^{{{{\bar{r}}}^{*}}={{r}^{*}}}{{{r}^{2}}{{\left| AD\phi  \right|}^{2}}}\left( t,r,\theta,\varphi \right)d{{{\bar{r}}}^{*}}d{{\sigma }^{2}}}~. \notag\\ 
\end{eqnarray}
In order to have an estimate of (\ref{secondderivatifphi}), we have to specify the form of the scalar potential $P\left(\phi,\bar{\phi}\right)$. It is worth mentioning that there four known examples in the field theories, namely, the mass term, the ${{\phi }^{4}}$-theory, the sine-Gordon, and the  Toda field theories, respectively.
\begin{assumption}\label{asumsiP}
    The scalar potential $P\left(\phi,\bar{\phi}\right)$ has to be either of the following form
    \begin{eqnarray}
    P\left( |\phi| \right)&=&{{c}_{1}}{{|\phi| }^{2}}\label{mass}\\
        P\left( |\phi| \right)&=&{{c}_{2}}{{|\phi| }^{4}}\label{pi4}\\
        P\left( |\phi|\right)&=&{{c}_{3}}\left( 1-\cos \eta |\phi|\right)\label{sine}\\
       P\left( |\phi| \right)&=&{{c}_{4}}{{e}^{-\lambda |\phi| }}\label{toda}
    \end{eqnarray}
$c_i,\eta,$ and $\lambda$ are real constants with  $\lambda >0$.
\end{assumption}
By applying the given method, we can obtain the estimate for the first term of (\ref{secondderivatifphi}), 
\begin{equation}
\tilde{P}=\int_{{{S}^{2}}}{\int_{{{{\bar{r}}}^{*}}=r_{\phi }^{*}}^{{{{\bar{r}}}^{*}}={{r}^{*}}}{{{r}^{2}}{{\left| {{\nabla }_{\phi }}\mathcal{P} \right|}^{2}}}d{{{\bar{r}}}^{*}}d{{\sigma }^{2}}} \lesssim 
\begin{cases}
\frac{\left(E^{k}(t)\right)^{3/2}}{t^3} ~ , &\text{if $P$ is of the form \eqref{pi4}} ~ ,\\
\frac{E^{k}(t)}{t} ~ , &\text{if $P$ is is either of the form \eqref{mass},\eqref{sine},\eqref{toda}}.
\end{cases}
\end{equation}
Using the Holder inequality, the estimate for the residual term can be obtained, which is given as follows,
\begin{eqnarray}
	\int_{{{S}^{2}}}{\int_{{{{\bar{r}}}^{*}}=r_{\phi }^{*}}^{{{{\bar{r}}}^{*}}={{r}^{*}}}{{{r}^{2}}{{\left| \nabla A\phi  \right|}^{2}}}\left( t,r,\theta,\varphi \right)d{{{\bar{r}}}^{*}}d{{\sigma }^{2}}}&\lesssim& {{\left| \phi  \right|}_{\infty }}\left( \int_{{{S}^{2}}}{\int_{{{{\bar{r}}}^{*}}=r_{\phi }^{*}}^{{{{\bar{r}}}^{*}}={{r}^{*}}}{{{r}^{2}}{{\left| \nabla A \right|}^{2}}}\left( t,r,\theta,\varphi \right)d{{{\bar{r}}}^{*}}d{{\sigma }^{2}}} \right) \notag\\ 
	&\lesssim& \frac{{{E}^{K}}\left( t \right)}{{{t}^{2}}}~,
\end{eqnarray}
where
\begin{eqnarray}
	{{\left| \phi  \right|}_{\infty }}=\underset{{{{\bar{r}}}^{*}},\theta ,\varphi }{\mathop{\text{suf}}}\,\left| \phi \left( {{{\bar{r}}}^{*}},\theta ,\varphi  \right) \right|<\infty~,
\end{eqnarray}
\begin{eqnarray}
	\int_{{{S}^{2}}}{\int_{{{{\bar{r}}}^{*}}=r_{\phi }^{*}}^{{{{\bar{r}}}^{*}}={{r}^{*}}}{{{r}^{2}}{{\left| A\nabla \phi  \right|}^{2}}}\left( t,r,\theta,\varphi \right)d{{{\bar{r}}}^{*}}d{{\sigma }^{2}}}&\lesssim& {{\left| A \right|}_{\infty }}\int_{{{S}^{2}}}{\int_{{{{\bar{r}}}^{*}}=r_{\phi }^{*}}^{{{{\bar{r}}}^{*}}={{r}^{*}}}{{{r}^{2}}{{\left| \nabla \phi  \right|}^{2}}}\left( t,r,\theta,\varphi \right)d{{{\bar{r}}}^{*}}d{{\sigma }^{2}}} \notag\\ 
	&\lesssim& \frac{{{E}^{K}}\left( t \right)}{{{t}^{2}}}~,
\end{eqnarray}
and
\begin{eqnarray}
	\int_{{{S}^{2}}}{\int_{{{{\bar{r}}}^{*}}=r_{\phi }^{*}}^{{{{\bar{r}}}^{*}}={{r}^{*}}}{{{r}^{2}}{{\left| A D \phi  \right|}^{2}}}\left( t,r,\theta,\varphi \right)d{{{\bar{r}}}^{*}}d{{\sigma }^{2}}}&\lesssim& {{\left| A \right|}_{\infty }}\int_{{{S}^{2}}}{\int_{{{{\bar{r}}}^{*}}=r_{\phi }^{*}}^{{{{\bar{r}}}^{*}}={{r}^{*}}}{{{r}^{2}}{{\left| D \phi  \right|}^{2}}}\left( t,r,\theta,\varphi \right)d{{{\bar{r}}}^{*}}d{{\sigma }^{2}}} \notag\\ 
	&\lesssim&\frac{{{E}^{K}}\left( t \right)}{{{t}^{2}}}~.
\end{eqnarray}
Therefore, one can derive
\begin{eqnarray}
	\int_{{{S}^{2}}}{\int_{{{{\bar{r}}}^{*}}=r_{\phi }^{*}}^{{{{\bar{r}}}^{*}}={{r}^{*}}}{{{r}^{2}}{{\left| \nabla \nabla \phi  \right|}^{2}}}\left( t,r,\theta,\varphi \right)d{{{\bar{r}}}^{*}}d{{\sigma }^{2}}}\lesssim \frac{{{E}^{K}}}{{{t}^{2}}}+\tilde{P}~.
\end{eqnarray}
Consequently, 
\begin{eqnarray}
	\int_{{{S}^{2}}}{\int_{{{{\bar{r}}}^{*}}=r_{\phi }^{*}}^{{{{\bar{r}}}^{*}}={{r}^{*}}}{{{r}^{2}}{{\nabla }_{{{r}^{*}}}}{{\left| \nabla \phi  \right|}^{2}}}\left( t,r,\theta,\varphi \right)d{{{\bar{r}}}^{*}}d{{\sigma }^{2}}}&\lesssim& \frac{{{E}^{K}}\left( t \right)}{{{t}^{2}}}\left(\frac{{{E}^{K}}}{{{t}^{2}}}+\tilde{P}\right)~,
\end{eqnarray}
Hence, it can be shown that the estimate for the first term of \eqref{curvaturpi} satisfy the following inequality
\begin{eqnarray}
	\int_{{{S}^{2}}}{{{r}^{2}}{{\left| \nabla \phi  \right|}^{2}}\left( t,r,\theta,\varphi \right)d{{\sigma }^{2}}}\lesssim\frac{{E^K\left(t\right)}}{{{t}^{2}}}+\frac{{E^K\left(t\right)}}{{{w}^{2}}}+\frac{{{\left( {{E}^{K}} \right)}^{2}}}{{{t}^{4}}}+\frac{{{E}^{K}}}{{{t}^{2}}}\tilde{P}~.	
\end{eqnarray}
In similar way, one can derive the estimate for the second and last term of \eqref{curvaturpi} as follows,
\begin{eqnarray}
	\int_{{{S}^{2}}}{{{r}^{2}}{{\left| \slashed{\mathcal{L}}\left| \nabla \phi  \right| \right|}^{2}}d{{\sigma }^{2}}}\lesssim \frac{\sum\limits_{j=1}^{3}{E_{{{\mathcal{L}}_{{{\Omega }_{j}}}}}^{K}}}{{{t}^{2}}}+\frac{\sum\limits_{j=1}^{3}{E_{{{\mathcal{L}}_{{{\Omega }_{j}}}}}^{K}}}{{{w}^{2}}}+\frac{\sum\limits_{j=1}^{3}{{{\left( E_{{{\mathcal{L}}_{{{\Omega }_{j}}}}}^{K} \right)}^{2}}}}{t^4}+\frac{\sum\limits_{j=1}^{3}{{{\left( E_{{{\mathcal{L}}_{{{\Omega }_{j}}}}}^{K} \right)}^{}}}}{t^2}\tilde{P}~,
\end{eqnarray}
\begin{eqnarray}
	\int_{{{S}^{2}}}{{{r}^{2}}{{\left| \slashed{\mathcal{L}}\slashed{\mathcal{L}}\left| \nabla \phi  \right| \right|}^{2}}d{{\sigma }^{2}}}&&\lesssim \frac{\sum\limits_{i=1}^{3}{\sum\limits_{j=1}^{3}{E_{{{\mathcal{L}}_{{{\Omega }_{j}}}}{{\mathcal{L}}_{{{\Omega }_{j}}}}}^{K}}}}{{{t}^{2}}}+\frac{\sum\limits_{i=1}^{3}{\sum\limits_{j=1}^{3}{E_{{{\mathcal{L}}_{{{\Omega }_{j}}}}{{\mathcal{L}}_{{{\Omega }_{j}}}}}^{K}}}}{{{w}^{2}}}+\frac{\sum\limits_{i=1}^{3}{\sum\limits_{j=1}^{3}{{{\left( E_{{{\mathcal{L}}_{{{\Omega }_{j}}}}{{\mathcal{L}}_{{{\Omega }_{j}}}}}^{K} \right)}^{2}}}}}{t^4}\notag\\
&&+\frac{\sum\limits_{i=1}^{3}{\sum\limits_{j=1}^{3}{{{\left( E_{{{\mathcal{L}}_{{{\Omega }_{j}}}}{{\mathcal{L}}_{{{\Omega }_{j}}}}}^{K} \right)}^{}}}}}{t^2}\tilde{P}~.
\end{eqnarray}
Thus, for $R$ fixed, consider region $t\ge 1$ such that we have $r+t\le rt$. After some computation using the same method as in the previous section, we have
\begin{eqnarray}
	\left| \nabla \phi  \right|\lesssim \frac{{{\left( {\hat{E}^{MH}} \right)}^{5/4}}}{v}~,
\end{eqnarray}
and for the bounded compact region $t\le 1$, $w\ge 1,r\ge R$, it can be demonstrated that
\begin{eqnarray}
	\left| \nabla\phi  \right|\lesssim E_{{{\mathcal{L}}_{{{\Omega }_{j}}}},{{\mathcal{L}}_{{{\Omega }_{i}}}}{{\mathcal{L}}_{{{\Omega }_{j}}}}}^{\hat{t}}\left( t={{t}_{0}} \right)~.
\end{eqnarray}
Finally, we have
\begin{eqnarray}\label{nabphi}
	\left| \nabla \phi  \right|\lesssim\frac{{{\left( {\hat{E}^{MH}} \right)}^{5/4}}}{1+\left|v\right|}~.
\end{eqnarray}
Thus, by using the result in \eqref{nabphi},\eqref{Aa}, and \eqref{decayscalar}, we can write down the estimate for $\left|D\phi\right|$ as in \eqref{decayDphi}.

This complete the proof of theorem \ref{teorema1}.

\section{Proof of Theorem 2}
\label{sec:decaynear}
In this section, we prove the decay estimate in the photon sphere region $2m<r<3m$. This will be done by integrating the Penrose diagram into rectangles representing the exterior of the Schwarzschild black holes, one side of which is enclosed by the horizon.

\begin{definition}\label{definitionH}
    A vector field $H$ is defined as
\begin{eqnarray}\label{H}
    H\equiv-\frac{h\left( {{r}^{*}} \right)}{\left( 1-\frac{2m}{r} \right)}\frac{\partial }{\partial w}-h\left( {{r}^{*}} \right)\frac{\partial }{\partial v}~.
\end{eqnarray}

\end{definition}
Then, one can show
\begin{eqnarray}
  {{T}_{\alpha \beta }}{{\pi }^{\alpha \beta }}\left( H \right)&=&\left\{ \frac{1}{{{r}^{2}}}{{\left| {{F }_{w\theta }} \right|}^{2}}+\frac{1}{{{r}^{2}}{{\sin }^{2}}\theta }{{\left| {{F }_{w\varphi }} \right|}^{2}}+{{\left| {{D}_{w}}\phi  \right|}^{2}} \right\}\frac{1}{{{\left( 1-\mu  \right)}^{2}}}\left[ {h}'-\frac{\mu }{r}h \right] \notag\\ 
 &-&\left\{ \frac{1}{{{r}^{2}}}{{\left| {{F }_{v\theta }} \right|}^{2}}+\frac{1}{{{r}^{2}}{{\sin }^{2}}\theta }{{\left| {{F }_{v\varphi }} \right|}^{2}}+{{\left| {{D}_{v}}\phi  \right|}^{2}} \right\}\frac{{{h}'}}{\left( 1-\mu  \right)}-\left\{ \frac{1}{{{(1-\mu )}^{2}}}{{\left| {{F}_{vw}} \right|}^{2}} \right. \notag\\ 
 &+&\left. \frac{1}{4{{r}^{4}}{{\sin }^{2}}\theta }{{\left| {{F }_{\theta \varphi }} \right|}^{2}} \right\}\left( \frac{1}{\left( 1-\mu  \right)}\left[ {h}'-\frac{\mu }{r}h \right]-{h}'+\frac{\left( 3\mu -2 \right)}{r}\left[ \frac{h}{\left( 1-\mu  \right)}-h \right] \right) \notag\\ 
 &-&\left\{ \frac{2}{\left( 1-\mu  \right)}{{D}_{v}}\phi \overline{{{D}_{w}}\phi }+\frac{1}{2{{r}^{2}}}{{\left| {{D}_{\theta }}\phi  \right|}^{2}}+\frac{1}{2{{r}^{2}}{{\sin }^{2}}\theta }{{\left| {{D}_{\varphi }}\phi  \right|}^{2}}+P \right\} \notag\\ 
 &\times& \left( \frac{1}{\left( 1-\mu  \right)}\left[ {h}'-\frac{\mu }{r}h \right]-{h}'+\frac{\mu }{r}\left[ \frac{h}{\left( 1-\mu  \right)}-h \right] \right) \notag\\ 
 &+&\frac{\left( 1-\mu  \right)\left( 2+{{\sin }^{2}}\theta  \right)}{4r}P\left[ \frac{h}{\left( 1-\mu  \right)}-h \right] ~. 
\end{eqnarray}
\begin{definition}\label{defnear1}
The functional energy $\mathcal{E}^{H}$ are the flux of $H$ passing through the hypersurfaces given by
\begin{eqnarray}\label{eH1}
  {{\mathcal{E}}^{H}}\left( w={{w}_{i}} \right)\left( {{v}_{i}}\le v\le {{v}_{i+1}} \right)&\equiv&\int_{{{v}_{i}}}^{{{v}_{i+1}}}{\int_{{{S}^{2}}}{{{J}_{\alpha }}\left( H \right){{n}^{\alpha }}\text{dVo}{{\text{l}}_{w={{w}_{i}}}}}} \notag\\ 
 &=&-2\int_{{{v}_{i}}}^{{{v}_{i+1}}}{\int_{{{S}^{2}}}{\left[ \frac{h}{\left( 1-\mu  \right)}\left( \frac{1}{{{(1-\mu )}^{2}}}{{\left| {{F }_{vw}} \right|}^{2}}+\frac{(1-\mu )}{4{{r}^{4}}{{\sin }^{2}}\theta }{{\left| {{F }_{\varphi \theta }} \right|}^{2}} \right. \right.}} \notag\\ 
 &+&\left. {{D}_{v}}\phi \overline{{{D}_{w}}\phi }+\frac{\left( 1-\mu  \right)}{4{{r}^{2}}}{{\left| {{D}_{\theta }}\phi  \right|}^{2}}+\frac{\left( 1-\mu  \right)}{4{{r}^{2}}{{\sin }^{2}}\theta }{{\left| {{D}_{\varphi }}\phi  \right|}^{2}}+\frac{\left( 1-\mu  \right)}{2}P \right) \notag\\ 
 &+&\left. h\left( \frac{1}{{{r}^{2}}}{{\left| {{F }_{v\theta }} \right|}^{2}}+\frac{1}{{{r}^{2}}{{\sin }^{2}}\theta }{{\left| {{F }_{v\varphi }} \right|}^{2}}+{{\left| {{D}_{v}}\phi  \right|}^{2}} \right) \right]{{r}^{2}}d{{\sigma }^{2}}dv~, 
\end{eqnarray}
\end{definition}

\begin{eqnarray}\label{eH2}
  {{\mathcal{E}}^{H}}\left( v={{v}_{i}} \right)\left( {{w}_{i}}\le w\le {{w}_{i+1}} \right)&\equiv&-2\int_{{{w}_{i}}}^{{{w}_{i+1}}}{\int_{{{S}^{2}}}{\left[ \frac{h}{\left( 1-\mu  \right)}\left( \frac{1}{{{(1-\mu )}^{2}}}{{\left| {{F }_{vw}} \right|}^{2}}+\frac{(1-\mu )}{4{{r}^{4}}{{\sin }^{2}}\theta }{{\left| {{F }_{\varphi \theta }} \right|}^{2}} \right. \right.}} \notag\\ 
 &+&\left. {{D}_{v}}\phi \overline{{{D}_{w}}\phi }+\frac{\left( 1-\mu  \right)}{4{{r}^{2}}}{{\left| {{D}_{\theta }}\phi  \right|}^{2}}+\frac{\left( 1-\mu  \right)}{4{{r}^{2}}{{\sin }^{2}}\theta }{{\left| {{D}_{\varphi }}\phi  \right|}^{2}}+\frac{\left( 1-\mu  \right)}{2}P \right) \notag\\ 
 &+&\left. h\left( \frac{1}{{{r}^{2}}}{{\left| {{F }_{w\theta }} \right|}^{2}}+\frac{1}{{{r}^{2}}{{\sin }^{2}}\theta }{{\left| {{F }_{w\varphi }} \right|}^{2}}+{{\left| {{D}_{w}}\phi  \right|}^{2}} \right) \right]{{r}^{2}}d{{\sigma }^{2}}dv ~. 
\end{eqnarray}
Applying the divergence theorem for the field in a rectangle in the Penrose diagram representing the exterior of the Schwarzschild spacetimes of which one side contains the horizon, say in the region $\left[ {{w}_{i}},\infty  \right] \times \left[ {{v}_{i}},{{v}_{i+1}} \right]$, we obtain
\begin{eqnarray}\label{eqIH}
  {{I}^{H}}\left( {{v}_{i}}\le v\le {{v}_{i+1}} \right)\left( {{w}_{i}}\le w\le \infty  \right)&+&{{\mathcal{E}}^{H}}\left( w=\infty  \right)\left( {{v}_{i}}\le v\le {{v}_{i+1}} \right)={{\mathcal{E}}^{H}}\left( w={{w}_{i}} \right)\left( {{v}_{i}}\le v\le {{v}_{i+1}} \right) \notag\\ 
 &+&{{\mathcal{E}}^{H}}\left( v={{v}_{i}} \right)\left( {{w}_{i}}\le w\le \infty  \right)-{{\mathcal{E}}^{H}}\left( v={{v}_{i+1}} \right)\left( {{w}_{i}}\le w\le \infty  \right)~,\notag\\  
\end{eqnarray}
where
\begin{eqnarray}\label{IH}
  {{I}^{H}}&=&\int_{{{v}_{i}}}^{{{v}_{i+1}}}{\int_{{{w}_{i}}}^{{{w}_{i+1}}}{\int_{{{S}^{2}}}{\left\{ \frac{1}{{{r}^{2}}}{{\left| {{F }_{w\theta }} \right|}^{2}}+\frac{1}{{{r}^{2}}{{\sin }^{2}}\theta }{{\left| {{F }_{w\varphi }} \right|}^{2}}+{{\left| {{D}_{w}}\phi  \right|}^{2}} \right\}\frac{1}{{{\left( 1-\mu  \right)}^{2}}}\left[ {h}'-\frac{\mu }{r}h \right]}}} \notag\\ 
 &-&\left\{ \frac{1}{{{r}^{2}}}{{\left| {{F }_{v\theta }} \right|}^{2}}+\frac{1}{{{r}^{2}}{{\sin }^{2}}\theta }{{\left| {{F }_{v\varphi }} \right|}^{2}}+{{\left| {{D}_{v}}\phi  \right|}^{2}} \right\}\frac{{{h}'}}{\left( 1-\mu  \right)}-\left\{ \frac{1}{{{(1-\mu )}^{2}}}{{\left| {{F }_{vw}} \right|}^{2}} \right. \notag\\ 
 &+&\left. \frac{1}{4{{r}^{4}}{{\sin }^{2}}\theta }{{\left| {{F }_{\theta \varphi }} \right|}^{2}} \right\}\left( \frac{1}{\left( 1-\mu  \right)}\left[ {h}'-\frac{\mu }{r}h \right]-{h}'+\frac{\left( 3\mu -2 \right)}{r}\left[ \frac{h}{\left( 1-\mu  \right)}-h \right] \right) \notag\\ 
 &-&\left\{ \frac{2}{\left( 1-\mu  \right)}{{D}_{v}}\phi \overline{{{D}_{w}}\phi }+\frac{1}{2{{r}^{2}}}{{\left| {{D}_{\theta }}\phi  \right|}^{2}}+\frac{1}{2{{r}^{2}}{{\sin }^{2}}\theta }{{\left| {{D}_{\varphi }}\phi  \right|}^{2}}+P \right\} \notag\\ 
 &\times& \left( \frac{1}{\left( 1-\mu  \right)}\left[ {h}'-\frac{\mu }{r}h \right]-{h}'+\frac{\mu }{r}\left[ \frac{h}{\left( 1-\mu  \right)}-h \right] \right) \notag\\ 
 &+& \left. \frac{\left( 1-\mu  \right)\left( 2+{{\sin }^{2}}\theta  \right)}{4r}P\left[ \frac{h}{\left( 1-\mu  \right)}-h \right] \right\}{{r}^{2}}d{{\sigma }^{2}}\left( 1-\mu  \right)dwdv~.  
\end{eqnarray}
Furthermore, let $h$ be supported in the region $2m\le r\le \left( 1.2 \right){{r}_{1}}$ for ${{r}_{1}}$ chosen such that $2m<{{r}_{0}}\le {{r}_{1}}<1.2{{r}_{1}}<3m$. For all $r\le {{r}_{1}}$,  we have 
\begin{eqnarray}
  h&\ge& 0, \label{h1}\\ 
 {h}'&\ge& 0, \label{h2}\\ 
 \frac{\mu }{r}h-{h}'&\ge&0, \label{h3}\\ 
 \frac{3}{r}h-\frac{1}{\left( 1-\mu  \right)}{h}'&\ge& 0 ~.\label{h4}
\end{eqnarray}
For example, we may simply choose $h\left( {{r}^{*}}\right)=r^{*}$ which satisfies \eqref{h1}-\eqref{h4}. Here, we choose $r_1$ small enough such that $1.2 r_1  < 3m$.
Additionally, we have the following estimates: 
\begin{preposition}\label{prep4}
Let define ${{\tilde{\mathcal{E}}}^{H}}=-{{\mathcal{E}}^{H}}\left( w={{w}_{i}} \right)\left( {{v}_{i}}\le v\le {{v}_{i+1}} \right)$. For $\left( {{w}_{i}},{{v}_{i}} \right)$ such that $r\left( {{w}_{i}},{{v}_{i}} \right)={{r}_{1}}$ and for ${{v}_{i+1}}\ge {{v}_{i}}$, we have
\begin{eqnarray}\label{estimate1}
    {{\tilde{\mathcal{E}}}^{H}}\lesssim {{\mathcal{E}}^{\left( {}^{\partial }/{}_{\partial t} \right)}}\left( w={{w}_{i}} \right)\left( {{v}_{i}}\le v\le {{v}_{i+1}} \right)
\end{eqnarray}
\end{preposition}
\begin{proof}
It is worth noting that within the photon sphere region, the energy generated from a time-like vector field, denoted by ${{\mathcal{E}}^{\left( {}^{\partial }/{}_{\partial t} \right)}}$, is given by the following expression
\begin{eqnarray}\label{et}
  {{\mathcal{E}}^{\left( {}^{\partial }/{}_{\partial t} \right)}}\left( w={{w}_{i}} \right)\left( {{v}_{i}}\le v\le {{v}_{i+1}} \right)&=&2\int_{{{v}_{i}}}^{{{v}_{i+1}}}{\int_{{{S}^{2}}}{\left\{ \frac{1}{{{r}^{2}}}{{\left| {{F }_{v\theta }} \right|}^{2}}+\frac{1}{{{r}^{2}}{{\sin }^{2}}\theta }{{\left| {{F }_{v\varphi }} \right|}^{2}}+\frac{1}{{{(1-\mu )}^{2}}}{{\left| {{F }_{vw}} \right|}^{2}} \right.}} \notag\\ 
 &+&\frac{(1-\mu )}{4{{r}^{4}}{{\sin }^{2}}\theta }{{\left| {{F }_{\varphi \theta }} \right|}^{2}}+{{\left| {{D}_{v}}\phi  \right|}^{2}}+{{\left| {{D}_{w}}\phi  \right|}^{2}}+\frac{\left( 1-\mu  \right)}{4{{r}^{2}}}{{\left| {{D}_{\theta }}\phi  \right|}^{2}} \notag\\ 
 &+&\left. \frac{\left( 1-\mu  \right)}{4{{r}^{2}}{{\sin }^{2}}\theta }{{\left| {{D}_{\varphi }}\phi  \right|}^{2}}+\frac{\left( 1-\mu  \right)}{2}P \right\}{{r}^{2}}d{{\sigma }^{2}}dv~.  
\end{eqnarray}
Considering \eqref{eH2}, the point $\left( {{w}_{i}},{{v}_{i}} \right)$, and $\left( {{v}_{i}}\le v\le {{v}_{i+1}} \right)$ are in the region $r\ge {{r}_{1}}$ with $r\left( {{w}_{i}},{{v}_{i}} \right)={{r}_{1}}$, and ${{v}_{i+1}}\ge {{v}_{i}}$. Thus, in this region
\begin{eqnarray}
    \frac{h}{\left( 1-\mu  \right)}\lesssim 1~,
\end{eqnarray}
which gives immediately \eqref{estimate1}. 
\end{proof}

\begin{preposition}\label{prep5}
    For ${{t}_{i+1}}=1.1{{t}_{i}}$, let us define the functional energy
\begin{eqnarray}
    {{\mathcal{J}}^{C}}\left( {{t}_{i}}\le t\le {{t}_{i+1}} \right)\left( {{r}_{0}}<r<{{R}_{0}} \right)\equiv\int_{{{t}_{i}}}^{{{t}_{i+1}}}{\int_{{{r}^{*}}=r_{0}^{*}}^{{{r}^{*}}=R_{0}^{*}}{\int\limits_{{{S}^{2}}}{\left[ {{\left| {{F }_{\hat{v}\hat{w}}} \right|}^{2}}+\frac{1}{4}{{\left| {{F }_{\hat{\varphi }\hat{\theta }}} \right|}^{2}} \right]\left| {{r}^{*}}-{{\left( 3m \right)}^{*}} \right|d{{r}^{*}}d{{\sigma }^{2}}dt}}}~,\notag\\
\end{eqnarray}
For $2m<{{r}_{0}}\le {{r}_{1}}<1.2{{r}_{1}}<3m$, it can be demonstrated that
\begin{eqnarray}\label{prepo5}
  &&\int_{{{r}^{*}}=-\infty }^{{{r}^{*}}=-\infty }{\int_{{{S}^{2}}}{{{r}^{2}}\left( 1-\mu  \right)d{{r}^{*}}d{{\sigma }^{2}}\left( {{\left| {{{\hat{F }}}_{\hat{t}\hat{\theta }}} \right|}^{2}}+{{\left| {{{\hat{F }}}_{\hat{t}\hat{\varphi }}} \right|}^{2}}+{{\left| {{{\hat{F }}}_{{{{\hat{r}}}^{*}}\hat{\theta }}} \right|}^{2}}+{{\left| {{{\hat{F }}}_{{{{\hat{r}}}^{*}}\hat{\varphi }}} \right|}^{2}} \right.}} \notag\\ 
 &&\left. +{{\left| {{D}_{{\hat{t}}}}\phi  \right|}^{2}}+{{\left| {{D}_{{\hat{\theta }}}}\phi  \right|}^{2}}+{{\left| {{D}_{{\hat{\varphi }}}}\phi  \right|}^{2}}+P\left( \phi  \right) \right)\cdot {{\chi }_{\left[ r_{1}^{*},\left( 1.2 \right)r_{1}^{*} \right]}} \notag\\ 
 &&\lesssim \left| {E}^{\hat{t}}\left( -0.85{{t}_{i}}\le {{r}^{*}}\le 0.85{{t}_{i}} \right)\left( t={{t}_{i}} \right) \right| ~.
\end{eqnarray}
\end{preposition}
\begin{proof}
We expect the angular momentum derivatives, or any other Killing derivatives, in the \eqref{asumeq} as in far away from horizon case to come out due to the presence of the surface on $r=3m$. Hence, using inequality (13) in \cite{Ghanem1}, we have
\begin{eqnarray}\label{es2}
    {{\mathcal{J}}^{C}}\left( {{t}_{i}}\le t\le {{t}_{i+1}} \right)\left( {{r}_{0}}<r<{{R}_{0}} \right)\lesssim\left| {{\Tilde{E}}}^{\hat{t}}\left( {{t}_{i+1}} \right) \right|+\left| {{\Tilde{E}}}^{\hat{t}}\left( {{t}_{i}} \right) \right|~.
\end{eqnarray}
We define a function $f(r^{*})$ as
    \begin{eqnarray}\label{cut}
        f\left( {{r}^{*}} \right)\equiv\int_{-\infty }^{{{r}^{*}}}{{{\chi }_{\left[ r_{1}^{*},\left( 1.2 \right)r_{1}^{*} \right]}}\left( {{r}^{*}} \right)}d{{r}^{*}}~,
    \end{eqnarray}
    where $\chi$ is the sharp cut-off function such that
\begin{equation}
f\left( {{r}^{*}} \right)=
\begin{cases}
1 ~ , for~ r_{1}^{*}<{{r}^{*}}< 1.2 r_{1}^{*}  ~ ,\\
0 ~ , for~ {{r}^{*}}<r_{1}^{*}\cup {{r}^{*}}> 1.2 r_{1}^{*} .
\end{cases}
\end{equation}
In the previous result, we have  equation \eqref{Ttime}. Then, applying  \eqref{cut} for the vector $G$ as we derive \eqref{Ttime}, one can show that
\begin{eqnarray}
  {{T}_{\alpha \beta }}{{\pi }^{\alpha \beta }}\left( G \right)&=&\left\{ \frac{1}{{{r}^{2}}}{{\left| {{F }_{w\theta }} \right|}^{2}}+\frac{1}{{{r}^{2}}{{\sin }^{2}}\theta }{{\left| {{F }_{w\varphi }} \right|}^{2}}+\frac{1}{{{r}^{2}}}{{\left| {{F }_{v\theta }} \right|}^{2}} \right. \notag\\ 
 &+&\left. \frac{1}{{{r}^{2}}{{\sin }^{2}}\theta }{{\left| {{F }_{v\varphi }} \right|}^{2}}+{{\left| {{D}_{v}}\phi  \right|}^{2}}+{{\left| {{D}_{w}}\phi  \right|}^{2}} \right\}\frac{{{\chi }_{\left[ r_{1}^{*},\left( 1.2 \right)r_{1}^{*} \right]}}}{\left( 1-\mu  \right)} \notag\\ 
 &-&2\left\{ \frac{1}{{{(1-\mu )}^{2}}}{{\left| {{F }_{vw}} \right|}^{2}}+\frac{1}{4{{r}^{4}}{{\sin }^{2}}\theta }{{\left| {{F }_{\theta \varphi }} \right|}^{2}} \right\} \notag\\ 
 &\times& \left( {{\chi }_{\left[ r_{1}^{*},\left( 1.2 \right)r_{1}^{*} \right]}}+\frac{\left( 3\mu -2 \right)}{r}\int_{-\infty }^{{{r}^{*}}}{{{\chi }_{\left[ r_{1}^{*},\left( 1.2 \right)r_{1}^{*} \right]}}\left( {{r}^{*}} \right)}d{{r}^{*}} \right) \notag\\ 
 &-&\left\{ \frac{2}{\left( 1-\mu  \right)}{{D}_{v}}\phi \overline{{{D}_{w}}\phi }+\frac{1}{2{{r}^{2}}}{{\left| {{D}_{\theta }}\phi  \right|}^{2}}+\frac{1}{2{{r}^{2}}{{\sin }^{2}}\theta }{{\left| {{D}_{\varphi }}\phi  \right|}^{2}}+P \right\} \notag\\ 
 &\times& \left( {{\chi }_{\left[ r_{1}^{*},\left( 1.2 \right)r_{1}^{*} \right]}}+\frac{\mu }{r}\int_{-\infty }^{{{r}^{*}}}{{{\chi }_{\left[ r_{1}^{*},\left( 1.2 \right)r_{1}^{*} \right]}}\left( {{r}^{*}} \right)}d{{r}^{*}} \right) \notag\\ 
 &+&\frac{\left( 1-\mu  \right)\left( 2+{{\sin }^{2}}\theta  \right)}{2r}P\int_{-\infty }^{{{r}^{*}}}{{{\chi }_{\left[ r_{1}^{*},\left( 1.2 \right)r_{1}^{*} \right]}}\left( {{r}^{*}} \right)}d{{r}^{*}}~.  
\end{eqnarray}
Applying the divergence theorem between the two hypersurfaces $\left\{ t={{t}_{i}} \right\}$ and $\left\{ t={{t}_{i+1}} \right\}$, we attain the following expression
\begin{eqnarray}\label{esIH}
  &&\int_{{{r}^{*}}=-\infty }^{{{r}^{*}}=\infty }{\int_{{{S}^{2}}}{\left\{ \frac{1}{{{r}^{2}}}{{\left| {{F }_{w\theta }} \right|}^{2}}+\frac{1}{{{r}^{2}}{{\sin }^{2}}\theta }{{\left| {{F }_{w\varphi }} \right|}^{2}}+\frac{1}{{{r}^{2}}}{{\left| {{F}_{v\theta }} \right|}^{2}} \right.}} \notag\\ 
 &&+ \left.\frac{1}{{{r}^{2}}{{\sin }^{2}}\theta }{{\left| {{F }_{v\varphi }} \right|}^{2}}+{{\left| {{D}_{v}}\phi  \right|}^{2}}+{{\left| {{D}_{w}}\phi  \right|}^{2}} \right\}{{\chi }_{\left[ r_{1}^{*},\left( 1.2 \right)r_{1}^{*} \right]}}{{r}^{2}}d{{r}^{*}}d{{\sigma }^{2}}dt \notag\\ 
 &&+\frac{\left( 1-\mu  \right)\left( 2+{{\sin }^{2}}\theta  \right)}{2r}P\int_{-\infty }^{{{r}^{*}}}{{{\chi }_{\left[ r_{1}^{*},\left( 1.2 \right)r_{1}^{*} \right]}}\left( {{r}^{*}} \right)}d{{r}^{*}} \notag\\ 
 &\lesssim& {{\mathcal{J}}^{C}}\left( {{t}_{i}}\le t\le {{t}_{i+1}} \right)\left( {{r}_{0}}<r<{{R}_{0}} \right)+{{E}^{G}}\left( {{t}_{i+1}} \right)-{{E}^{G}}\left( {{t}_{i}} \right)~.  
\end{eqnarray}
By replacing $F$ in equation \eqref{esIH} with $\hat{F}$, we can prove equation \eqref{prepo5}.
\end{proof}

\begin{preposition}\label{prep6}
    For ${{w}_{i}}={{t}_{i}}-r_{1}^{*},{{v}_{i}}={{t}_{i}}+r_{1}^{*}$, one can prove that
\begin{eqnarray}\label{es3}
  {{{\tilde{I}}}^{H}}\left( {{v}_{i}}\le v\le {{v}_{i+1}} \right)\left( {{w}_{i}}\le w\le \infty  \right)\left( r\le {{r}_{1}} \right)&\lesssim& {{\mathcal{E}}^{H}}\left( w=\infty  \right)\left( {{v}_{i}}\le v\le {{v}_{i+1}} \right) \notag\\ 
 &+&{{\mathcal{E}}^{H}}\left( v={{v}_{i+1}} \right)\left( {{w}_{i}}\le w\le \infty  \right) \notag\\ 
 &+&{{\mathcal{E}}^{\left( ^{\partial }{{/}_{\partial t}} \right)}}\left( w={{w}_{i}} \right)\left( {{v}_{i}}\le v\le {{v}_{i+1}} \right) \notag\\ 
 &+&{{\mathcal{E}}^{H}}\left( v={{v}_{i}} \right)\left( {{w}_{i}}\le w\le \infty  \right) \notag\\ 
 &+&\left| {{E}^{{\hat{t}}}}\left( -0.85{{t}_{i}}\le {{r}^{*}}\le 0.85{{t}_{i}} \right)\left( t={{t}_{i}} \right) \right|~,\notag\\ 
\end{eqnarray}
 where $\Tilde{I}^H=-I^H$. 
\end{preposition}

\begin{proof}
It can be shown that form \eqref{IH} and \eqref{estimate1}, using those equation we can prove \eqref{es3}.
\end{proof}

\begin{preposition}\label{prep7}
    For ${{v}_{i+1}}\ge {{v}_{i}}$, one can demonstrate
 \begin{eqnarray}\label{es6}
      \underset{{{v}_{i}}\le v\le {{v}_{i+1}}}{\mathop{\inf }}\,\left\{ {{{\tilde{\mathcal{E}}}}^{H}}\left( v={{v}_{i+1}} \right)\left( {{w}_{i}}\le w\le \infty  \right) \right\}&\lesssim& \frac{{{{\tilde{I}}}^{H}}\left( {{v}_{i}}\le v\le {{v}_{i+1}} \right)\left( {{w}_{i}}\le w\le \infty  \right)\left( r\le {{r}_{1}} \right)}{\left( {{v}_{i+1}}-{{v}_{i}} \right)} \notag\\ 
 &+&\underset{{{v}_{i}}\le v\le {{v}_{i+1}}}{\mathop{\sup }}\,{{\mathcal{E}}^{\left( {}^{\partial }/{}_{\partial t} \right)}}\left( v \right)\left( {{w}_{i}}\le w\le \infty  \right)\left( r\ge {{r}_{1}} \right)\notag\\
 \end{eqnarray}
\end{preposition}
\begin{proof}
For $v\geq v_i$, given the expression of $h,h'$ in the region $r\le r_1$, one can get
\begin{eqnarray}
  {{{\tilde{\mathcal{E}}}}^{H}}\left( v \right)\left( r\le {{r}_{1}} \right)&\lesssim& \int_{{{w}_{i}}}^{\infty }{\int_{{{S}^{2}}}{\left[ \left( \frac{1}{{{(1-\mu )}^{2}}}{{\left| {{F}_{vw}} \right|}^{2}}+\frac{(1-\mu )}{4{{r}^{4}}{{\sin }^{2}}\theta }{{\left| {{F}_{\varphi \theta }} \right|}^{2}} \right. \right.}}+{{D}_{v}}\phi \overline{{{D}_{w}}\phi } \notag\\ 
 &+& \left. \frac{\left( 1-\mu  \right)}{4{{r}^{2}}}{{\left| {{D}_{\theta }}\phi  \right|}^{2}}+\frac{\left( 1-\mu  \right)}{4{{r}^{2}}{{\sin }^{2}}\theta }{{\left| {{D}_{\varphi }}\phi  \right|}^{2}}+\frac{\left( 1-\mu  \right)}{2}P \right)\mu \left( \frac{1}{\left( 1-\mu  \right)}{h}'-\frac{3}{r}h \right) \notag\\ 
 &+&\frac{{{h}'}}{\left( 1-\mu  \right)}\left( \frac{1}{{{r}^{2}}}{{\left| {{F}_{v\theta }} \right|}^{2}}+\frac{1}{{{r}^{2}}{{\sin }^{2}}\theta }{{\left| {{F}_{v\varphi }} \right|}^{2}}+{{\left| {{D}_{v}}\Phi  \right|}^{2}} \right) \notag\\ 
 &+& \left.\frac{1}{{{\left( 1-\mu  \right)}^{2}}}\left( -\frac{\mu }{r}h+{h}' \right)\left( \frac{1}{{{r}^{2}}}{{\left| {{F}_{w\theta }} \right|}^{2}}+\frac{1}{{{r}^{2}}{{\sin }^{2}}\theta }{{\left| {{F}_{w\phi }} \right|}^{2}}+{{\left| {{D}_{w}}\Phi  \right|}^{2}} \right) \right]{{r}^{2}}d{{\sigma }^{2}}dw~. \notag\\ 
\end{eqnarray}
On the other hand, we also have \eqref{et}, thus, using the boundedness of $h$ and $h'$, for $r\geq r_1$, it can be shown that
\begin{eqnarray}
  {{{\tilde{\mathcal{E}}}}^{H}}\left( v \right)\left( r\ge {{r}_{1}} \right)&\lesssim& \int_{{{w}_{i}}}^{\infty }{\int_{{{S}^{2}}}{\left\{ \frac{1}{{{r}^{2}}}{{\left| {{F}_{w\theta }} \right|}^{2}}+\frac{1}{{{r}^{2}}{{\sin }^{2}}\theta }{{\left| {{F}_{v\varphi }} \right|}^{2}}+\frac{1}{(1-\mu )}{{\left| {{F}_{vw}} \right|}^{2}} \right.}} \notag\\ 
 &+&\frac{(1-\mu )}{4{{r}^{4}}{{\sin }^{2}}\theta }{{\left| {{F}_{\varphi \theta }} \right|}^{2}}+{{\left| {{D}_{v}}\phi  \right|}^{2}}+{{\left| {{D}_{w}}\phi  \right|}^{2}}+\frac{\left( 1-\mu  \right)}{4{{r}^{2}}}{{\left| {{D}_{\theta }}\phi  \right|}^{2}} \notag\\ 
 &+& \left. \frac{\left( 1-\mu  \right)}{4{{r}^{2}}{{\sin }^{2}}\theta }{{\left| {{D}_{\phi }}\phi  \right|}^{2}}+\frac{\left( 1-\mu  \right)}{2}P \right\}{{r}^{2}}d{{\sigma }^{2}}dv \notag\\ 
 &\lesssim& {{\mathcal{E}}^{\left( {}^{\partial }/{}_{\partial t} \right)}} \left( v \right)\left( r\ge {{r}_{1}} \right) ~.
\end{eqnarray}
Hence,
\begin{eqnarray}
  {{{\tilde{\mathcal{E}}}}^{H}}\left( v \right)\left( r\le {{r}_{1}} \right)&\lesssim& \int_{{{w}_{i}}}^{\infty }{\int_{{{S}^{2}}}{\left[ \left( \frac{1}{{{(1-\mu )}^{2}}}{{\left| {{F}_{vw}} \right|}^{2}}+\frac{(1-\mu )}{4{{r}^{4}}{{\sin }^{2}}\theta }{{\left| {{F}_{\varphi \theta }} \right|}^{2}} \right. \right.}}+{{D}_{v}}\phi \overline{{{D}_{w}}\phi } \notag\\ 
 &+& \left. \frac{\left( 1-\mu  \right)}{4{{r}^{2}}}{{\left| {{D}_{\theta }}\phi  \right|}^{2}}+\frac{\left( 1-\mu  \right)}{4{{r}^{2}}{{\sin }^{2}}\theta }{{\left| {{D}_{\varphi }}\phi  \right|}^{2}}+\frac{\left( 1-\mu  \right)}{2}P \right)\mu \left( \frac{1}{\left( 1-\mu  \right)}{h}'-\frac{3}{r}h \right) \notag\\ 
 &+&\frac{{{h}'}}{\left( 1-\mu  \right)}\left( \frac{1}{{{r}^{2}}}{{\left| {{F}_{v\theta }} \right|}^{2}}+\frac{1}{{{r}^{2}}{{\sin }^{2}}\theta }{{\left| {{F}_{v\varphi }} \right|}^{2}}+{{\left| {{D}_{v}}\Phi  \right|}^{2}} \right) \notag\\ 
 &+& \left.\frac{1}{{{\left( 1-\mu  \right)}^{2}}}\left( -\frac{\mu }{r}h+{h}' \right)\left( \frac{1}{{{r}^{2}}}{{\left| {{F}_{w\theta }} \right|}^{2}}+\frac{1}{{{r}^{2}}{{\sin }^{2}}\theta }{{\left| {{F}_{w\phi }} \right|}^{2}}+{{\left| {{D}_{w}}\Phi  \right|}^{2}} \right) \right]{{r}^{2}}d{{\sigma }^{2}}dw \notag\\ 
 &+& {{\mathcal{E}}^{\left( {}^{\partial }/{}_{\partial t} \right)}} \left( v \right)\left( r\ge {{r}_{1}} \right) 
\end{eqnarray}
We can therefore calculate the inequality
\begin{eqnarray}
  \left( {{v}_{i+1}}-{{v}_{i}} \right)\underset{{{v}_{i}}\le v\le {{v}_{i+1}}}{\mathop{\inf }}\,\left\{ {{{\tilde{\mathcal{E}}}}^{H}}\left( v \right)\left( {{w}_{i}}\le w\le \infty  \right) \right\}&\lesssim& \int_{{{v}_{i}}}^{{{v}_{i+1}}}{{{{\tilde{\mathcal{E}}}}^{H}}\left( v \right)\left( {{w}_{i}}\le w\le \infty  \right)dv} \notag\\ 
 & \lesssim& \frac{{{{\tilde{I}}}^{H}}\left( {{v}_{i}}\le v\le {{v}_{i+1}} \right)\left( {{w}_{i}}\le w\le \infty  \right)\left( r\le {{r}_{1}} \right)}{\left( {{v}_{i+1}}-{{v}_{i}} \right)} \notag\\ 
 & +&\underset{{{v}_{i}}\le v\le {{v}_{i+1}}}{\mathop{\sup }}\,{{\mathcal{E}}^{\left( {}^{\partial }/{}_{\partial t} \right)}}\left( v \right)\left( {{w}_{i}}\le w\le \infty  \right)\left( r\ge {{r}_{1}} \right) ~.\notag\\ 
\end{eqnarray}
\end{proof}

\begin{preposition}\label{prep8}
    For ${{w}_{i}}={{t}_{i}}-r_{1}^{*},{{v}_{i}}={{t}_{i}}+r_{1}^{*}$, we have
    \begin{eqnarray}\label{es7}
        0\le {{\tilde{I}}^{H}}\left( {{v}_{i}}\le v\le {{v}_{i+1}} \right)\left( {{w}_{i}}\le w\le \infty  \right)\left( r\le {{r}_{1}} \right)\lesssim \left| {E}^{\hat{t}} \right|+{{E}^{\#\left( {}^{\partial }/{}_{\partial t} \right)}}\left( t={{t}_{0}} \right)
    \end{eqnarray}
where
\begin{eqnarray}
  {{E}^{\#\left( {}^{\partial }/{}_{\partial t} \right)}}\left( t={{t}_{0}} \right)&=&\int_{{{r}^{*}}=-\infty }^{{{r}^{*}}=\infty }{\int_{{{S}^{2}}}{\left\{ \frac{1}{{{r}^{2}}\left( 1-\mu  \right)}{{\left| {{F }_{w\theta }} \right|}^{2}}+\frac{1}{{{r}^{2}}{{\sin }^{2}}\theta \left( 1-\mu  \right)}{{\left| {{F }_{w\varphi }} \right|}^{2}} \right.}} \notag\\ 
 &+&\frac{1}{{{r}^{2}}}{{\left| {{F }_{v\theta }} \right|}^{2}}+\frac{1}{{{r}^{2}}{{\sin }^{2}}\theta }{{\left| {{F }_{v\varphi }} \right|}^{2}}+\frac{(1-\mu )}{4{{r}^{4}}{{\sin }^{2}}\theta }{{\left| {{F }_{\varphi \theta }} \right|}^{2}}+\frac{1}{{{\left( 1-\mu  \right)}^{2}}}{{\left| {{F }_{vw}} \right|}^{2}} \notag\\ 
 &+& \left. {{\left| {{D}_{v}}\phi  \right|}^{2}}+{{\left| {{D}_{w}}\phi  \right|}^{2}}+\frac{\left( 1-\mu  \right)}{4{{r}^{2}}}{{\left| {{D}_{\theta }}\phi  \right|}^{2}}+\frac{\left( 1-\mu  \right)}{4{{r}^{2}}{{\sin }^{2}}\theta }{{\left| {{D}_{\varphi }}\phi  \right|}^{2}}+\frac{\left( 1-\mu  \right)}{2}P \right\}{{r}^{2}}d{{\sigma }^{2}}dv  \notag\\
\end{eqnarray}
\end{preposition}

\begin{proof}
First, we need to define the functional energy from vector field $H$
\begin{definition}
The functional energy $E^H\left(t\right)$ is defined by
\begin{eqnarray}\label{EHH}
    {{E}^{H}}\left( t \right)=\int_{-\infty }^{\infty }{\int_{{{S}^{2}}}{{{J}_{\alpha }}\left( H \right){{n}^{\alpha }}}}d\text{vo}{{\text{l}}_{t}}~.
\end{eqnarray}
\end{definition}
It can be derive that
\begin{eqnarray}
  {{E}^{H}}\left( t \right)&=&\int_{-\infty }^{\infty }{\int_{{{S}^{2}}}{\left\{ \frac{h}{\left( 1-\mu  \right)}\left( \frac{1}{{{r}^{2}}}{{\left| {{F}_{w\theta }} \right|}^{2}}+\frac{1}{{{r}^{2}}{{\sin }^{2}}\theta }{{\left| {{F}_{w\varphi }} \right|}^{2}}+{{\left| {{D}_{w}}\phi  \right|}^{2}} \right) \right.}} \notag\\ 
 & +&\left( \frac{h}{\left( 1-\mu  \right)}+h \right)\left( \frac{1}{{{(1-\mu )}^{2}}}{{\left| {{F}_{vw}} \right|}^{2}}+\frac{(1-\mu )}{4{{r}^{4}}{{\sin }^{2}}\theta }{{\left| {{F}_{\varphi \theta }} \right|}^{2}}+{{D}_{v}}\phi \overline{{{D}_{w}}\phi } \right. \notag\\ 
 &+& \left.\frac{\left( 1-\mu  \right)}{4{{r}^{2}}}{{\left| {{D}_{\theta }}\phi  \right|}^{2}}+\frac{\left( 1-\mu  \right)}{4{{r}^{2}}{{\sin }^{2}}\theta }{{\left| {{D}_{\varphi }}\phi  \right|}^{2}}+\frac{\left( 1-\mu  \right)}{2}P \right) \notag\\ 
 &+& \left. h\left( \frac{1}{{{r}^{2}}}{{\left| {{F}_{v\theta }} \right|}^{2}}+\frac{1}{{{r}^{2}}{{\sin }^{2}}\theta }{{\left| {{F}_{v\varphi }} \right|}^{2}}+{{\left| {{D}_{v}}\phi  \right|}^{2}} \right) \right\}{{r}^{2}}d{{\sigma }^{2}}d{{r}^{*}}~.  
\end{eqnarray}
By using the divergence theorem in the region $\left( v\le {{v}_{0}} \right)\left( {{t}_{0}}\le t\le \infty  \right)\left( r\le {{r}_{1}} \right)$, it follows that
\begin{eqnarray}\label{eq:tildeIH}
  &&{{{\tilde{I}}}^{H}}\left( v\le {{v}_{0}} \right)\left( {{t}_{0}}\le t\le \infty  \right)\left( r\le {{r}_{1}} \right)+{{{\tilde{\mathcal{E}}}}^{H}}\left( v={{v}_{0}} \right)\left( {{w}_{0}}\le w\le \infty  \right) \notag\\ 
 &+&{{{\tilde{\mathcal{E}}}}^{H}}\left( w=\infty  \right)\left( -\infty \le v\le {{v}_{0}} \right)=E^H\left( {{t}_{0}} \right)  ~.
\end{eqnarray}
Since the left hand side of \eqref{eq:tildeIH}, we get
\begin{eqnarray}
  {{{\tilde{\mathcal{E}}}}^{H}}\left( v={{v}_{0}} \right)\left( {{w}_{0}}\le w\le \infty  \right)&\lesssim& {{E}^{H}}\left( {{t}_{0}} \right) \notag\\ 
 &\lesssim& {{E}^{\#\left( {}^{\partial }/{}_{\partial t} \right)}}\left( {{t}_{0}} \right)~.  
\end{eqnarray}
It is easy to see from the divergence theorem and the fact that $\frac{\partial}{\partial t}$ is a Killing vector, by integrating in the region $\left( v\le {{v}_{0}} \right)\left( {{t}_{0}}\le t\le \infty  \right)\left( r\le {{r}_{1}} \right)$, and using the positivity of the energy \eqref{energytimelike}, we can obtain
\begin{eqnarray}\label{es5}
    {{\mathcal{E}}^{\left( {}^{\partial }/{}_{\partial t} \right)}}\left( w={{w}_{i}} \right)\left( {{v}_{i}}\le v\le {{v}_{i+1}} \right)\le \left| {{E}^{{\hat{t}}}}\left( -0.85{{t}_{i}}\le {{r}^{*}}\le 0.85{{t}_{i}} \right) \right|~.
\end{eqnarray}
Hence, applying again the inequality \eqref{es3} and using \eqref{es5}, we prove \eqref{es7}.
\end{proof}

\begin{preposition}\label{prep9}
    For all $v$, let ${{w}_{0}}\left( v \right)=v-2r_{1}^{*}$, ${{v}_{+}}=\max \left\{ 1,v \right\}$, we have the following inequality 
    \begin{eqnarray}
        {{\tilde{\mathcal{E}}}^{H}}\left( v \right)\left( {{w}_{0}}\left( v \right)\le w\le \infty  \right)\lesssim \frac{\left| {{E}^{\hat{t}}} \right|+{{E}^{\#\left( {}^{\partial }/{}_{\partial t} \right)}}\left( t={{t}_{0}} \right)+{{E}^{MH}}}{v_{+}^{2}}~,
    \end{eqnarray}
    and
    \begin{eqnarray}
        {{\tilde{\mathcal{E}}}^{H}}\left( w \right)\left( v-1\le \bar{v}\le v \right)\lesssim \frac{\left| {{E}^{\hat{t}}} \right|+{{E}^{\#\left( {}^{\partial }/{}_{\partial t} \right)}}\left( t={{t}_{0}} \right)+{{E}^{MH}}}{v_{+}^{2}}~.
    \end{eqnarray}
\end{preposition}
\begin{proof}
Remark that from preposition \ref{prep7}, using the divergence theorem and integrating in a well-chosen region, we have
\begin{eqnarray}\label{es8}
    \underset{{{v}_{i}}\le v\le {{v}_{i+1}}}{\mathop{\sup }}\,{{\mathcal{E}}^{\left( {}^{\partial }/{}_{\partial t} \right)}}\left( v \right)\left( {{w}_{i}}\le w\le \infty  \right)\left( r\ge {{r}_{1}} \right)\le \frac{{{E}^{K}}\left( {{t}_{i}} \right)}{t_{i}^{2}}
\end{eqnarray}
Thus, by using \eqref{es7} and \eqref{es8}, there always exist $v_{i}^{\#}\in \left[ {{v}_{i}},{{v}_{i+1}} \right]$ such that \eqref{es6} satisfy
\begin{eqnarray}
    \underset{{{v}_{i}}\le v\le {{v}_{i+1}}}{\mathop{\inf }}\,\left\{ {{{\tilde{\mathcal{E}}}}^{H}}\left( v \right)\left( {{w}_{i}}\le w\le \infty  \right) \right\}\lesssim \frac{{{E}^{{\hat{t}}}}+{{E}^{\#\left( {}^{\partial }/{}_{\partial t} \right)}}\left( t={{t}_{0}} \right)}{\left( {{v}_{i+1}}-{{v}_{i}} \right)^2}+\frac{{{E}^{K}}\left( {{t}_{i}} \right)}{t_{i}^{2}}~.
\end{eqnarray}
Now, let $w_{i}^{\#}=v_{i}^{\#}-2{{r}_{1}}$, then
\begin{eqnarray}\label{es9}
    \tilde{\mathcal{E}}^{H}\left( v_{i}^{\#} \right)\left( w_{i}^{\#}\le w\le \infty  \right)\lesssim\frac{{{E}^{{\hat{t}}}}+{{E}^{\#\left( {}^{\partial }/{}_{\partial t} \right)}}\left( t={{t}_{0}} \right)}{\left( t_i \right)^2}+\frac{{{E}^{K}}\left( {{t}_{i}} \right)}{t_{i}^{2}}~.
\end{eqnarray}
Then, the above preposition can be proved by substituting equation \eqref{EK} to \eqref{es9}. 
\end{proof}

\subsection{Decay for \texorpdfstring{${{F}_{\hat{v}\hat{w}}},{{F}_{{{e}_{1}}{{e}_{2}}}}$}{}}
In the case at hand, we have the Sobolev inequality
\begin{eqnarray}
  {{\left| {{F}_{\hat{v}\hat{w}}} \right|}^{2}}&\lesssim& \int_{\bar{v}=v-1}^{\bar{v}=v}{\int_{{{S}^{2}}}{\left( {{\left| {{F}_{\hat{v}\hat{w}}} \right|}^{2}}+{{\left| {{\mathcal{L}}_{v}}{{F}_{\hat{v}\hat{w}}} \right|}^{2}}+{{\left| \slashed{\mathcal{L}}{{F}_{\hat{v}\hat{w}}} \right|}^{2}}+{{\left| \slashed{\mathcal{L}}{{\mathcal{L}}_{v}}{{F}_{\hat{v}\hat{w}}} \right|}^{2}} \right.}} \notag\\ 
 && + \left. {{\left| \slashed{\mathcal{L}}\slashed{\mathcal{L}}{{F}_{\hat{v}\hat{w}}} \right|}^{2}}+{{\left| \slashed{\mathcal{L}}\slashed{\mathcal{L}}{{\mathcal{L}}_{v}}{{F}_{\hat{v}\hat{w}}} \right|}^{2}} \right)d{{\sigma }^{2}}d\bar{v}~.  
\end{eqnarray}
First, we will demonstrate that the following inequality holds
\begin{eqnarray}\label{Fnear1}
\int_{\bar{v}=v-1}^{\bar{v}=v}{\int_{{{S}^{2}}}{{{\left| {{\mathcal{L}}_{v}}{{F }_{\hat{v}\hat{w}}} \right|}^{2}}}d{{\sigma }^{2}}d\bar{v}}&\lesssim& \int_{\bar{v}=v-1}^{\bar{v}=v}{\int_{{{S}^{2}}}{{{\left| {{\nabla }_{{{{\hat{e}}}_{a}}}}{{F }_{\hat{v}{{{\hat{e}}}_{a}}}}+{{\nabla }_{{{{\hat{e}}}_{b}}}}{{F }_{\hat{v}{{{\hat{e}}}_{b}}}}+{{D}_{{\hat{v}}}}\phi +\phi  \right|}^{2}}d{{\sigma }^{2}}d\bar{v}}} \notag\\ 
 &\lesssim& \int_{\bar{v}=v-1}^{\bar{v}=v}{\int_{{{S}^{2}}}{\left\{ {{\left| \slashed{\mathcal{L}}{{F }_{\hat{v}{{{\hat{e}}}_{a}}}} \right|}^{2}}+{{\left| \slashed{\mathcal{L}}{{F }_{\hat{v}{{{\hat{e}}}_{b}}}} \right|}^{2}}+{{\left| {{F }_{\hat{v}{{{\hat{e}}}_{a}}}} \right|}^{2}} \right.}} \notag\\ 
 && + \left. {{\left| {{F }_{\hat{v}{{{\hat{e}}}_{b}}}} \right|}^{2}}+{{\left| {{F }_{\hat{v}\hat{w}}} \right|}^{2}}+{{\left| {{D}_{{\hat{v}}}}\phi  \right|}^{2}}+{{\left| \phi  \right|}^{2}} \right\}d{{\sigma }^{2}}d\bar{v} ~. 
\end{eqnarray}
To estimate \eqref{Fnear1}, we need to estimate the term $\int_{\bar{v}=v-1}^{\bar{v}=v}{\int_{{{S}^{2}}}{{{\left| \phi  \right|}^{2}}}d{{\sigma }^{2}}d\bar{v}}$.

Using Holder inequality, we see that
\begin{eqnarray}
  \frac{\partial }{\partial w}\left\| \phi  \right\|_{{{L}^{2}}}^{2}&=&\int_{\bar{v}=v-1}^{\bar{v}=v}{\int_{{{S}^{2}}}{\frac{\partial }{\partial w}{{\left| \phi  \right|}^{2}}}d{{\sigma }^{2}}d\bar{v}} \notag\\ 
 2{{\left\| \phi  \right\|}_{{{L}^{2}}}}\frac{\partial {{\left\| \phi  \right\|}_{{{L}^{2}}}}}{\partial w}&=&\int_{\bar{v}=v-1}^{\bar{v}=v}{\int_{{{S}^{2}}}{2\phi {{\partial }_{w}}}\phi d{{\sigma }^{2}}d\bar{v}} \notag\\ 
 &\lesssim& 2{{\left\| \phi  \right\|}_{{{L}^{2}}}}{{\left( \int_{\bar{v}=v-1}^{\bar{v}=v}{\int_{{{S}^{2}}}{\left| {{\partial }_{w}}\phi  \right|}d{{\sigma }^{2}}d\bar{v}} \right)}^{1/2}} \notag\\ 
\frac{\partial {{\left\| \phi  \right\|}_{{{L}^{2}}}}}{\partial w}&\lesssim& {{\left( \int_{\bar{v}=v-1}^{\bar{v}=v}{\int_{{{S}^{2}}}{\left| {{\partial }_{w}}\phi  \right|}d{{\sigma }^{2}}d\bar{v}} \right)}^{1/2}} \notag\\ 
 &\lesssim& {{\left\{ {{{\tilde{\mathcal{E}}}}^{H}}\left( w \right)\left( v-1\le \bar{v}\le v \right) \right\}}^{1/2}}~,  
\end{eqnarray}
where 
\begin{eqnarray}
    \left\| \phi  \right\|_{{{L}^{2}}}^{2}=\int_{\bar{v}=v-1}^{\bar{v}=v}{\int_{{{S}^{2}}}{{{\left| \phi  \right|}^{2}}}d{{\sigma }^{2}}d\bar{v}}~.
\end{eqnarray}
It can be shown that
\begin{eqnarray}\label{scalares}
    \int_{\bar{v}=v-1}^{\bar{v}=v}{\int_{{{S}^{2}}}{{{\left| \phi  \right|}^{2}}}d{{\sigma }^{2}}d\bar{v}}\lesssim {{\left( \frac{w}{{{v}_{+}}} \right)}^{2}}{{E}_{1}}~.
\end{eqnarray}
Thus, we can write down the estimate for \eqref{Fnear1} as
\begin{eqnarray}
    \int_{\bar{v}=v-1}^{\bar{v}=v}{\int_{{{S}^{2}}}{{{\left| {{\mathcal{L}}_{v}}{{F }_{\hat{v}\hat{w}}} \right|}^{2}}}d{{\sigma }^{2}}d\bar{v}}\lesssim {{\left( \frac{w}{{{v}_{+}}} \right)}^{2}}{{E}_{1}}
\end{eqnarray}
where
\begin{eqnarray}\label{E1}
    {{E}_{1}}={{\left[ \sum\limits_{j=0}^{6}{E_{{{r}^{j}}{{\left( \slashed{\mathcal{L}} \right)}^{j}}F}^{\left( {}^{\partial }/{}_{\partial t} \right)}\left( t={{t}_{0}} \right)}+\sum\limits_{j=0}^{5}{E_{{{r}^{j}}{{\left( \slashed{\mathcal{L}} \right)}^{j}}F}^{\left( K \right)}\left( t={{t}_{0}} \right)}+\sum\limits_{j=1}^{3}{E_{{{r}^{j}}{{\left( \slashed{\mathcal{L}} \right)}^{j}}F}^{\#\left( {}^{\partial }/{}_{\partial t} \right)}\left( t={{t}_{0}} \right)} \right]}^{1/2}}
\end{eqnarray}
Therefore, we have
\begin{eqnarray}
    {{\left| {{F}_{\hat{v}\hat{w}}} \right|}^{2}}\lesssim {{\left( \frac{w}{{{v}_{+}}} \right)}^{2}}{{E}_{1}}~.
\end{eqnarray}
Focusing on the component ${{F}_{{{e}_{1}}{{e}_{2}}}}$, we also have the Sobolev inequality,
\begin{eqnarray}
  {{\left| {{F}_{{{e}_{1}}{{e}_{2}}}} \right|}^{2}}&\lesssim& \int_{\bar{v}=v-1}^{\bar{v}=v}{\int_{{{S}^{2}}}{\left( {{\left| {{F}_{{{e}_{1}}{{e}_{2}}}} \right|}^{2}}+{{\left| {{\mathcal{L}}_{v}}{{F}_{{{e}_{1}}{{e}_{2}}}} \right|}^{2}}+{{\left| \slashed{\mathcal{L}}{{F}_{{{e}_{1}}{{e}_{2}}}} \right|}^{2}}+{{\left| \slashed{\mathcal{L}}{{\mathcal{L}}_{v}}{{F}_{{{e}_{1}}{{e}_{2}}}} \right|}^{2}} \right.}} \notag\\ 
 &+&\left.{{\left| \slashed{\mathcal{L}}\slashed{\mathcal{L}}{{F}_{{{e}_{1}}{{e}_{2}}}} \right|}^{2}}+{{\left| \slashed{\mathcal{L}}\slashed{\mathcal{L}}{{\mathcal{L}}_{v}}{{F}_{{{e}_{1}}{{e}_{2}}}} \right|}^{2}} \right)d{{\sigma }^{2}}d\bar{v}~.  
\end{eqnarray}
Computing 
\begin{eqnarray}
  \int_{\bar{v}=v-1}^{\bar{v}=v}{\int_{{{S}^{2}}}{{{\left| {{\mathcal{L}}_{v}}{{F }_{{{e}_{1}}{{e}_{2}}}} \right|}^{2}}}d{{\sigma }^{2}}d\bar{v}}&\lesssim& \int_{\bar{v}=v-1}^{\bar{v}=v}{\int_{{{S}^{2}}}{\left\{ {{\left| {{\mathcal{L}}_{{{e}_{1}}}}{{F }_{\hat{v}{{e}_{2}}}} \right|}^{2}}+{{\left| {{\mathcal{L}}_{{{e}_{2}}}}{{F }_{\hat{v}{{e}_{1}}}} \right|}^{2}}+{{\left| {{F }_{{{e}_{1}}{{e}_{2}}}} \right|}^{2}} \right.}} \notag\\ 
 &+&\left. {{\left| {{F }_{\hat{v}{{e}_{2}}}} \right|}^{2}}+{{\left| {{D}_{{\hat{v}}}}\phi  \right|}^{2}}+{{\left| \phi  \right|}^{2}} \right\}d{{\sigma }^{2}}d\bar{v} \notag\\ 
 &\lesssim& {{\left( \frac{w}{{{v}_{+}}} \right)}^{2}}{{E}_{1}}~.  
\end{eqnarray}
it follows that 
\begin{eqnarray}
    {{\left| {{F}_{{{e}_{1}}{{e}_{2}}}} \right|}^{2}}\lesssim {{\left( \frac{w}{{{v}_{+}}} \right)}^{2}}{{E}_{1}}~.
\end{eqnarray}

\subsection{Decay for \texorpdfstring{${{F}_{\hat{v}{{e}_{1}}}},{{F}_{\hat{v}{{e}_{2}}}}$}{}}
We have the Sobolev inequality
\begin{eqnarray}\label{eq:Fvteta}
  {{\left| {{F}_{\hat{v}\theta }} \right|}^{2}}&\lesssim& \int_{\bar{v}=v-1}^{\bar{v}=v}{\int_{{{S}^{2}}}{\left( {{\left| {{F}_{\hat{v}\theta }} \right|}^{2}} + {{\left| {{\mathcal{L}}_{v}}{{F}_{\hat{v}\theta }} \right|}^{2}} +{{\left| \slashed{\mathcal{L}}{{F}_{\hat{v}\theta }} \right|}^{2}}+{{\left| \slashed{\mathcal{L}}{{\mathcal{L}}_{v}}{{F}_{\hat{v}\theta }} \right|}^{2}} \right.}} \notag\\ 
 &+&\left.{{\left| \slashed{\mathcal{L}}\slashed{\mathcal{L}}{{F}_{\hat{v}\theta }} \right|}^{2}}+{{\left| \slashed{\mathcal{L}}\slashed{\mathcal{L}}{{\mathcal{L}}_{v}}{{F}_{\hat{v}\theta }} \right|}^{2}} \right)d{{\sigma }^{2}}d\bar{v} ~ . 
\end{eqnarray}
Here, we want to estimate the fourth term of the right hand side in \eqref{eq:Fvteta}. First, using the equation of motion \eqref{eom1}, we show that
\begin{eqnarray}
{{\mathcal{L}}_{v}}{{F }_{\hat{v}\theta }}&=&{{\nabla }_{t}}{{F }_{\hat{v}\hat{\theta }}}-\frac{1}{4}\left( {{\nabla }_{{\hat{\theta }}}}{{F}_{r*t}}+\left( 1-\mu  \right){{\nabla }_{{\hat{\varphi }}}}{{F }_{\hat{\varphi }\hat{\theta }}} \right) \notag\\ 
 &+&\frac{1}{4}i\left( -{{D}_{{\hat{\theta }}}}\phi \text{ }\bar{\phi }+\phi \overline{{{D}_{{\hat{\theta }}}}\phi } \right)+\frac{\mu }{2r}{{F }_{\hat{v}\theta }}~.
\end{eqnarray}
Then, applying the Bianchi identities 
\begin{eqnarray}
    {{\nabla }_{\lambda }}{{F }_{\alpha \beta }}+{{\nabla }_{\alpha }}{{F }_{\beta \lambda }}+{{\nabla }_{\beta }}{{F }_{\lambda \alpha }}=0~,
\end{eqnarray}
we also have
\begin{eqnarray}
  {{\mathcal{L}}_{t}}{{F }_{\hat{v}\hat{\theta }}}&=&{{\nabla }_{t}}{{F }_{\hat{v}\hat{\theta }}}+F \left( {{\nabla }_{t}}\hat{v},\hat{\theta } \right)+F\left( \hat{v},{{\nabla }_{t}}\hat{\theta } \right) \notag\\ 
 {{\mathcal{L}}_{{\hat{\varphi }}}}{{F }_{\hat{\varphi }\hat{\theta }}}&=&{{\nabla }_{{\hat{\varphi }}}}{{F }_{\hat{\varphi }\hat{\theta }}}+F \left( {{\nabla }_{{\hat{\varphi }}}}\hat{\varphi },\hat{\theta } \right)+F \left( \hat{\varphi },{{\nabla }_{{\hat{\varphi }}}}\hat{\theta } \right) \notag\\
 {{\mathcal{L}}_{{\hat{\theta }}}}{{F }_{r*t}}&=&{{\nabla }_{{\hat{\theta }}}}{{F }_{r*t}}+F \left( {{\nabla }_{{\hat{\theta }}}}{{r}^{*}},t \right)+F \left( {{r}^{*}},{{\nabla }_{{\hat{\theta }}}}t \right) ~.
\end{eqnarray}
So, one can get
\begin{eqnarray}
{{\mathcal{L}}_{v}}{{F }_{\hat{v}\theta }}&=&{{\mathcal{L}}_{t}}{{F }_{\hat{v}\hat{\theta }}}-\frac{1}{2}\left( 1-\mu  \right){{\mathcal{L}}_{{\hat{\theta }}}}{{F }_{\hat{v}\hat{w}}}-\frac{1}{4}\left( 1-\mu  \right){{\mathcal{L}}_{{\hat{\varphi }}}}{{F }_{\hat{\varphi }\hat{\theta }}} \notag\\ 
 &-&\frac{1}{2}\frac{\left( 1-\mu  \right)}{r}{{F }_{\hat{v}\hat{\theta }}}-\frac{1}{4}i\left( -{{D}_{{\hat{\theta }}}}\phi \text{ }\bar{\phi }+\phi \overline{{{D}_{{\hat{\theta }}}}\phi } \right)  ~.
\end{eqnarray}
Similarly, we obtain
\begin{eqnarray}
  {{\mathcal{L}}_{v}}{{F }_{\hat{v}\hat{\varphi }}}&=&{{\mathcal{L}}_{t}}{{F }_{\hat{v}\hat{\varphi }}}-\frac{1}{2}\left( 1-\mu  \right){{\mathcal{L}}_{{\hat{\varphi }}}}{{F }_{\hat{v}\hat{w}}}-\frac{1}{4}\left( 1-\mu  \right){{\mathcal{L}}_{{\hat{\theta }}}}{{F }_{\hat{\theta }\hat{\varphi }}} \notag\\ 
 &-&\frac{1}{2}\frac{\left( 1-\mu  \right)}{r}{{F }_{\hat{v}\hat{\varphi }}}-\frac{1}{4}i\left( -{{D}_{{\hat{\varphi }}}}\phi \text{ }\bar{\phi }+\phi \overline{{{D}_{{\hat{\varphi }}}}\phi } \right)  ~.
\end{eqnarray}
Thus, for $r\ge 2m>0$
\begin{eqnarray}
    {{\left| {{\mathcal{L}}_{v}}{{F }_{\hat{v}\theta }} \right|}^{2}}\lesssim {{\left| {{\mathcal{L}}_{t}}{{F }_{\hat{v}\hat{\theta }}} \right|}^{2}}+{{\left| {{\mathcal{L}}_{{\hat{\theta }}}}{{F }_{\hat{v}\hat{w}}} \right|}^{2}}+{{\left| {{\mathcal{L}}_{{\hat{\varphi }}}}{{F }_{\hat{\varphi }\hat{\theta }}} \right|}^{2}}+{{\left| {{F }_{\hat{v}\hat{\theta }}} \right|}^{2}}+{{\left| {{D}_{{\hat{\theta }}}}\phi  \right|}^{2}}+{{\left| \phi  \right|}^{2}}~.
\end{eqnarray}
Finally, after some computation, we have
\begin{eqnarray}
    \left| {{F}_{\hat{v}\hat{\theta }}} \right|\lesssim {{\left( \frac{w}{{{v}_{+}}} \right)}^{2}}{{E}_{2}}~,
\end{eqnarray}
\begin{eqnarray}
    \left| {{F}_{\hat{v}\hat{\varphi }}} \right|\lesssim {{\left( \frac{w}{{{v}_{+}}} \right)}^{2}}{{E}_{2}}~,
\end{eqnarray}
where
\begin{eqnarray}\label{E2}
    {{E}_{2}}={{\left[ E_{F}^{2}+\sum\limits_{i=0}^{1}{\sum\limits_{j=1}^{2}{E_{{{r}^{j}}{{\left( \slashed{\mathcal{L}} \right)}^{j}}F{{\left( {{\mathcal{L}}_{t}} \right)}^{i}}F}^{\#\left( {}^{\partial }/{}_{\partial t} \right)}\left( t={{t}_{0}} \right)}}+E_{{{r}^{3}}{{\left( \slashed{\mathcal{L}} \right)}^{3}}F}^{\#\left( {}^{\partial }/{}_{\partial t} \right)}\left( t={{t}_{0}} \right) \right]}^{1/2}}~.
\end{eqnarray}

\subsection{Decay for \texorpdfstring{$\sqrt{1-\frac{2m}{r}}{{F}_{\hat{w}{{e}_{1}}}},\sqrt{1-\frac{2m}{r}}{{F}_{\hat{w}{{e}_{2}}}}$}{}}
As in the preceding section, we also have
\begin{eqnarray}
  {{\left| \sqrt{1-\mu }{{F}_{\hat{w}\hat{\theta }}} \right|}^{2}}&\lesssim& \int_{\bar{w}={{w}_{0}}}^{\bar{w}=\infty }{\int_{{{S}^{2}}}{\left( {{\left| \sqrt{1-\mu }{{F}_{\hat{w}\hat{\theta }}} \right|}^{2}}+{{\left| {{\mathcal{L}}_{w}}\sqrt{1-\mu }{{F}_{\hat{w}\hat{\theta }}} \right|}^{2}}+{{\left| \slashed{\mathcal{L}}{{\mathcal{L}}_{w}}\left( \sqrt{1-\mu }{{F}_{\hat{w}\hat{\theta }}} \right) \right|}^{2}}\right.}} \notag\\ 
 &+& \left. {{\left| \slashed{\mathcal{L}}\sqrt{1-\mu }{{F}_{\hat{w}\hat{\theta }}} \right|}^{2}}+{{\left| \slashed{\mathcal{L}}\slashed{\mathcal{L}}\left( \sqrt{1-\mu }{{F}_{\hat{w}\hat{\theta }}} \right) \right|}^{2}}+{{\left| \slashed{\mathcal{L}}\slashed{\mathcal{L}}{{\mathcal{L}}_{w}}\left( \sqrt{1-\mu }{{F}_{\hat{w}\hat{\theta }}} \right) \right|}^{2}} \right)d{{\sigma }^{2}}d\bar{w}\notag\\  
\end{eqnarray}
Applying the same procedure, after some computation, it follows that
\begin{eqnarray}
    {{\left| \sqrt{1-\mu }{{F}_{\hat{w}\hat{\theta }}} \right|}^{2}}\lesssim {{\left( \frac{w}{{{v}_{+}}} \right)}^{2}}{{E}_{2}}~,
\end{eqnarray}
\begin{eqnarray}
    {{\left| \sqrt{1-\mu }{{F}_{\hat{w}\hat{\varphi }}} \right|}^{2}}\lesssim {{\left( \frac{w}{{{v}_{+}}} \right)}^{2}}{{E}_{2}}~.
\end{eqnarray}

\subsection{Decay estimate for \texorpdfstring{$\phi$}{}}
Let us consider again the Sobolev inequality for the scalar field
\begin{eqnarray}
{{\left| \phi  \right|}^{2}}\lesssim \int_{\bar{v}=v-1}^{\bar{v}=v}{\int_{{{S}^{2}}}{\left( {{\left| \phi  \right|}^{2}}+{{\left| {{\mathcal{L}}_{v}}\phi  \right|}^{2}}+{{\left| \slashed{\mathcal{L}}\phi  \right|}^{2}}+{{\left| \slashed{\mathcal{L}}{{\mathcal{L}}_{v}}\phi  \right|}^{2}}+{{\left| \slashed{\mathcal{L}}\slashed{\mathcal{L}}\phi  \right|}^{2}}+{{\left| \slashed{\mathcal{L}}\slashed{\mathcal{L}}{{\mathcal{L}}_{v}}\phi  \right|}^{2}} \right)}}d{{\sigma }^{2}}d\bar{v} ~ .
\end{eqnarray}
Considering \eqref{eH2}, we can get
\begin{eqnarray}
      \tilde{\mathcal{E}}_{{{\mathcal{L}}_{{{\Omega }_{j}}}}}^{H}&=&2\int_{\bar{v}=v-1}^{\bar{v}=v}{\int_{{{S}^{2}}}{\left[ h\left( \frac{1}{{{(1-\mu )}^{2}}}{{\left| \slashed{\mathcal{L}}{{\Psi }_{vw}} \right|}^{2}}+(1-\mu ){{\left| \slashed{\mathcal{L}}{{\Psi }_{{{e}_{1}}{{e}_{2}}}} \right|}^{2}} \right. \right.}} \notag\\ 
 &+&\slashed{\mathcal{L}}\left( {{D}_{v}}\Phi \overline{{{D}_{w}}\Phi } \right)+\frac{\left( 1-\mu  \right)}{4}{{\left| \slashed{\mathcal{L}}\left( {{D}_{{{e}_{1}}}}\Phi  \right) \right|}^{2}} \notag\\ 
 &+& \left. \frac{\left( 1-\mu  \right)}{4}{{\left| \slashed{\mathcal{L}}\left( {{D}_{{{e}_{2}}}}\Phi  \right) \right|}^{2}}+\frac{\left( 1-\mu  \right)}{2}\slashed{\mathcal{L}}\mathcal{P} \right) \notag\\ 
 &+&\left. \frac{h}{\left( 1-\mu  \right)}\left( {{\left| \slashed{\mathcal{L}}{{\Psi }_{v{{e}_{1}}}} \right|}^{2}}+{{\left| \slashed{\mathcal{L}}{{\Psi }_{v{{e}_{2}}}} \right|}^{2}}+{{\left| \slashed{\mathcal{L}}\left( {{D}_{v}}\Phi  \right) \right|}^{2}} \right) \right]{{r}^{2}}d{{\sigma }^{2}}d\bar{v}  ~,
\end{eqnarray}
\begin{eqnarray}
\tilde{\mathcal{E}}_{{{\mathcal{L}}_{{{\Omega }_{j}}}}{{\mathcal{L}}_{{{\Omega }_{k}}}}}^{H}&=&2\int_{\bar{v}=v-1}^{\bar{v}=v}{\int_{{{S}^{2}}}{\left[ h\left( \frac{1}{{{(1-\mu )}^{2}}}{{\left| \slashed{\mathcal{L}}\slashed{\mathcal{L}}{{\Psi }_{vw}} \right|}^{2}}+(1-\mu ){{\left| \slashed{\mathcal{L}}\slashed{\mathcal{L}}{{\Psi }_{{{e}_{1}}{{e}_{2}}}} \right|}^{2}} \right. \right.}} \notag\\ 
 & +&\slashed{\mathcal{L}}\slashed{\mathcal{L}}\left( {{D}_{v}}\Phi \overline{{{D}_{w}}\Phi } \right)+\frac{\left( 1-\mu  \right)}{4}{{\left| \slashed{\mathcal{L}}\slashed{\mathcal{L}}\left( {{D}_{{{e}_{1}}}}\Phi  \right) \right|}^{2}} \notag\\ 
 &+& \left.\frac{\left( 1-\mu  \right)}{4}{{\left| \slashed{\mathcal{L}}\slashed{\mathcal{L}}\left( {{D}_{{{e}_{2}}}}\Phi  \right) \right|}^{2}}+\frac{\left( 1-\mu  \right)}{2}\mathcal{L}\mathcal{P} \right) \notag\\ 
 &+&\left. \frac{h}{\left( 1-\mu  \right)}\left( {{\left| \slashed{\mathcal{L}}\slashed{\mathcal{L}}{{\Psi }_{v{{e}_{1}}}} \right|}^{2}}+{{\left| \slashed{\mathcal{L}}\slashed{\mathcal{L}}{{\Psi }_{v{{e}_{2}}}} \right|}^{2}}+{{\left| \slashed{\mathcal{L}}\slashed{\mathcal{L}}\left( {{D}_{v}}\Phi  \right) \right|}^{2}} \right) \right]{{r}^{2}}d{{\sigma }^{2}}d\bar{v}~. \notag\\ 
\end{eqnarray}
Suppose that
\begin{eqnarray}\label{E3}
    {{E}_{3}}={{\left[ \left| E_{_{,{{\mathcal{L}}_{{{\Omega }_{j}}}},{{\mathcal{L}}_{{{\Omega }_{j}}}}{{\mathcal{L}}_{{{\Omega }_{k}}}}}}^{\hat{t}} \right|+E_{_{,{{\mathcal{L}}_{{{\Omega }_{j}}}},{{\mathcal{L}}_{{{\Omega }_{j}}}}{{\mathcal{L}}_{{{\Omega }_{k}}}}}}^{\#\left( {}^{\partial }/{}_{\partial t} \right)}\left( t={{t}_{0}} \right)+{{E}^{MH}} \right]}^{1/2}}~,
\end{eqnarray}
from \eqref{scalares}, we have
\begin{eqnarray}
    \int_{\bar{v}=v-1}^{\bar{v}=v}{\int_{{{S}^{2}}}{\left( {{\left| \phi  \right|}^{2}}+{{\left| \slashed{\nabla }\phi  \right|}^{2}}+{{\left| \slashed{\nabla }\slashed{\nabla }\phi  \right|}^{2}} \right)}}d{{\sigma }^{2}}d\bar{v}\lesssim {{\left( \frac{w}{{{v}_{+}}} \right)}^{2}}{{E}_{3}}~,
\end{eqnarray}
\begin{eqnarray}
    \int_{\bar{v}=v-1}^{\bar{v}=v}{\int_{{{S}^{2}}}{\left( {{\left| {{\nabla }_{v}}\phi  \right|}^{2}}+{{\left| \slashed{\nabla }{{\nabla }_{v}}\phi  \right|}^{2}}+{{\left| \slashed{\nabla }\slashed{\nabla }{{\nabla }_{v}}\phi  \right|}^{2}} \right)}}d{{\sigma }^{2}}d\bar{v}\lesssim {{E}_{3}}~.
\end{eqnarray}
Finally, one can have
\begin{eqnarray}\label{nearphi}
    \left| \phi  \right|\lesssim {{\left( 1+{{\left( \frac{w}{{{v}_{+}}} \right)}^{2}} \right)}^{1/2}}{{E}_{3}}^{1/2}~.
\end{eqnarray}

\subsection{Decay Estimate for \texorpdfstring{$D\phi$}{} near horizon}
Now, focusing on the estimate for $\left| D\phi  \right|$, where
\begin{eqnarray}\label{Depi}
    \left| D\phi  \right|\le \left| \nabla \phi  \right|+\left| A \right|\left| \phi  \right|~.
\end{eqnarray}
we have the Sobolev inequality
\begin{eqnarray}\label{eq:nablaphi2}
  {{\left| \nabla \phi  \right|}^{2}}&\lesssim& \int_{\bar{v}=v-1}^{\bar{v}=v}{\int_{{{S}^{2}}}{\left\{ {{\left| \nabla \phi  \right|}^{2}}+{{\left| \slashed{\nabla }\nabla \phi  \right|}^{2}}+{{\left| \slashed{\nabla }\slashed{\nabla }\nabla \phi  \right|}^{2}}+ {{\left| {{\nabla }_{v}}\left( \nabla \phi  \right) \right|}^{2}} \right.}} \notag\\ 
 &+&\left. {{\left| \slashed{\nabla }{{\nabla }_{v}}\left( \nabla \phi  \right) \right|}^{2}}+{{\left| \slashed{\nabla }\slashed{\nabla }{{\nabla }_{v}}\left( \nabla \phi  \right) \right|}^{2}} \right\}d{{\sigma }^{2}}d\bar{v}~.
\end{eqnarray}
The estimate of the fourth term of the right hand side in \eqref{eq:nablaphi2} has the form 
\begin{eqnarray}\label{poteq}
      \int_{\bar{v}=v-1}^{\bar{v}=v}{\int_{{{S}^{2}}}{{{\left| {{\nabla }_{v}}\left( \nabla \phi  \right) \right|}^{2}}}}d{{\sigma }^{2}}d\bar{v}&\lesssim& \int_{\bar{v}=v-1}^{\bar{v}=v}{\int_{{{S}^{2}}}{{{\left| \nabla \nabla \phi  \right|}^{2}}}}d{{\sigma }^{2}}d\bar{v} \notag\\ 
 &\lesssim& \int_{\bar{v}=v-1}^{\bar{v}=v}{\int_{{{S}^{2}}}{{{\left| {{\nabla }_{\phi }}P \right|}^{2}}}}d{{\sigma }^{2}}d\bar{v}~.  
\end{eqnarray}
As in the far away from horizon case, in order to have an estimate of \eqref{poteq}, we have to specify the form of the scalar potential. We also consider the form as in assumption \ref{asumsiP}.
Then, we have
\begin{equation}
\hat{P}=\int_{\bar{v}=v-1}^{\bar{v}=v}{\int_{{{S}^{2}}}{{{\left| {{\nabla }_{\phi }}\mathcal{P} \right|}^{2}}}}d{{\sigma }^{2}}d\bar{v}\lesssim 
\begin{cases}
E_{3}^{3} ~ , &\text{if $P$ is of the form \eqref{pi4}} ~ ,\\
E_3 + 1 ~ , &\text{if $P$ is is either of the form \eqref{mass},\eqref{sine},\eqref{toda}}.
\end{cases}
\end{equation}

Next, we need to estimate
\begin{eqnarray}
  \int_{\bar{v}=v-1}^{\bar{v}=v}{\int_{{{S}^{2}}}{{{\left| \slashed{\nabla }{{\nabla }_{v}}\left( \nabla \phi  \right) \right|}^{2}}}}d{{\sigma }^{2}}d\bar{v}&\lesssim& \int_{\bar{v}=v-1}^{\bar{v}=v}{\int_{{{S}^{2}}}{{{\left| \nabla \nabla \left( \nabla \phi  \right) \right|}^{2}}}}d{{\sigma }^{2}}d\bar{v}~. \notag\\ 
 &\lesssim& \int_{\bar{v}=v-1}^{\bar{v}=v}{\int_{{{S}^{2}}}{{{\left| \nabla \left( {{\nabla }_{\phi }}P \right) \right|}^{2}}}}d{{\sigma }^{2}}d\bar{v}  
\end{eqnarray}
Considering the potential form as in assumption \ref{asumsiP}, it follows that 
\begin{equation}
\int_{\bar{v}=v-1}^{\bar{v}=v}{\int_{{{S}^{2}}}{{{\left| \nabla \left( {{\nabla }_{\phi }}\mathcal{P} \right) \right|}^{2}}}}d{{\sigma }^{2}}d\bar{v}\lesssim 
\begin{cases}
E_{3}^{4} ~ , &\text{if $P$ is of the form \eqref{pi4}} ~ ,\\
E_3 + 1 ~ , &\text{if $P$ is is either of the form \eqref{mass},\eqref{sine},\eqref{toda}}.
\end{cases}
\end{equation}
Similarly, we have
\begin{eqnarray}
      \int_{\bar{v}=v-1}^{\bar{v}=v}{\int_{{{S}^{2}}}{{{\left| \slashed{\nabla }\slashed{\nabla }{{\nabla }_{v}}\left( \nabla \phi  \right) \right|}^{2}}}}d{{\sigma }^{2}}d\bar{v}&\lesssim& \int_{\bar{v}=v-1}^{\bar{v}=v}{\int_{{{S}^{2}}}{{{\left| \nabla \nabla \nabla \nabla \phi  \right|}^{2}}}}d{{\sigma }^{2}}d\bar{v} \notag\\ 
 &\lesssim& \int_{\bar{v}=v-1}^{\bar{v}=v}{\int_{{{S}^{2}}}{{{\left| \nabla \nabla \left( {{\nabla }_{\phi }}\mathcal{P} \right) \right|}^{2}}}}d{{\sigma }^{2}}d\bar{v}~.  
\end{eqnarray}
Then, one can get
\begin{equation}
\int_{\bar{v}=v-1}^{\bar{v}=v}{\int_{{{S}^{2}}}{{{\left| \nabla \nabla \left( {{\nabla }_{\phi }}\mathcal{P} \right) \right|}^{2}}}}d{{\sigma }^{2}}d\bar{v}\lesssim 
\begin{cases}
E_{3}^{7}+E_3 ~ , &\text{if $P$ is of the form \eqref{pi4}} ~ ,\\
E_3 + 1 ~ , &\text{if $P$ is is either of the form \eqref{mass},\eqref{sine}, or \eqref{toda}}.
\end{cases}
\end{equation}
Let us define
\begin{equation}\label{E4}
 E_4\equiv
\begin{cases}
{{\left( \frac{w}{{{v}_{+}}} \right)}^{2}}{{E}_{3}}+E_{3}^{7}+E_{3}^{4}+E_{3}^{3}+{{E}_{3}} ~ , &\text{if $P$ is of the form \eqref{pi4}} ~ ,\\
 {{\left( \frac{w}{{{v}_{+}}} \right)}^{2}}{{E}_{3}}+{{E}_{3}}+1 ~ , &\text{if $P$ is is either of the form \eqref{mass},\eqref{sine}, or \eqref{toda}}.
\end{cases}
\end{equation}
Then, we can demonstrate
\begin{eqnarray}\label{nabnearphi}
    \left| \nabla \phi  \right|\lesssim E_{4}^{1/2}~.
\end{eqnarray}
To get the complete estimate of \eqref{Depi}, one need to estimate $\left| A \right|$. As in the preceding section, we have the Sobolev inequality
\begin{eqnarray}
      {{\left| A \right|}^{2}}&\lesssim &\int_{\bar{v}=v-1}^{\bar{v}=v}{\int_{{{S}^{2}}}{\left( {{\left| A \right|}^{2}}+{{\left| {{\mathcal{L}}_{v}}A \right|}^{2}}+{{\left| \slashed{\mathcal{L}}A \right|}^{2}}+{{\left| \slashed{\mathcal{L}}{{\mathcal{L}}_{v}}A \right|}^{2}} \right.}} \notag\\ 
 &+&\left. {{\left| \slashed{\mathcal{L}}\slashed{\mathcal{L}}A \right|}^{2}}+{{\left| \slashed{\mathcal{L}}\slashed{\mathcal{L}}{{\mathcal{L}}_{v}}A \right|}^{2}} \right)d{{\sigma }^{2}}d\bar{v}~.  
\end{eqnarray}
By using the same procedure as in the far away from horizon case for computing the estimate of \eqref{Aint}, one can have
\begin{eqnarray}
    \int_{\bar{v}=v-1}^{\bar{v}=v}{\int_{{{S}^{2}}}{{{\left| A \right|}^{2}}}d{{\sigma }^{2}}d\bar{v}}\lesssim {{\left( \frac{w}{{{v}_{+}}} \right)}^{2}}{{E}_{1}}~,
\end{eqnarray}
which leads to
\begin{eqnarray}\label{Aes}
    \left| A \right|\lesssim {{\left\{ {{\left( \frac{w}{{{v}_{+}}} \right)}^{2}}{{E}_{1}}+{{E}_{1}} \right\}}^{1/2}}~.
\end{eqnarray}
Finally, combining the result in \eqref{nearphi}, \eqref{nabnearphi}, and \eqref{Aes} to \eqref{Depi}, we obtain \eqref{Dpi}.

This complete the proof of theorem \ref{teorema2}.

\section{Acknowledgments}

The work in this paper is supported by  Riset KK ITB 2023 and Riset Fundamental Kemendikbudristek-ITB 2023.

\end{document}